\documentclass[openacc]{rsproca_new}

\usepackage{amsfonts}
\usepackage{graphicx}
\usepackage{amsfonts}
\usepackage{graphicx}
\usepackage[fixamsmath,disallowspaces]{mathtools}
\usepackage{amsmath, amssymb, bm}
\usepackage{fixmath}

\newtheorem{proposition}{\bf Proposition}[section]

\newtheorem{lemma}{\bf Lemma}[section]
\newtheorem{remark}{\bf Remark}[section]

\begin{document}

\citearticle{M{\"u}ller A. 2023 Hamel's equations and geometric mechanics of constrained and floating multibody and space systems}{20220732}{479, 2023}

\title{
Hamel's Equations and Geometric Mechanics of Constrained and Floating Multibody and Space Systems
}

\author{%%%% Author details
Andreas M\"uller $^{1}$}

%%%%%%%%% Insert author address here
\address{$^{1}$Johannes Kepler University, Linz, Austria}

%%%% Subject entries to be placed here %%%%
\subject{Geometric mechanics, Lie groups, computational mechanics, robotics}

%%%% Keyword entries to be placed here %%%%

%%%% Insert corresponding author and its email address}
\corres{Andreas M\"uller\\
\email{a.mueller@jku.at}}

%%%% Abstract text to be placed here %%%%%%%%%%%%
\begin{abstract}
Modern geometric approaches to analytical mechanics rest on a bundle
structure of the configuration space. The connection on this bundle allows
for an intrinsic splitting of the reduced Euler-Lagrange equations. Hamel's
equations, on the other hand, provide a universal approach to non-holonomic
mechanics in local coordinates. The link between Hamel's formulation and
geometric approaches in local coordinates has not been discussed
sufficiently. \newline
The reduced Euler-Lagrange equations as well as the curvature of the
connection, are derived with Hamel's original formalism. Intrinsic splitting
into Euler-Lagrange and Euler-Poincar\'{e} equations, and inertial
decoupling is achieved by means of the locked velocity. Various aspects of
this method are discussed.
\end{abstract}

%%%%%%%%%%%%%%% End of first page %%%%%%%%%%%%%%%%%%%%%
\begin{fmtext}
\section{Introduction}
Many dynamical systems and controlled multibody systems possess symmetry
invariants, and can be modeled on a principle bundle. The bundle formulation
(for Lagrangian systems) was developed in \cite%
{MontgomeryMarsdenRatiu1984,Montgomery1986} as key concept in geometric
mechanics, where the configuration space is regarded as principle bundle $%
Q=G\times Q/G$ with symmetry group $G$. Central is the notion of a
connection as it allows encoding specific symmetries of the system \cite%
{BlochBook2003,KoonMarsden1997,Sreenath1992,ShammasChosetRizzi2007}. The 
\emph{(natural) mechanical connection}, deduced from the system momentum,
was introduced in \cite{Guichardet1984,Marsden1992,MarsdenScheurle1993}. 
Given a $G$-invariant Lagrangian $\ell $, in bundle coordinates $(\Omega
^{\alpha },\dot{r}^{I})$ of a (left) trivialization, where $\Omega ^{\alpha }
$ is the \emph{locked velocity}, the dynamics of unconstrained floating-base
systems 
%TCIMACRO{\TeXButton{black}{\color{black}}}%
%BeginExpansion
\color{black}%
%EndExpansion
is governed by the Lagrange-Poincar\'{e} equations \cite{BlochBook2003}%
\begin{align}
\frac{d}{dt}\frac{\partial \ell }{\partial \Omega ^{\alpha }}& =\frac{%
\partial \ell }{\partial \Omega ^{\beta }}\left( -\mathcal{E}_{I\alpha
}^{\beta }\dot{r}^{I}+c_{\alpha \lambda }^{\beta }\Omega ^{\lambda }\right) 
\label{LP1} \\
\frac{d}{dt}\frac{\partial \ell }{\partial \dot{r}^{I}}-\frac{\partial \ell 
}{\partial r^{I}}& =\frac{\partial \ell }{\partial \Omega ^{\beta }}\left( -%
\color[rgb]{0,0,0}\mathcal{B}_{IJ}^{\beta }\color{black}\dot{r}^{J}+%
\color[rgb]{0,0,0}\mathcal{E}_{I\alpha }^{\beta }\color{black}\Omega
^{\alpha }\right)   \label{LP2}
\end{align}%
with connection coefficients $\mathcal{A}_{I}^{\alpha }$ and its curvature $%
\mathcal{B}_{IJ}^{\alpha }$.
\end{fmtext}

\maketitle

%TCIMACRO{\TeXButton{red}{\color[rgb]{0,0,0}}}%
%BeginExpansion
\color[rgb]{0,0,0}%
%EndExpansion
The coefficients $\mathcal{E}_{\beta I}^{\alpha }=c_{\beta \lambda }^{\alpha
}\mathcal{A}_{I}^{\lambda }$, and the curvature coefficients%
\begin{equation}
\mathcal{B}_{IJ}^{\alpha }=\frac{\partial \mathcal{A}_{I}^{\alpha }}{%
\partial r^{J}}-\frac{\partial \mathcal{A}_{J}^{\alpha }}{\partial r^{I}}\pm
c_{\beta \lambda }^{\alpha }\mathcal{A}_{I}^{\beta }\mathcal{A}_{J}^{\lambda
}  \label{curvature}
\end{equation}%
%TCIMACRO{\TeXButton{black}{\color{black}}}%
%BeginExpansion
\color{black}%
%EndExpansion
are determined by the structure constants $c_{\beta \lambda }^{\alpha }$ of
the Lie algebra $\mathfrak{g}$ of $G$. 
%TCIMACRO{\TeXButton{red}{\color[rgb]{0,0,0}}}%
%BeginExpansion
\color[rgb]{0,0,0}%
%EndExpansion
The symmetry group encodes invariances of the Lagrangian. Most prominent
examples are floating-base multibody systems\footnote{%
%TCIMACRO{\TeXButton{red}{\color[rgb]{0,0,0}}}%
%BeginExpansion
\color[rgb]{0,0,0}%
%EndExpansion
Throughout the paper, the term floating-base systems covers a large class of
mechanical (control) systems characterized by a base body free to move in
space to which further bodies (links) are geometrically connected, e.g. by
joints.}, where the motion $g$ of a base body evolves in a subgroup $G$ of
the group of Euclidean motions $SE\left( 3\right) $, with body-velocity $%
\hat{%
%TCIMACRO{\TeXButton{xi}{\bm{\xi}}}%
%BeginExpansion
\bm{\xi}%
%EndExpansion
}=g^{-1}\dot{g}\in \mathfrak{g}$, and $Q/G$ is the 'joint space' (which will
be identified with ${\mathbb{V}}^{n}$) with joint coordinates $r^{I}$. The
kinetic energy is then $G$-invariant, and the locked velocity is defined in
terms of the mechanical connection as $\Omega ^{\alpha }:=u^{\alpha }$ (see
Sec. \ref{secUnConBundle}), with%
\begin{equation}
u^{\alpha }=\xi ^{\alpha }+\mathcal{A}_{I}^{\alpha }(r^{I})\dot{r}^{I}.
\label{w}
\end{equation}%
The curvature vanishes if and only if the constraints defined by momentum
conservation, and thus the Pfaffian system (\ref{w}), is integrable.
Similarly for constrained systems, a \emph{kinematic connection} was
introduced, and the scleronomic constraints are expressed as $u^{\alpha }=0$%
. Mechanical systems whose spatial motion is constrained are typical
examples, where the constraints restrict the motion of a base body in $G$.
Now the symmetry group accounts for the invariance of the constraints, and
the curvature vanishes if and only if the constraints are integrable. The
dynamics equations for constrained systems are obviously obtained from (\ref%
{LP2}) by setting $\Omega ^{\alpha }=0$. Thus, in both cases, a kinematic
relation of the form (\ref{w}) applies, and the curvature appears in the
equations of motion, and is a central object in motion planning and control.

Geometric mechanics provides an intrinsic and coordinate-free framework for
modeling, analysis, and control of finite-dimensional (discrete) as well as
infinite-dimensional (continua) systems. Local coordinate formulations, as
the one above, are used for computations, where it may be necessary to
switch between different local coordinate charts. In this context, the fact
that (\ref{LP1},\ref{LP2}) is a specific form of the Hamel equations \cite%
{Hamel-MathAnal1924,HamelBook} for finite-dimensional systems in local
coordinates deserves recognition, which is the topic of this paper. In
geometric mechanics, Hamel's equations are now introduced in an elegant
modern form, e.g. \cite{Zenkov2016}, on the expense that the relation to the
original Hamel equations is lost, however. Moreover, how the geometric
framework and the Hamel formalism are related, and how connection and
curvature of the configuration space (bundle) are related to the Hamel
coefficients, is not discussed in the literature. Without making explicit
reference to the original formulation, (\ref{LP1},\ref{LP2}) are referred to
as Hamel equations \cite{Marsden1992}. 
%TCIMACRO{\TeXButton{red}{\color[rgb]{0,0,0}}}%
%BeginExpansion
\color[rgb]{0,0,0}%
%EndExpansion
Also in robotics and multibody system dynamics%
%TCIMACRO{\TeXButton{black}{\color{black}}}%
%BeginExpansion
\color{black}%
%EndExpansion
, the Euler-Poincar\'{e} equations (\ref{LP1}) on $SE\left( 3\right) $ are
often interchangeably referred to as Hamel equations or Euler-Poincar\'{e}
equations, e.g. \cite{Mishra2020}.

It is shown in this paper that the coefficients in (\ref{LP1},\ref{LP2}) are
naturally derived as the Hamel coefficients. This provides a link between
the original Hamel formalism and the bundle formulation. The explicit
derivation also admits consolidating the different 
%TCIMACRO{\TeXButton{red}{\color[rgb]{0,0,0}}}%
%BeginExpansion
\color[rgb]{0,0,0}%
%EndExpansion
coordinate expressions of local curvature 
%TCIMACRO{\TeXButton{black}{\color{black}}}%
%BeginExpansion
\color{black}%
%EndExpansion
found in the literature, which is crucial for applying the above equations
(see Rem. \ref{remSign}). There are also various aspects that need to be
taken into account when using the geometric formalism. One is the concept of
locked velocity \cite{BlochBook2003}, and the related concept of average
velocity \cite{OrinGoswami2008,OrinGoswamiLee2013} that proved to be
powerful tools for control of floating-base systems. It is discussed that
there is no frame which can be associated with this velocity whose motion is
a function of $g\in G$ and $\mathbf{r}\in {\mathbb{V}}^{n}$, which has
consequences for control of floating-base systems. Another aspect discussed
in this paper is that the mechanical connection on $Q$ induced by the locked
velocity intrinsically splits the reduced Euler-Lagrange equations in
horizontal and vertical. As an important consequence, the equations (\ref%
{LP1},\ref{LP2}) are inertially decoupled, which is relevant for
computational multibody dynamics. 
%TCIMACRO{\TeXButton{red}{\color[rgb]{0,0,0}}}%
%BeginExpansion
\color[rgb]{0,0,0}%
%EndExpansion
Throughout the paper, all constraints are assumed to be linear (i.e.
catastatic) and scleronomic.%
%TCIMACRO{\TeXButton{black}{\color{black}}}%
%BeginExpansion
\color{black}%
%EndExpansion

The paper is organized as follows. In Sec. \ref{secBHEqu}, the classical
Hamel equations are summarized for unconstrained and constrained systems. 
%TCIMACRO{\TeXButton{red}{\color[rgb]{0,0,0}}}%
%BeginExpansion
\color[rgb]{0,0,0}%
%EndExpansion
Hamel's formulation is related to the formulation on a (locally) trivial
bundle in Sec. \ref{secKinConstrBundle}. The obtained relations for the
Hamel coefficients in local coordinates are the basis for the derivations in
all subsequent sections. In Sec. \ref{secSymConstr}, kinematically
constrained systems with symmetries are treated, where the configuration
space $Q$ is a principal bundle, and the curvature coefficients $\mathcal{B}%
_{IJ}^{\alpha }$ are obtained as the Hamel coefficients in local bundle
coordinates. This principle bundle approach is adopted in Sec. \ref%
{secUnConBundle} for unconstrained systems with $G$-invariant Lagrangian,
where the coefficients $\mathcal{E}_{\beta I}^{\alpha }$ are obtained
immediately as the Hamel coefficients. A rolling sphere and a floating
satellite are used to demonstrate application of the equations. Finally,
unconstrained floating-base systems with conserved momentum are considered
in Sec. \ref{secUnConDymConserved} with a note on geometric phase and
pseudo-holonomic motions. Numerical simulation results for a satellite and a
space robot are reported in Sec. \ref{secUnConDymConserved} and in the
supplement \cite{Supplement}. 
%TCIMACRO{\TeXButton{black}{\color{black}}}%
%BeginExpansion
\color{black}%
%EndExpansion
For background material, an excellent introduction to geometric mechanics
can be found in the text books \cite{BlochBook2003,Holm2011Book1} and the
overview articles \cite%
{MarsdenScheurle1993,BlochKrishnaprasadMarsdenMurray1996,KoonMarsden1997}. 
%TCIMACRO{\TeXButton{red}{\color[rgb]{0,0,0}}}%
%BeginExpansion
\color[rgb]{0,0,0}%
%EndExpansion
Relevant concepts from differential geometry and on bundles can be found in 
\cite{BaezMuniain1994,Frankel2004,BlochBook2003}. For all derivations,
Hamel's original approach is the point of departure. This relies on
coordinates on $Q$, and it is necessary to introduce local coordinates also
on the symmetry group (yet the final formulation (\ref{LP1},\ref{LP2}) only
needs local coordinates on $Q/G$ with parameter space ${\mathbb{V}}^{n}$).
It applies to finite-dimensional systems for which always exist local
coordinates on $Q$, and canonical coordinates on $G$ such as multibody
systems. This is in contrast to modern geometric mechanics, where globally
valid equations are derived coordinate free, and local coordinates are
introduced when needed. However, the coordinate formulations allow to relate
the geometric approach to Hamel's formulation. It is assumed throughout the
paper that $G$ is a finite-dimensional Lie group possessing local
coordinates. Let the Lie algebra $\mathfrak{g}$ be isomorphic to the vector
space ${\mathbb{R}}^{n}$. Then $\hat{\mathbf{x}}\in \mathfrak{g}$ denotes
the Lie algebra element corresponding to the vector $\mathbf{x}\in {\mathbb{R%
}}^{n}\cong \mathfrak{g}$. Ricci's summation convention is used, e.g. $%
B_{I}^{a}u^{I}=\sum_{I}B_{I}^{a}u^{I}$ implies summation over index $I$. The
notation is summarized in appendix A.%
%TCIMACRO{\TeXButton{black}{\color{black}}}%
%BeginExpansion
\color{black}%
%EndExpansion

\section{The Hamel Equations%
%TCIMACRO{\TeXButton{secBHEqu}{\label{secBHEqu}}}%
%BeginExpansion
\label{secBHEqu}%
%EndExpansion
}

\subsection{Boltzmann-Hamel Equations in Quasi-Velocities}

%TCIMACRO{\TeXButton{red}{\color[rgb]{0,0,0}}}%
%BeginExpansion
\color[rgb]{0,0,0}%
%EndExpansion
The configuration is described in terms of $n$ generalized coordinates $%
q^{a},a=1,\ldots ,n$, with parameter manifold $Q={\mathbb{V}}^{n}$, which
serves as configuration space. In case of a multibody system with $n_{%
\mathrm{t}}$ translation and $n_{\mathrm{r}}$ revolute joints, for instance, 
${\mathbb{V}}^{n}=T^{n_{\mathrm{r}}}\times {\mathbb{R}}^{n_{\mathrm{t}}}$,
with $n_{\mathrm{r}}$-torus $T^{n_{\mathrm{r}}}$ and $n=n_{\mathrm{r}}+n_{%
\mathrm{t}}$. For such systems, the generalized coordinates may only be
locally valid. 
%TCIMACRO{\TeXButton{black}{\color{black}}}%
%BeginExpansion
\color{black}%
%EndExpansion
Quasi-velocity coordinates $u^{a},a=1,\ldots ,n$, are introduced that are
related to the generalized speeds $\dot{q}^{a}$ by%
\begin{equation}
u^{a}=A_{b}^{a}\dot{q}^{b},\ \ \ \dot{q}^{a}=B_{b}^{a}u^{b}  \label{qdu}
\end{equation}%
%TCIMACRO{\TeXButton{red}{\color[rgb]{0,0,0}}}%
%BeginExpansion
\color[rgb]{0,0,0}%
%EndExpansion
where $A_{b}^{a}$ and $B_{b}^{a}$ are smooth functions of $\mathbf{q}$. 
%TCIMACRO{\TeXButton{black}{\color{black}}}%
%BeginExpansion
\color{black}%
%EndExpansion
With vectors of generalized coordinates $\mathbf{q}\in {\mathbb{V}}^{n}$ and
quasi-velocities $\mathbf{u}\in {\mathbb{R}}^{n}$, these relations are
written in matrix form as%
\begin{equation}
\mathbf{u}=\mathbf{A}\left( \mathbf{q}\right) \dot{\mathbf{q}},\ \ \ \dot{%
\mathbf{q}}=\mathbf{B\left( \mathbf{q}\right) u}.  \label{u2}
\end{equation}%
%TCIMACRO{\TeXButton{red}{\color[rgb]{0,0,0}}}%
%BeginExpansion
\color[rgb]{0,0,0}%
%EndExpansion
It is assumed that $q^{a}$ are valid local coordinates so that $\mathbf{A}$
is regular, and $\mathbf{B}=\mathbf{A}^{-1}$. 
%TCIMACRO{\TeXButton{black}{\color{black}}}%
%BeginExpansion
\color{black}%
%EndExpansion
The relation of $u^{a}$ to the corresponding quasi-coordinates $\pi
^{a},a=1,\ldots ,n$ is described by the differential forms $d\pi
^{a}=A_{b}^{a}dq^{b}$, and the inverse relation by $dq^{a}=B_{b}^{a}d\pi
^{b} $. The quasi-coordinates are generalized coordinates if and only if the
differential forms are exact, otherwise $u^{a}$ are called non-holonomic
velocities \cite{Whittaker1988,PapastavridisBook2002}%
%TCIMACRO{\TeXButton{red}{\color[rgb]{0,0,0}}}%
%BeginExpansion
\color[rgb]{0,0,0}%
%EndExpansion
, following \cite[p. 473]{HamelBook},\cite[p. 218]{Sommerfeld1948}.%
%TCIMACRO{\TeXButton{black}{\color{black}}}%
%BeginExpansion
\color{black}%
%EndExpansion

Denote with $L\left( q^{a},u^{a}\right) $ a Lagrangian in terms of
quasi-velocity coordinates $u^{a}$, which for mechanical systems is defined
as kinetic energy $T\left( q^{a},u^{a}\right) $ minus potential energy $%
U\left( q^{a}\right) $. Then the variational form of the Hamel equations 
\cite{Hamel-MathAnal1904,Hamel-MathAnal1924,HamelBook,PapastavridisBook2002}
is%
\begin{equation}
\left( \frac{d}{dt}\frac{\partial L}{\partial u^{a}}-\frac{\partial L}{%
\partial \pi ^{a}}-Q_{a}\right) \delta \pi ^{a}+\frac{\partial L}{\partial
u^{a}}\left( \frac{d\delta \pi ^{a}-\delta d\pi ^{a}}{dt}\right) =0
\label{BHLagrange}
\end{equation}%
where $Q_{a}$ are generalized forces dual to $u^{a}$. The explicit form of
the \emph{Boltzmann-Hamel equations }%
\begin{equation}
\frac{d}{dt}\frac{\partial L}{\partial u^{a}}-\frac{\partial L}{\partial
q^{b}}B_{a}^{b}+\frac{\partial L}{\partial u^{b}}\gamma _{ac}^{b}u^{c}=Q_{a}
\label{BHEqu2}
\end{equation}%
is obtained after inserting the \emph{transitory relation} 
\begin{equation}
\frac{d\delta \pi ^{a}-\delta d\pi ^{a}}{dt}=\gamma _{cb}^{a}u^{b}\delta \pi
^{c}  \label{transpos}
\end{equation}%
in which $\gamma _{cb}^{a}$ are the \emph{Hamel coefficients} defined as%
\begin{equation}
\gamma _{ab}^{c}:=\left( \frac{\partial A_{r}^{c}}{\partial q^{s}}-\frac{%
\partial A_{s}^{c}}{\partial q^{r}}\right) B_{a}^{r}B_{b}^{s}.  \label{gamma}
\end{equation}%
The Hamel coefficients vanish identically if and only if (\ref{qdu}) are
integrable, i.e. if $u^{a}$ are holonomic velocities. The equations (\ref%
{BHEqu2}) along with the kinematic equations (\ref{qdu}) govern the dynamics
in the non-holonomic tangent bundle 
%TCIMACRO{\TeXButton{red}{\color[rgb]{0,0,0}}}%
%BeginExpansion
\color[rgb]{0,0,0}%
%EndExpansion
(i.e. tangent space defined by non-holonomic constraints). 
%TCIMACRO{\TeXButton{black}{\color{black}}}%
%BeginExpansion
\color{black}%
%EndExpansion
They are referred to as the Boltzmann-Hamel equations (e.g. in \cite%
{Maisser1997,MaruskinBloch2007}) as they where (in a very similar form)
presented by Boltzmann in \cite{Boltzmann1902,Boltzmann1904} and by Hamel in 
\cite{Hamel1904,Hamel-MathAnal1924,HamelBook}. It was Hamel, however, who
generalized them to systems an a Lie group \cite%
{Hamel1904,Hamel-MathAnal1904}.

The motion of many mechanical systems can be described on a $n$-dimensional
Lie-group $G$, so that quasi-velocities belong the corresponding Lie algebra 
$\mathfrak{g}$. If $g\left( t\right) \in G$, such quasi-velocities can be
introduced as left- or right-invariant vector fields, $\hat{\mathbf{u}}%
=g^{-1}\dot{g}$ or $\hat{\mathbf{u}}=\dot{g}g^{-1}$, respectively. Let $%
q^{a} $ be canonical coordinates on $G$, and let $\mathbf{u}\in {\mathbb{R}}%
^{n}$ be the vector representation of $\hat{\mathbf{u}}\in \mathfrak{g}$,
then (\ref{u2}) is a map from ${\mathbb{R}}^{n}$ to ${\mathbb{R}}^{n}\cong 
\mathfrak{g}$. When using the left-invariant definition of quasi-velocities,
the Hamel coefficients (\ref{gamma}) are identical to the structure
constants $c_{ab}^{c}$ of the Lie group. This was already shown by Hamel 
\cite[p. 428]{Hamel-MathAnal1904} using the transitory relations $d\delta
\pi ^{a}-\delta d\pi ^{a}=\gamma _{cb}^{a}d\pi ^{c}\delta \pi ^{c}$. Using
right-invariant quasi-velocities leads to a change of sign: $\gamma
_{ab}^{c}=-c_{ab}^{c}$. Typical example for such quasi-velocities are the
angular velocity (where $G=SO\left( 3\right) $) or rigid body twists (where $%
G=SE\left( 3\right) $). Then left-invariance implies body-fixed
representation of angular velocity or twists, and right-invariance implies
spatial representation \cite{Murray,MUBOScrews1}. The explicit derivation of
the Hamel-coefficients for $SO\left( 3\right) $ and $SE\left( 3\right) $
using the definition of Hamel coefficients can be found in \cite%
{MuellerJNLS2021}. Clearly, $\pi ^{a}$ are holonomic coordinates if and only
if $G$ is Abelian. 
%TCIMACRO{\TeXButton{red}{\color[rgb]{0,0,0}}}%
%BeginExpansion
\color[rgb]{0,0,0}%
%EndExpansion
It must be emphasized that canonical coordinates on $G$ are only locally
valid in general. This applies in particular to $SO\left( 3\right) $, and
thus to $SE\left( 3\right) $, since it is not simply connected, which leads
to the well-known parameterization singularity of rotations.%
%TCIMACRO{\TeXButton{black}{\color{black}}}%
%BeginExpansion
\color{black}%
%EndExpansion

In summary, the Hamel equations on a Lie group are the forced Euler-Poincar%
\'{e} equations for Lagrangian $L\left( q^{a},u^{a}\right) $ 
\begin{equation}
\frac{d}{dt}\frac{\partial L}{\partial u^{a}}\pm c_{ac}^{b}\frac{\partial L}{%
\partial u^{b}}u^{c}=\frac{\partial L}{\partial q^{b}}B_{a}^{b}+Q_{a}
\label{Poincare1}
\end{equation}%
where the positive sign applies to left-invariant, and the negative sign to
right-invariant quasi-velocities. These are the reduced Euler-Lagrange
equations for systems whose Lagrangian is (left or right) invariant under
action of a symmetry group $G$ \cite{MarsdenScheurle1993,MarsdenBook1995}.

\subsection{Hamel Equations for Constrained Systems in Quasi-Velocities}

The original Hamel equations for constrained systems where presented in \cite%
{Hamel-MathAnal1904}. The velocities $\dot{\mathbf{q}}\in \mathbb{R}^{n}$
are now subjected to $\bar{m}$ Pfaffian constraints, written as $u^{\alpha
}=0,\alpha =1,\ldots ,\bar{m}$, with%
\begin{equation}
u^{\alpha }:=A_{a}^{\alpha }\left( q^{a}\right) \dot{q}^{a}.
\label{PfaffianConstr}
\end{equation}%
%TCIMACRO{\TeXButton{red}{\color[rgb]{0,0,0}}}%
%BeginExpansion
\color[rgb]{0,0,0}%
%EndExpansion
It is assumed that the system of Pfaffian constraints is regular, i.e. the $%
\bar{m}$ constraints are independent. 
%TCIMACRO{\TeXButton{black}{\color{black}}}%
%BeginExpansion
\color{black}%
%EndExpansion
Then, $\bar{\delta}:=n-\bar{m}$ independent quasi-velocity coordinates are
introduced as%
\begin{equation}
u^{I}=A_{a}^{I}\left( q^{a}\right) \dot{q}^{a},\ \ I=\bar{m}+1,\ldots ,n
\label{PfaffForm}
\end{equation}%
%TCIMACRO{\TeXButton{red}{\color[rgb]{0,0,0}}}%
%BeginExpansion
\color[rgb]{0,0,0}%
%EndExpansion
where $\bar{\delta}$ is the differential DOF of the system (also called
instantaneous DOF) \cite{IFToMM_Terminology_MMT2003,CISM_Book_2019}. 
%TCIMACRO{\TeXButton{black}{\color{black}}}%
%BeginExpansion
\color{black}%
%EndExpansion
The overbar of $\bar{\delta}$ and $\bar{m}$ indicate that the constraints
are generally non-holonomic. If they are integrable, there are $m=\bar{m}$
geometric constraints, and $\delta =n-m=\bar{\delta}$ is the finite DOF. The 
$n-\bar{m}$ independent coordinates are indexed with capital lattin letters $%
I,J,K$. The Pfaffian system (\ref{PfaffianConstr}) and the solution (\ref%
{PfaffForm}) are summarized as $u^{a}=A_{b}^{a}\left( q^{a}\right) \dot{q}%
^{b}$, as in (\ref{qdu}) with index set $\{a\}=\{\alpha ,I\}$. 
%TCIMACRO{\TeXButton{red}{\color[rgb]{0,0,0}}}%
%BeginExpansion
\color[rgb]{0,0,0}%
%EndExpansion
The independent coordinates are only locally valid. Moreover, the
configuration space of a constrained system is in general not globally a
manifold but possesses singularities. This strictly limits the global
validity all coordinate formulations.%
%TCIMACRO{\TeXButton{black}{\color{black}}}%
%BeginExpansion
\color{black}%
%EndExpansion

Denote with $\mathcal{L}(q^{a},u^{b}):=L(q^{a},B_{b}^{a}u^{b})$ the
Lagrangian in which $\dot{q}^{a}$ is replaced by $u^{a}$, by means of (\ref%
{qdu}). The Hamel equations for the constrained system, in terms of the
independent velocities $u^{I}$, are then obtained from (\ref{BHEqu2}) as%
\begin{equation}
\frac{d}{dt}\frac{\partial \mathcal{L}}{\partial u^{I}}-\frac{\partial 
\mathcal{L}}{\partial q^{a}}B_{I}^{a}+\frac{\partial \mathcal{L}}{\partial
u^{a}}\gamma _{IJ}^{a}u^{J}=Q_{I},\ \ I=\bar{m}+1,\ldots ,n  \label{BHconstr}
\end{equation}%
where $u^{\alpha },\alpha =1,\ldots \bar{m}$ are set to zero after taking
the derivatives. The generalized forces are $Q_{I}=B_{I}^{a}Q_{a}$. The
Hamel coefficients 
%TCIMACRO{\TeXButton{red}{\color[rgb]{0,0,0}}}%
%BeginExpansion
\color[rgb]{0,0,0}%
%EndExpansion
in (\ref{BHconstr}) are obtained by restricting (\ref{gamma}) to indices $%
I,J $ as%
%TCIMACRO{\TeXButton{black}{\color{black}}}%
%BeginExpansion
\color{black}%
%EndExpansion
\begin{equation}
\gamma _{IJ}^{a}=\left( \frac{\partial A_{b}^{a}}{\partial q^{c}}-\frac{%
\partial A_{c}^{a}}{\partial q^{b}}\right) B_{I}^{b}B_{J}^{c}.
\label{gammaNH}
\end{equation}%
The $n-\bar{m}$ equations (\ref{BHconstr}) complemented with the $n$
kinematic equations%
\begin{equation}
\dot{q}^{a}=B_{I}^{a}\left( q^{a}\right) u^{I},\ \ I=\bar{m}+1,\ldots ,n
\label{qadot}
\end{equation}%
govern the system dynamics in terms of state variables $(q^{a},u^{I})$. Not
all of the $q^{a}$ may be independent if the constraints are not completely
non-holonomic. 
%TCIMACRO{\TeXButton{red}{\color[rgb]{0,0,0}}}%
%BeginExpansion
\color[rgb]{0,0,0}%
%EndExpansion
Relation (\ref{qadot}) is obtained from (\ref{u2}) assuming $u^{I}$ are
locally valid coordinates on the tangent space and constraints are regular.%
%TCIMACRO{\TeXButton{black}{\color{black}}}%
%BeginExpansion
\color{black}%
%EndExpansion

Quasi-velocities $u^{I}$ are integrable if and only if $\gamma
_{IJ}^{K}\equiv 0$. The constraints (\ref{PfaffianConstr}), and thus the
co-distribution $D^{\ast }\subset T_{\mathbf{q}}^{\ast }Q$ with $D_{\mathbf{q%
}}^{\ast }:=\mathrm{span}~(\mathbf{A}\left( \mathbf{q}\right) )$ defined by
the constraints, are integrable (in Pfaffian sense) if and only if $\gamma
_{IJ}^{\alpha }\equiv 0$. The \emph{constraint distribution} $D$ on $Q$,
defined as $D_{\mathbf{q}}:=\ker \mathbf{A}\left( \mathbf{q}\right) \subset
T_{\mathbf{q}}Q$, is thus integrable (in Cartan sense) if and only if the
Hamel coefficients vanish and the constraints are regular 
%TCIMACRO{\TeXButton{red}{\color[rgb]{0,0,0}}}%
%BeginExpansion
\color[rgb]{0,0,0}%
%EndExpansion
($\mathrm{\mathrm{rank}}\,\mathbf{A}$ is constant)%
%TCIMACRO{\TeXButton{black}{\color{black}}}%
%BeginExpansion
\color{black}%
%EndExpansion
. This may not apply to non-regular constraints.

\section{Kinematically Constrained Systems on a Trivial Bundle%
%TCIMACRO{\TeXButton{secKinConstrBundle}{\label{secKinConstrBundle}}}%
%BeginExpansion
\label{secKinConstrBundle}%
%EndExpansion
}

Many kinematic control problems can be formulated on a trivial bundle.
Trivial because there is a global splitting into independent and dependent
velocities. 
%TCIMACRO{\TeXButton{red}{\color[rgb]{0,0,0}}}%
%BeginExpansion
\color[rgb]{0,0,0}%
%EndExpansion
The independent velocities serve as control inputs. 
%TCIMACRO{\TeXButton{black}{\color{black}}}%
%BeginExpansion
\color{black}%
%EndExpansion
Moreover, many control ~system are in Chaplygin form, i.e. the kinematic
relations only depend on the independent coordinates.

\subsection{Constrained Hamel Equations in Terms of Holonomic Velocities}

Consider (mechanical) systems described by coordinates $q^{a}$ and their
time derivatives $\dot{q}^{a}$, rather then non-holonomic velocities $u^{I}$%
, subjected to 
%TCIMACRO{\TeXButton{red}{\color[rgb]{0,0,0}}}%
%BeginExpansion
\color[rgb]{0,0,0}%
%EndExpansion
scleronomic 
%TCIMACRO{\TeXButton{black}{\color{black}}}%
%BeginExpansion
\color{black}%
%EndExpansion
non-holonomic Pfaffian constraints (\ref{PfaffianConstr}). A set of
(locally) independent velocity coordinates can be selected.

A particular choice of independent velocities is to use time derivatives of $%
\bar{\delta}:=n-\bar{m}$ coordinates. To this end, the coordinates are
partitioned as $\mathbf{q}=(s^{\alpha },r^{I})\in {\mathbb{V}}^{\bar{m}%
}\times {\mathbb{V}}^{\bar{\delta}}=:Q$, where $\dot{s}^{\alpha },1,\ldots ,%
\bar{m}$ are the dependent, and the remaining $\dot{r}^{I},I=\bar{m}%
+1,\ldots ,n$ are independent velocity coordinates, i.e. $u^{I}:=\dot{r}^{I}$%
. 
%TCIMACRO{\TeXButton{red}{\color[rgb]{0,0,0}}}%
%BeginExpansion
\color[rgb]{0,0,0}%
%EndExpansion
This presumes that $\dot{r}^{I}$ are valid local coordinates on $Q$. 
%TCIMACRO{\TeXButton{black}{\color{black}}}%
%BeginExpansion
\color{black}%
%EndExpansion
Notice that for non-holonomic constraints, this dependency does not hold
true for the coordinates $s^{\alpha },r^{I}$, and $Q$ serves as $n$%
-dimensional configuration space. The constraints (\ref{PfaffianConstr}) are
then written as%
\begin{align}
u^{\alpha } & := A_{\beta }^{\alpha }\left( q^{a}\right) \dot{s}^{\beta
}+A_{I}^{\alpha }\left( q^{a}\right) \dot{r}^{I},\alpha =1,\ldots ,\bar{m}
\label{PfaffianConstr2} \\
\mathbf{u} &=\mathbf{A}\left( \mathbf{q}\right) \dot{\mathbf{q}}=\left( 
\begin{array}{cc}
\mathbf{A}_{1} & \mathbf{A}_{2} \\ 
\mathbf{0} & \mathbf{I}%
\end{array}%
\right) \left( 
\begin{array}{c}
\dot{\mathbf{s}} \\ 
\dot{\mathbf{r}}%
\end{array}%
\right) ,\ \ \mathrm{with\ }\mathbf{A}_{1}=\left( A_{\beta }^{\alpha
}\right) ,\mathbf{A}_{2}=\left( A_{I}^{\alpha }\right) .  \label{sys1}
\end{align}%
where (\ref{sys1}) resembles the matrix form (\ref{u2}). The inverse
relation of (\ref{PfaffianConstr2}) and (\ref{sys1}) are, respectively,%
\begin{align}
\dot{s}^{\alpha } &=B_{\beta }^{\alpha }\left( q^{a}\right) u^{\beta
}+B_{I}^{\alpha }\left( q^{a}\right) \dot{r}^{I}  \label{salpha} \\
\left( 
\begin{array}{c}
\dot{\mathbf{s}} \\ 
\dot{\mathbf{r}}%
\end{array}%
\right) &=\mathbf{B\left( \mathbf{q}\right) u},\ \ \mathrm{with\ }\mathbf{B}%
=\left( 
\begin{array}{cc}
\mathbf{B}_{1} & \mathbf{B}_{2} \\ 
\mathbf{0} & \mathbf{I}%
\end{array}%
\right) ,\ \ \mathrm{with\ }%
%TCIMACRO{\TeXButton{red}{\color[rgb]{0,0,0}}}%
%BeginExpansion
\color[rgb]{0,0,0}%
%EndExpansion
\mathbf{B}_{1}=\mathbf{A}_{1}^{-1},\mathbf{B}_{2}=-\mathbf{A}_{1}^{-1}%
\mathbf{A}_{2}.  \label{sys2}
\end{align}%
With $\mathbf{A}$ in (\ref{sys1}) and $\mathbf{B}$ in (\ref{sys2}), the
expression (\ref{gammaNH}) gives rise to the Hamel coefficients%
\begin{align}
\gamma _{\beta \lambda }^{\alpha }& =\left( \frac{\partial A_{\mu }^{\alpha }%
}{\partial s^{\nu }}-\frac{\partial A_{\nu }^{\alpha }}{\partial s^{\mu }}%
\right) B_{\beta }^{\mu }B_{\lambda }^{\nu }  \notag \\
\gamma _{IJ}^{\alpha }& =\left( \frac{\partial A_{r}^{\alpha }}{\partial
q^{s}}-\frac{\partial A_{s}^{\alpha }}{\partial q^{r}}\right)
B_{I}^{r}B_{J}^{s}  \label{gammaHolo} \\
& =\frac{\partial A_{I}^{\alpha }}{\partial r^{J}}-\frac{\partial
A_{J}^{\alpha }}{\partial r^{I}}+\left( \frac{\partial A_{\beta }^{\alpha }}{%
\partial s^{\mu }}-\frac{\partial A_{\mu }^{\alpha }}{\partial s^{\beta }}%
\right) B_{I}^{\beta }B_{J}^{\mu }+\left( \frac{\partial A_{I}^{\alpha }}{%
\partial s^{\beta }}-\frac{\partial A_{\beta }^{\alpha }}{\partial r^{I}}%
\right) B_{J}^{\beta }+\left( \frac{\partial A_{\beta }^{\alpha }}{\partial
r^{J}}-\frac{\partial A_{J}^{\alpha }}{\partial s^{\beta }}\right)
B_{I}^{\beta }  \notag \\
\gamma _{\beta J}^{\alpha }& =\left( \frac{\partial A_{r}^{\alpha }}{%
\partial q^{s}}-\frac{\partial A_{s}^{\alpha }}{\partial q^{r}}\right)
B_{\beta }^{r}B_{J}^{s}=\left( \frac{\partial A_{\nu }^{\alpha }}{\partial
s^{\mu }}-\frac{\partial A_{\mu }^{\alpha }}{\partial s^{\nu }}\right)
B_{\beta }^{\nu }B_{J}^{\mu }+\left( \frac{\partial A_{\mu }^{\alpha }}{%
\partial r^{J}}-\frac{\partial A_{J}^{\alpha }}{\partial s^{\beta }}\right)
B_{\beta }^{\mu },\ \ \alpha =1,\ldots ,\bar{m}.  \notag
\end{align}%
Since the velocities are integrable, it holds true that $\gamma
_{ab}^{K}\equiv 0$. The Hamel coefficients vanish if and only if the
constraints are holonomic. The expressions (\ref{gammaHolo}) will be central
throughout the paper as the individual terms in (\ref{gammaHolo}) allow
deriving the coordinate form of the reduced Euler-Lagrange equations for
systems with symmetry directly from the Hamel formulation.

The Lagrangian is written as $L(s^{\alpha },r^{I},\dot{s}^{\alpha },\dot{r}%
^{I})$ to indicate the coordinate partitioning. As in (\ref{BHconstr}),
denote with $\mathcal{L}(s^{\alpha },r^{I},u^{a}):=L(s^{\alpha
},r^{I},B_{\beta }^{\alpha }u^{\beta }+B_{I}^{\alpha }u^{I},u^{I})$ the
Lagrangian with $\dot{r}^{I}=u^{I}$ and $\dot{s}^{\alpha }$ replaced by (\ref%
{salpha}). The Hamel equations follow from (\ref{BHconstr}). Noting that $%
\frac{\partial \mathcal{L}}{\partial u^{\alpha }}=\frac{\partial L}{\partial 
\dot{s}^{\beta }}B_{\alpha }^{\beta }$, and (with slight abuse of notation)
identifying $u^{I}=\dot{r}^{I}$, yields the Hamel equations in independent
holonomic velocities $\dot{r}^{I}$ 
\begin{equation}
\frac{d}{dt}\frac{\partial \mathcal{L}}{\partial \dot{r}^{I}}-\frac{\partial 
\mathcal{L}}{\partial r^{I}}+\frac{\partial L}{\partial \dot{s}^{\alpha }}%
\gamma _{IJ}^{\alpha }\dot{r}^{J}=Q_{I},\ \ I=\bar{m}+1,\ldots ,n
\label{HE-holo}
\end{equation}%
where $u^{\alpha }$ are set to zero, i.e. the solution%
\begin{equation}
\dot{s}^{\alpha }=B_{I}^{\alpha }\left( q^{a}\right) \dot{r}^{I}
\label{solConstr}
\end{equation}%
of the constraints is imposed, after taking the derivatives. The dynamic
equations (\ref{HE-holo}) along with the kinematic equations (\ref{solConstr}%
) govern the dynamics of the non-holonomically constrained system in terms
of the state $\left( q^{a},\dot{q}^{a}\right) $ evolving on the
non-holonomic tangent bundle.

\subsection{Kinematic Constraints in Terms of a Bundle Connection}

The kinematic constraints $u^{\alpha }=0$ are now formulated with%
\begin{equation}
u^{\alpha }:=\dot{s}^{\alpha }+\mathcal{A}_{I}^{\alpha }\left( q^{a}\right) 
\dot{r}^{I},\alpha =1,\ldots ,\bar{m}  \label{PfaffianConstr3}
\end{equation}%
where $\left( \mathcal{A}_{I}^{\alpha }\right) :=\mathbf{A}_{1}^{-1}\mathbf{A%
}_{2}=-\mathbf{B}_{2}$, with $\mathbf{B}_{2}=\left( B_{I}^{\alpha }\right) $
in (\ref{salpha}). Combined with $u^{I}:=\dot{r}^{I}$, this is written in
matrix form as%
\begin{equation}
\mathbf{u}=\bar{\mathbf{A}}\left( \mathbf{q}\right) \dot{\mathbf{q}}=\left( 
\begin{array}{cc}
\mathbf{I} & \ \ -\mathbf{B}_{2} \\ 
\mathbf{0} & \mathbf{I}%
\end{array}%
\right) \left( 
\begin{array}{c}
\dot{\mathbf{s}} \\ 
\dot{\mathbf{r}}%
\end{array}%
\right) .  \label{Abar}
\end{equation}%
which possesses the obvious inverse relation $\dot{s}^{\alpha }=u^{\alpha }-%
\mathcal{A}_{I}^{\alpha }\left( q^{a}\right) u^{I}$, analogously to (\ref%
{solConstr}), and thus%
\begin{equation}
\dot{\mathbf{q}}=\bar{\mathbf{B}}\left( \mathbf{q}\right) \mathbf{u}=\left( 
\begin{array}{cc}
\mathbf{I} & \ \mathbf{B}_{2} \\ 
\mathbf{0} & \mathbf{I}%
\end{array}%
\right) \mathbf{u}.  \label{Bbar}
\end{equation}%
%TCIMACRO{\TeXButton{red}{\color[rgb]{0,0,0}}}%
%BeginExpansion
\color[rgb]{0,0,0}%
%EndExpansion
Noting the specific structure of (\ref{Abar}) and (\ref{Bbar}), and that
only $\mathcal{A}_{I}^{\alpha }$ depends on $\mathbf{q}$, the corresponding
Hamel coefficients are found from (\ref{gammaNH}) as $\gamma _{\beta
J}^{\alpha }=-\gamma _{J\beta }^{\alpha }=-\frac{\partial \mathcal{A}%
_{J}^{\alpha }}{\partial s^{\beta }},$ and%
%TCIMACRO{\TeXButton{black}{\color{black}} }%
%BeginExpansion
\color{black}
%EndExpansion
\begin{equation}
\gamma _{IJ}^{\alpha }=\frac{\partial \mathcal{A}_{I}^{\alpha }}{\partial
r^{J}}-\frac{\partial \mathcal{A}_{J}^{\alpha }}{\partial r^{I}}+\frac{%
\partial \mathcal{A}_{J}^{\alpha }}{\partial s^{\beta }}\mathcal{A}%
_{I}^{\beta }-\frac{\partial \mathcal{A}_{I}^{\alpha }}{\partial s^{\beta }}%
\mathcal{A}_{J}^{\beta },\ \ \alpha =1,\ldots ,\bar{m};\ I,J=\bar{m}%
+1,\ldots ,n.  \label{gamma2}
\end{equation}

%TCIMACRO{\TeXButton{red}{\color[rgb]{0,0,0}}}%
%BeginExpansion
\color[rgb]{0,0,0}%
%EndExpansion
The configuration space $Q={\mathbb{V}}^{\bar{m}}\times {\mathbb{V}}^{\bar{%
\delta}}$ is regarded as a trivial bundle\footnote{%
%TCIMACRO{\TeXButton{red}{\color[rgb]{0,0,0}}}%
%BeginExpansion
\color[rgb]{0,0,0}%
%EndExpansion
$Q=M\times F$ is a trivial bundle if it can be written as Cartesian product
of a manifold $M$ and $F$, and if there is a projection $\pi :Q\rightarrow M$
\cite{BaezMuniain1994}. $F$ is called the standard fiber. {}For the
considered systems, the base manifold is the coordinate subspace ${\mathbb{V}%
}^{\bar{\delta}}$ corresponding to the independent velocities, and the fiber
is the subspace ${\mathbb{V}}^{\bar{m}}$ corresponding to dependent
velocities. 
%TCIMACRO{\TeXButton{black}{\color{black}}}%
%BeginExpansion
\color{black}%
%EndExpansion
} over the base manifold ${\mathbb{V}}^{\bar{\delta}}$ with fiber ${\mathbb{V%
}}^{\bar{m}}$, and bundle coordinates $(s^{\alpha },r^{I})\in {\mathbb{V}}^{%
\bar{m}}\times {\mathbb{V}}^{\bar{\delta}}$. 
%TCIMACRO{\TeXButton{black}{\color{black}}}%
%BeginExpansion
\color{black}%
%EndExpansion
The horizontal space of this trivial bundle is the constraint distribution,
i.e. the vector space of velocities satisfying the constraints. The
homogenous kinematic constraints (\ref{PfaffianConstr3}) define a connection
on this bundle. Writing the constraints in terms of the Pfaffian forms $%
\omega ^{\alpha }:=u^{\alpha }dt=ds^{\alpha }+\mathcal{A}_{I}^{\alpha
}dr^{I} $, a connection is introduced as $\mathcal{A}=\omega ^{\alpha }\frac{%
\partial }{\partial s^{\alpha }}$. This is referred to as an Ehresmann
connection \cite%
{BlochKrishnaprasadMarsdenMurray1996,BlochBook2003,MarsdenBook1995} with
reference to the original publication \cite{Ehresmann1950}, and $\mathcal{A}%
_{I}^{\alpha }$ are the local coordinates of the connection. Since it arises
from the kinematic constraints, it is called the \emph{kinematic connection} 
\cite{BlochKrishnaprasadMarsdenMurray1996}. The connection relates
(independent) motions in the base manifold ${\mathbb{V}}^{\bar{\delta}}$ to
motions in the fiber. Whether this connection (i.e. the constraints) is
holonomic is revealed by its curvature, denoted $\mathcal{B}_{IJ}^{\alpha }$%
. Moreover, the curvature of the kinematic connection plays a key role in
the control of constrained mechanical systems \cite%
{BlochKrishnaprasadMarsdenMurray1996,BlochBook2003,BlochMarsdenZenkov2005} 
%TCIMACRO{\TeXButton{red}{\color[rgb]{0,0,0}}}%
%BeginExpansion
\color[rgb]{0,0,0}%
%EndExpansion
as well as in locomotion planning \cite{HattonChosetIJRS2011} as it encodes
how motions in the base manifold generate motions in the fiber. On the
trivial vector bundle, the Lie bracket in the curvature (\ref{curvature}) is
the Lie bracket $[\mathcal{A}_{J},\mathcal{A}_{I}]^{\alpha }=\frac{\partial 
\mathcal{A}_{J}^{\alpha }}{\partial s^{\beta }}\mathcal{A}_{I}^{\beta }-%
\frac{\partial \mathcal{A}_{I}^{\alpha }}{\partial s^{\beta }}\mathcal{A}%
_{J}^{\beta }$ of vector fields $\mathcal{A}_{I},\mathcal{A}_{J}$ on $Q$, so
that the local curvature is $\mathcal{B}_{IJ}^{\alpha }=\frac{\partial 
\mathcal{A}_{I}^{\alpha }}{\partial r^{J}}-\frac{\partial \mathcal{A}%
_{J}^{\alpha }}{\partial r^{I}}+[\mathcal{A}_{J},\mathcal{A}_{I}]^{\alpha }$ 
\cite[p. 32]{BlochKrishnaprasadMarsdenMurray1996},\cite[p. 108]%
{BlochBook2003}. The Hamel coefficients (\ref{gamma2}) are thus clearly
related to the coordinate form of the curvature as follows.

\begin{proposition}
The Hamel coefficients (\ref{gamma2}) are identical to the curvature
components of the kinematic connection $\mathcal{A}$, in bundle coordinates $%
(s^{I},r^{\alpha })$, induced by the constraints with (\ref{PfaffianConstr3}%
), i.e. $\mathcal{B}_{IJ}^{\alpha }=\gamma _{IJ}^{\alpha }$.%
%TCIMACRO{\TeXButton{black}{\color{black}}}%
%BeginExpansion
\color{black}%
%EndExpansion
\end{proposition}

%TCIMACRO{\TeXButton{red}{\color[rgb]{0,0,0}}}%
%BeginExpansion
\color[rgb]{0,0,0}%
%EndExpansion
Although local coordinates are used in this paper, it should be mentioned
that the curvature of a connection $\mathcal{A}$ is its covariant
derivative, written coordinate-free as $\mathcal{B}\left( X,Y\right) =%
\mathrm{d}\mathcal{A}\left( X,Y\right) -[\mathcal{A}\left( X\right) ,%
\mathcal{A}\left( Y\right) ]$, with horizontal vector fields $X,Y$, i.e. $%
\mathcal{B}^{\alpha }\left( X,Y\right) =\mathcal{B}_{IJ}^{\alpha
}X^{I}Y^{J}=\gamma _{IJ}^{\alpha }X^{I}Y^{J}$.%
%TCIMACRO{\TeXButton{black}{\color{black}}}%
%BeginExpansion
\color{black}%
%EndExpansion

\begin{remark}
%TCIMACRO{\TeXButton{remChaplygin}{\label{remChaplygin}}}%
%BeginExpansion
\label{remChaplygin}%
%EndExpansion
The constraints are holonomic if and only if the kinematic connection is
flat, i.e. the components $\mathcal{B}_{IJ}^{\alpha }$ of the curvature
2-form vanish identically. In this case, ${\mathbb{V}}^{\bar{\delta}}$
serves as configuration space. Constraints are said to be in Chaplygin form
if $\mathcal{A}=\mathcal{A}\left( \mathbf{r}\right) $, referring to
Chaplygin's publications \cite{Chaplygin1897a,Chaplygin1897b}. In this case,
the Hamel-coefficients (\ref{BHconstr}) reduce to $\gamma _{IJ}^{\alpha }=%
\frac{\partial A_{I}^{\alpha }}{\partial r^{J}}-\frac{\partial A_{J}^{\alpha
}}{\partial r^{I}}$, which implies the obvious condition $\frac{\partial
A_{I}^{\alpha }}{\partial r^{J}}\equiv \frac{\partial A_{J}^{\alpha }}{%
\partial r^{I}}$ for integrability of $ds^{\alpha }+\mathcal{A}_{I}^{\alpha
}(r^{J})dr^{I}=0$.
\end{remark}

\subsection{Hamel Equations on a Trivial Bundle, Lagrange--d'Alembert
equations}

%TCIMACRO{\TeXButton{red}{\color[rgb]{0,0,0}}}%
%BeginExpansion
\color[rgb]{0,0,0}%
%EndExpansion
The Hamel coefficients can now be identified with the components of the
local curvature. Then the Hamel equations of the constrained system in
holonomic bundle coordinates $(r^{I},s^{\alpha })$ follow from (\ref{HE-holo}%
) as%
%TCIMACRO{\TeXButton{black}{\color{black}}}%
%BeginExpansion
\color{black}%
%EndExpansion
\begin{equation}
\frac{d}{dt}\frac{\partial \mathcal{L}}{\partial \dot{r}^{I}}-\frac{\partial 
\mathcal{L}}{\partial r^{I}}+\frac{\partial L}{\partial \dot{s}^{\alpha }}%
\mathcal{B}_{IJ}^{\alpha }\dot{r}^{J}=Q_{I},\ \ I=\bar{m}+1,\ldots ,n
\label{HE-holo2}
\end{equation}%
with $Q_{I}=\mathcal{A}_{I}^{a}Q_{a}$, in which $u^{\alpha }$ is set to
zero, and $u^{I}$ is replaced by $\dot{r}^{I}$. The dynamic equations are
written in terms of the coordinates $r^{I}$ on the base manifold of the
bundle. The remaining coordinates, the fiber coordinates, are obtained as
solution of 
\begin{equation}
\dot{s}^{\alpha }=-\mathcal{A}_{I}^{\alpha }\left( q^{a}\right) \dot{r}^{I}.
\label{recs}
\end{equation}%
The constrained dynamics is governed by the Hamel equations (\ref{HE-holo2})
along with the kinematic equations (\ref{recs}). 
%TCIMACRO{\TeXButton{red}{\color[rgb]{0,0,0}}}%
%BeginExpansion
\color[rgb]{0,0,0}%
%EndExpansion
The equations (\ref{HE-holo2}) are obtained as constrained
Lagrange--d'Alembert equations with variations satisfying the constraints $%
0=\delta s^{\alpha }+\mathcal{A}_{I}^{\alpha }(r^{I})\delta r^{I}$ \cite%
{Marsden1992,BlochBook2003,BlochKrishnaprasadMarsdenMurray1996}, which is a
particular form of Hamel's equations when using a connection to introduce
constraints.%
%TCIMACRO{\TeXButton{black}{\color{black}}}%
%BeginExpansion
\color{black}%
%EndExpansion

\begin{remark}
Introduce the \emph{constrained Lagrangian} $L_{\mathrm{c}}(s^{\alpha
},r^{I},\dot{r}^{I}):=L(s^{\alpha },r^{I},-\mathcal{A}_{I}^{\alpha }\dot{r}%
^{I},\dot{r}^{I})=\mathcal{L}(s^{\alpha },r^{I},u^{\alpha }:=0,u^{I}:=\dot{r}%
^{I})$, i.e. the Lagrangian with the constraints resolved. The Hamel
equations (\ref{HE-holo2}) attain the instructive form%
\begin{equation}
\frac{d}{dt}\frac{\partial L_{\mathrm{c}}}{\partial \dot{r}^{I}}-\frac{%
\partial L_{\mathrm{c}}}{\partial r^{I}}+\frac{\partial L_{\mathrm{c}}}{%
\partial s^{\alpha }}\mathcal{A}_{I}^{\alpha }+\frac{\partial L}{\partial 
\dot{s}^{\alpha }}\gamma _{IJ}^{\alpha }\dot{r}^{J}=Q_{I},\ \ I=\bar{m}%
+1,\ldots ,n,  \label{BHconstrHolo}
\end{equation}%
which reveals the consequence of non-holonomicity of the constraints. They
were reported in \cite[p. 326]{BlochMarsdenZenkov2005} (setting $\gamma
_{IJ}^{\alpha }=\mathcal{B}_{IJ}^{\alpha }$), and in a similar form for
constraints independent of $s^{I}$ by Chaplygin \cite%
{Chaplygin1897a,PapastavridisBook2002}. The Hamel coefficients reveal the
consequence of non-holonomic constraints. Clearly, if the constraints are
completely holonomic, these are the classical Lagrange equations. A direct
calculation shows that (\ref{BHconstrHolo}) can be written in the form 
\begin{equation}
\frac{d}{dt}\frac{\partial L}{\partial \dot{r}^{I}}-\frac{\partial L}{%
\partial r^{I}}+\mathcal{A}_{I}^{a}\left( \frac{d}{dt}\frac{\partial L}{%
\partial \dot{s}^{\alpha }}-\frac{\partial L}{\partial s^{\alpha }}\right) =%
\mathcal{A}_{I}^{a}Q_{a},\ \ I=\bar{m}+1,\ldots ,n
\end{equation}%
which have been reported by Voronets \cite{Voronets1901,Woronetz1910}, and
are referred to as Voronets equations \cite%
{Maisser1997,SoltakhanovYushkovZegzhda2009}.
\end{remark}

\section{%
%TCIMACRO{\TeXButton{red}{\color[rgb]{0,0,0}}}%
%BeginExpansion
\color[rgb]{0,0,0}%
%EndExpansion
Kinematically Constrained Mechanical Systems with Symmetry\label%
{secSymConstr}}

Many kinematically constrained systems possess principal symmetries in the
sense that the kinematic constraints are invariant under the action of a
symmetry group $G$. The configuration space can then be regarded as a
principal bundle\footnote{%
%TCIMACRO{\TeXButton{red}{\color[rgb]{0,0,0}}}%
%BeginExpansion
\color[rgb]{0,0,0}%
%EndExpansion
A space $Q=G\times Q/G$ is a \emph{principle bundle}, where $G$ is a Lie
group, with Lie algebra $\mathfrak{g}$, acting free and proper on $Q$,
equipped with a projection $\pi :Q\rightarrow Q/G$ \cite{Frankel2004}. In a
local trivialization, the \emph{base space} $Q/G$ can be identified with a
manifold $B$ so that $Q=G\times B$. If this splitting is globally valid, $Q$
is a \emph{trivial principle bundle}. In a local trivialization, with local
coordinates $(\xi ^{\alpha },\dot{r}^{I})$, $B$ will be identified with ${%
\mathbb{V}}^{\bar{\delta}}$. For a given $\mathbf{r}\in {\mathbb{V}}^{\bar{%
\delta}}$, $\pi ^{-1}\left( \mathbf{r}\right) $ is the fiber over $\mathbf{r}
$. For mechanical systems, $\mathbf{r}$ describes the internal configuration
(shape) of the system, and fiber elements $g\in G$ represent the pose of a
base body.}. Moreover, this is a trivial principal bundle when the kinematic
constraints do not depend on group variables. Examples are mobile platforms
and manipulators, or locomotion systems, where $G$ is often a subgroup of $%
SE\left( 3\right) $, the group of rigid body (i.e. Euclidean) motions%
%TCIMACRO{\TeXButton{red}{\color[rgb]{0,0,0}}}%
%BeginExpansion
\color[rgb]{0,0,0}%
%EndExpansion
, and $g\in G$ describes the motion of a base body. 
%TCIMACRO{\TeXButton{black}{\color{black}}}%
%BeginExpansion
\color{black}%
%EndExpansion
The governing equations are the reduced Euler-Lagrange equations on the
principle bundle.\vspace{-2ex}

\subsection{Constraints in Terms of a Connection on a Trivial Principal
Bundle}

The (non-holonomic) constraints are assumed to be invariant under the action
of a Lie group $G$. Chaplygin systems are included as special case with
Abelian symmetry group. The \emph{configuration space} of the system is
regarded as a trivial principal bundle $Q=G\times {\mathbb{V}}^{\bar{\delta}%
} $ over the base manifold ${\mathbb{V}}^{\bar{\delta}}$ with fiber $G$. The
dimension of $G$ is assumed to be equal to the number of constraints, and is
denoted with $\bar{m}$ (to be consistent with the preceding section).
Typically, fiber elements $g\in G$ represent the overall configuration of
the system in ambient space, and are often called 'body coordinates' (or
'rigid coordinates'). The coordinates $\mathbf{r}=(r^{I})\in {\mathbb{V}}^{%
\bar{\delta}}$ represent the internal shape, and are called 'shape
coordinates' (or 'internal variables'), and ${\mathbb{V}}^{\bar{\delta}}$ is
called the \emph{shape space}. 
%TCIMACRO{\TeXButton{red}{\color[rgb]{0,0,0}}}%
%BeginExpansion
\color[rgb]{0,0,0}%
%EndExpansion
For multibody systems, $\mathbf{r}$ is the vector of joint variables. Notice
that $r^{I}$ are only locally valid coordinates in general. 
%TCIMACRO{\TeXButton{black}{\color{black}}}%
%BeginExpansion
\color{black}%
%EndExpansion
The case when the number of constraints is less then the dimension of the
symmetry group $G$ has been addressed for motion planning of
non-holonomically constrained mechanical control systems in \cite%
{ShammasChosetRizza2007}.

Kinematic constraints that are left-invariant under actions of $G$ are
expressed as Pfaffian system $u^{\alpha }=0,\alpha =1,\ldots ,\bar{m}$, with 
\begin{equation}
u^{\alpha }:=\xi ^{\alpha }+\mathcal{A}_{I}^{\alpha }(r^{I})\dot{r}^{I}
\label{KinConnect}
\end{equation}%
where $\hat{%
%TCIMACRO{\TeXButton{xi}{\bm{\xi}}}%
%BeginExpansion
\bm{\xi}%
%EndExpansion
}=g^{-1}\dot{g}\in \mathfrak{g}$, and $%
%TCIMACRO{\TeXButton{xi}{\bm{\xi}}}%
%BeginExpansion
\bm{\xi}%
%EndExpansion
=\left( \xi ^{\alpha }\right) \in {\mathbb{R}}^{n}$ are the fiber
coordinates in a left-trivialization. More precisely, for rigid body systems
(where $G=SE\left( 3\right) $), $%
%TCIMACRO{\TeXButton{xi}{\bm{\xi}}}%
%BeginExpansion
\bm{\xi}%
%EndExpansion
$ is the velocity (also called twist) of a reference body in 'body-fixed'
representation. The constraints (\ref{KinConnect}) give rise to a principal
connection on the trivial principal bundle, denoting $\mathcal{A}d\mathbf{r}%
=(\mathcal{A}_{I}^{\alpha }dr^{I})$,%
\begin{equation}
\mathcal{A}^{\mathrm{kin}}=Ad_{g}(g^{-1}dg+\mathcal{A}d\mathbf{r}%
)=dgg^{-1}+Ad_{g}(\mathcal{A}d\mathbf{r})  \label{Akin}
\end{equation}%
so that the horizontal subspace of the connection is the space of velocities
satisfying the constraints, and $\mathcal{A}^{\mathrm{kin}}$ is a $\mathfrak{%
g}$-valued one-form (which may be considered as a special type of Ehresmann
connection) called the \emph{kinematic connection} as it arises from (\ref%
{KinConnect}) by requiring it to be $G$-equivariant \cite{BlochBook2003}.
The name stems from the fact that it relates base and fiber motions
according to the kinematic constraints. It is sufficient to use the local
connection form $\mathcal{A}_{I}^{\alpha }$ in (\ref{KinConnect}) as it
encodes all relevant information. The symbol $\mathcal{A}^{\mathrm{kin}}$ is
used to distinguish it from the coefficients of the local form $\mathcal{A}%
_{I}^{\alpha }$. Next, the Hamel coefficients are derived and are identified
as the coefficients of the curvature of the kinematic connection.\vspace{-1ex%
}

\subsection{The Hamel Coefficients and the Kinematic Connection\label%
{secHamelCoeff}}

The dynamics of the constrained system is governed by the reduced 
%TCIMACRO{\TeXButton{red}{\color[rgb]{0,0,0}}}%
%BeginExpansion
\color[rgb]{0,0,0}%
%EndExpansion
Lagrange-d'Alembert-Poincar\'{e} equations%
%TCIMACRO{\TeXButton{black}{\color{black}}}%
%BeginExpansion
\color{black}%
%EndExpansion
, which have been derived from the variational principle \cite%
{MarsdenScheurle1993,MarsdenRatiuScheurle2000}. To derive them as the
constrained Hamel equations necessitates the corresponding Hamel
coefficients.

\begin{lemma}
\label{lemKiCon}The non-vanishing Hamel coefficients in (\ref{BHconstr}) for
the system subjected to left $G$-invariant constraints are%
\begin{equation}
\gamma _{IJ}^{\alpha }=\frac{\partial \mathcal{A}_{I}^{\alpha }}{\partial
r^{J}}-\frac{\partial \mathcal{A}_{J}^{\alpha }}{\partial r^{I}}+c_{\lambda
\mu }^{\alpha }\mathcal{A}_{I}^{\lambda }\mathcal{A}_{J}^{\mu },\ \ \alpha
=1,\ldots ,\bar{m}.  \label{gamma3}
\end{equation}
\end{lemma}

\begin{proof}
In order to apply the original definition (\ref{gamma}) of the Hamel
coefficients, 
%TCIMACRO{\TeXButton{red}{\color[rgb]{0,0,0}}}%
%BeginExpansion
\color[rgb]{0,0,0}%
%EndExpansion
local 
%TCIMACRO{\TeXButton{black}{\color{black}}}%
%BeginExpansion
\color{black}%
%EndExpansion
canonical coordinates $s^{\alpha },\alpha =1,\ldots ,\bar{m}$ are introduced
on $G$. The fiber coordinates are then expressed as $\xi ^{\alpha }=A_{\beta
}^{\alpha }\left( s^{\alpha }\right) \dot{s}^{\beta }$, with inverse
relation $\dot{s}^{\alpha }=B_{\beta }^{\alpha }\left( s^{\alpha }\right)
\xi ^{\beta }$, and (\ref{KinConnect}) is written as%
\begin{equation}
u^{\alpha }=A_{\beta }^{\alpha }\left( s^{\alpha }\right) \dot{s}^{\beta }+%
\mathcal{A}_{I}^{\alpha }(r^{J})\dot{r}^{I}.
\end{equation}%
This resembles the relation (\ref{PfaffianConstr2}) with $A_{\beta }^{\alpha
}=A_{\beta }^{\alpha }\left( s^{\alpha }\right) $ and $A_{I}^{\alpha
}=A_{I}^{\alpha }(r^{J})$. The inverse relation is 
\begin{equation}
\dot{s}^{\alpha }=B_{\beta }^{\alpha }\left( s^{a}\right) u^{\beta
}-B_{\beta }^{\alpha }\left( s^{a}\right) \mathcal{A}_{I}^{\beta
}(r^{J})u^{I}
\end{equation}%
with $B_{\beta }^{\alpha }$ defined in (\ref{sys2}). Noting that $A_{a}^{I}=%
\mathrm{const}$, the only non-zero Hamel coefficients for the constrained
system are $\gamma _{IJ}^{\alpha }$. They are immediately found from (\ref%
{gammaHolo}), by replacing $A_{I}^{\alpha }(r^{J})$ with $\mathcal{A}%
_{I}^{\alpha }(r^{J})$, and $B_{I}^{\alpha }$ with $-B_{\beta }^{\alpha
}\left( s^{a}\right) \mathcal{A}_{I}^{\beta }(r^{J})$, as%
\begin{equation}
\gamma _{IJ}^{\alpha }=\frac{\partial \mathcal{A}_{I}^{\alpha }}{\partial
r^{J}}-\frac{\partial \mathcal{A}_{J}^{\alpha }}{\partial r^{I}}+\left( 
\frac{\partial A_{\delta }^{\alpha }}{\partial s^{\lambda }}-\frac{\partial
A_{\lambda }^{\alpha }}{\partial s^{\delta }}\right) B_{\gamma }^{\delta
}B_{\mu }^{\lambda }\mathcal{A}_{I}^{\gamma }\mathcal{A}_{J}^{\mu }.
\label{gammatmp}
\end{equation}%
%TCIMACRO{\TeXButton{red}{\color[rgb]{0,0,0}}}%
%BeginExpansion
\color[rgb]{0,0,0}%
%EndExpansion
It was already shown by Hamel \cite[p. 428]{Hamel-MathAnal1904} that the
terms $\gamma _{\beta \lambda }^{\alpha }=\left( \frac{\partial A_{\mu
}^{\alpha }}{\partial s^{\nu }}-\frac{\partial A_{\nu }^{\alpha }}{\partial
s^{\mu }}\right) B_{\beta }^{\mu }B_{\lambda }^{\nu }=c_{\beta \lambda
}^{\alpha }$ are the structure constants of $G$. A more recent reference is 
\cite[p. 301]{McCauleyBook}. 
%TCIMACRO{\TeXButton{black}{\color{black}}}%
%BeginExpansion
\color{black}%
%EndExpansion
This was derived explicitly in \cite{MuellerJNLS2021} for $SO\left( 3\right) 
$ and$SE\left( 3\right) $ using the original definition of Hamel
coefficients.
\end{proof}

The expressions (\ref{gamma3}) are the coefficients of the local curvature
of the connection $\mathcal{A}^{\mathrm{kin}}$ on the principal bundle
possibly up to a change of sign \cite{BaezMuniain1994,Choquet1989} (for the
sign convention see Rem. \ref{remSign}). In context of geometric mechanics,
this can be stated as follows \cite[notice the correction on p. 44]%
{BlochKrishnaprasadMarsdenMurray1996}.

\begin{proposition}
%TCIMACRO{\TeXButton{propHC1}{\label{propHC1}}}%
%BeginExpansion
\label{propHC1}%
%EndExpansion
The Hamel coefficients (\ref{gamma3}) for left-invariant kinematic
constraints (\ref{KinConnect}) are the components of the curvature (\ref%
{curvature}) of the kinematic connection defined in bundle coordinates $(\xi
^{\alpha },\dot{r}^{I})$ on the corresponding left-trivialized principal
bundle: $\mathcal{B}_{IJ}^{\alpha }=\gamma _{IJ}^{\alpha }$. Noting that the
curvature coefficients are identical to the Hamel coefficients, the
connection is flat, if and only if the constraints are holonomic.
\end{proposition}

%TCIMACRO{\TeXButton{red}{\color[rgb]{0,0,0}}}%
%BeginExpansion
\color[rgb]{0,0,0}%
%EndExpansion
Also the curvature on the principle bundle can be defined coordinate-free as
covariant derivative $\mathcal{B}\left( X,Y\right) =\mathrm{d}\mathcal{A}%
\left( X,Y\right) -[\mathcal{A}\left( X\right) ,\mathcal{A}\left( Y\right) ]$%
, with horizontal vector fields $X,Y$, now with Lie bracket on $\mathfrak{g}$%
. In local coordinates, it is $\mathcal{B}^{\alpha }\left( X,Y\right)
=\gamma _{IJ}^{\alpha }X^{I}Y^{J}$.%
%TCIMACRO{\TeXButton{black}{\color{black}}}%
%BeginExpansion
\color{black}%
%EndExpansion

\begin{remark}
If 
%TCIMACRO{\TeXButton{red}{\color[rgb]{0,0,0}}}%
%BeginExpansion
\color[rgb]{0,0,0}%
%EndExpansion
right trivialization 
%TCIMACRO{\TeXButton{black}{\color{black}}}%
%BeginExpansion
\color{black}%
%EndExpansion
is used, i.e. $\hat{%
%TCIMACRO{\TeXButton{xi}{\bm{\xi}}}%
%BeginExpansion
\bm{\xi}%
%EndExpansion
}=\dot{g}g^{-1}\in \mathfrak{g}$ are right invariant vector fields (e.g.
using spatial velocities), then $\gamma _{\beta \lambda }^{\alpha }=\left( 
\frac{\partial A_{\mu }^{\alpha }}{\partial s^{\nu }}-\frac{\partial A_{\nu
}^{\alpha }}{\partial s^{\mu }}\right) B_{\beta }^{\mu }B_{\lambda }^{\nu
}=-c_{\beta \lambda }^{\alpha }$, and the Hamel coefficients (\ref{gamma3})
are%
\begin{equation}
\gamma _{IJ}^{\alpha }=\frac{\partial \mathcal{A}_{I}^{\alpha }}{\partial
r^{J}}-\frac{\partial \mathcal{A}_{J}^{\alpha }}{\partial r^{I}}-c_{\lambda
\mu }^{\alpha }\mathcal{A}_{I}^{\lambda }\mathcal{A}_{J}^{\mu },\ \ \alpha
=1,\ldots ,\bar{m}.  \label{gammaRight}
\end{equation}
\end{remark}

\begin{remark}
\label{remKinControl}%
%TCIMACRO{\TeXButton{red}{\color[rgb]{0,0,0}}}%
%BeginExpansion
\color[rgb]{0,0,0}%
%EndExpansion
The constraints (\ref{KinConnect}) give rise to a kinematic control system
on $G$ that can be written as $\xi ^{\alpha }=-\mathcal{A}_{I}^{\alpha
}(r^{I})\dot{r}^{I}$. Such control systems with symmetry were addressed in
various publications, e.g. \cite%
{KellyMurray1995,OstrowskiBurdick1996,ShammasChosetRizza2007,HattonChosetIJRS2011}%
. Controllability of this driftless control problem is encoded in the
control vector fields $\mathcal{A}_{I}^{\alpha }$ and the distribution on
the fiber defined by them. Written as $\dot{g}=-g\mathcal{A}d\dot{\mathbf{r}}
$ shows that the connection describes how a path in shape space is lifted to
a path in the group (called the horizontal lift of the curve), which is the
basis for kinematic control. The net change in the group variable as result
of the horizontal lift of a closed curve in shape space is the holonomy (in
this context called the \emph{geometric phase}). The latter can be related
to the curvature of the kinematic connection. Let $\phi $ be a closed path
in ${\mathbb{V}}^{\bar{\delta}}$. The geometric phase is then found as%
\begin{equation}
g\left( \phi \right) =-\oint_{\phi }g\mathcal{A}d\mathbf{r}=-\int_{\Phi }%
\mathcal{B}+\mathrm{hot}.  \label{int1}
\end{equation}%
using $\dot{g}=g\hat{%
%TCIMACRO{\TeXButton{xi}{\bm{\xi}}}%
%BeginExpansion
\bm{\xi}%
%EndExpansion
}=-g\mathcal{A}\dot{\mathbf{r}}$, where $\Phi $ is the area enclosed by the
path in ${\mathbb{V}}^{\bar{\delta}}$, and $\mathcal{B}=(\mathcal{B}%
_{IJ}^{\alpha }dr^{I}dr^{J})$. If the constraints are holonomic, i.e. the
curvature vanishes, the geometric phase is zero. This is an important
relation for locomotion planning, where the closed path $\phi $ represents a
gait, and the aim is to maximize the net motion in $G$ generated by a gait 
\cite{HattonChosetIJRS2011}. The geometric phase further discussed in Sec. %
\ref{secUnConDymConserved}.%
%TCIMACRO{\TeXButton{black}{\color{black}}}%
%BeginExpansion
\color{black}%
%EndExpansion
%TCIMACRO{\TeXButton{TeX field}{\vspace{-2ex}}}%
%BeginExpansion
\vspace{-2ex}%
%EndExpansion
\end{remark}

\subsection{The Hamel Equations%
%TCIMACRO{\TeXButton{red}{\color[rgb]{0,0,0}}}%
%BeginExpansion
\color[rgb]{0,0,0}%
%EndExpansion
, Lagrange--d'Alembert--Poincar\'{e} equations}

The Lagrangian $L(g,r^{I},\dot{g},\dot{r}^{I})$, defined on the
configuration space $Q$, is assumed to be (left or right) $G$-invariant so
that the reduced Lagrangian $l(r^{I},\xi ^{\alpha },\dot{r}^{I})$ can be
introduced on the corresponding bundle trivialization. Assume that the
motion in $G$ is completely determined by the motion in the base manifold,
i.e. the group orbits complement the constraints. This is called the
'principal kinematic case' \cite%
{BlochKrishnaprasadMarsdenMurray1996,BlochBook2003} since then there is no
momentum equation left on the fibre.

The Pfaffian equations (\ref{KinConnect}) can be resolved as $\xi ^{\alpha
}=u^{\alpha }-\mathcal{A}_{I}^{\alpha }\dot{r}^{I}$. In the Hamel formalism,
the $u^{\alpha }$ are regarded as intermediate coordinates, and setting $%
u^{\alpha }=0$ yields the constraint solution. Denote with $\ell
(r^{I},u^{\alpha },u^{I}):=l(r^{I},\xi ^{\alpha }:=u^{\alpha }-\mathcal{A}%
_{I}^{\alpha }\dot{r}^{I},\dot{r}^{I}:=u^{I})$ the Lagrangian in terms of $%
u^{\alpha }$ in (\ref{KinConnect}) and $\dot{r}^{I}:=u^{I}$. The Hamel
equations (\ref{BHconstr}) for the constrained system are expressed using
proposition \ref{propHC1} as%
\begin{equation}
\frac{d}{dt}\frac{\partial \ell }{\partial \dot{r}^{I}}-\frac{\partial \ell 
}{\partial r^{I}}=\frac{\partial l}{\partial u^{\beta }}\mathcal{B}%
_{IJ}^{\beta }\dot{r}^{J}+Q_{I},\ \ I=\bar{m}+1,\ldots ,n
\label{HE-kinBundle}
\end{equation}%
where $u^{\alpha }$ is set to zero, and $u^{I}=\dot{r}^{I}$. 
%TCIMACRO{\TeXButton{red}{\color[rgb]{0,0,0}}}%
%BeginExpansion
\color[rgb]{0,0,0}%
%EndExpansion
These are indeed equations (\ref{LP2}) when setting $\Omega ^{\alpha
}=u^{\alpha }=0$. The velocities in the fiber are obtained with the
kinematic connection as $\xi ^{\alpha }=-\mathcal{A}_{I}^{\alpha }\left( 
\mathbf{r}\right) \dot{r}^{I}$, which is the principle bundle equivalent of (%
\ref{recs}). In geometric mechanics, equations (\ref{HE-kinBundle}) are
obtained as Lagrange--d'Alembert--Poincar\'{e} equations with variations
satisfying the constraints $u^{\alpha }=0$ with (\ref{KinConnect}). 
%TCIMACRO{\TeXButton{black}{\color{black}}}%
%BeginExpansion
\color{black}%
%EndExpansion
The motion $g\left( t\right) $ in the fiber $G$ is obtained by solving the 
\emph{kinematic reconstruction equations}%
\begin{equation}
\dot{g}=g\hat{\bm{\xi}}\mathrm{\ \ (left\ trivialization)\ \ \ \ \ or\ \ \ \ 
}\dot{g}=\hat{\bm{\xi}}g\ \ \mathrm{(right\ trivialization)}  \label{KinRec}
\end{equation}%
where $\bm{\xi}=-\mathcal{A}\dot{\mathbf{r}}$. The Hamel equations (\ref%
{HE-kinBundle}) along with the reconstruction equations (\ref{KinRec})
govern the dynamics of the constrained system on the principle bundle $Q$.
Chaplygin systems (Rem. \ref{remChaplygin}) are special cases with Abelian
symmetry group $G={\mathbb{R}}^{\bar{m}}$, where $\mathcal{B}_{IJ}^{\alpha }=%
\frac{\partial \mathcal{A}_{J}^{\alpha }}{\partial r^{I}}-\frac{\partial 
\mathcal{A}_{I}^{\alpha }}{\partial r^{J}}$.

\begin{remark}
Equations (\ref{KinRec}) are the Poisson equations on $G$. They are known in
context of rigid body kinematics as the left- and right Poisson-Darboux
equations, referring to \cite{Darboux1887}, or as the generalized
Poisson-Darboux equations \cite{ConduracheAASAIAA2017}. In order to solve
these differential equations on $G$, 
%TCIMACRO{\TeXButton{red}{\color[rgb]{0,0,0}}}%
%BeginExpansion
\color[rgb]{0,0,0}%
%EndExpansion
$g$ is expressed as the exponential of a $\bm{\eta}\left( t\right) \in 
\mathfrak{g}$, and (\ref{KinRec}) are replaced by the following differential
equations on $\mathfrak{g}$%
%TCIMACRO{\TeXButton{black}{\color{black}} }%
%BeginExpansion
\color{black}
%EndExpansion
\begin{align}
%TCIMACRO{\TeXButton{red}{\color[rgb]{0,0,0}}}%
%BeginExpansion
\color[rgb]{0,0,0}%
%EndExpansion
\dot{\hat{\bm{\eta}}}& 
%TCIMACRO{\TeXButton{red}{\color[rgb]{0,0,0}}}%
%BeginExpansion
\color[rgb]{0,0,0}%
%EndExpansion
=\mathrm{dexp}_{-\hat{\bm{\eta}}}^{-1}(\hat{\bm{\xi}})=-\mathrm{dexp}_{-\hat{%
\bm{\eta}}}^{-1}(\mathcal{A}\dot{\mathbf{r}}),\ \mathrm{with}\ \ g=g_{0}\exp
(\hat{\bm{\eta}})\mathrm{\ \ (left\ trivialization)}  \label{KinRec2} \\
%TCIMACRO{\TeXButton{red}{\color[rgb]{0,0,0}}}%
%BeginExpansion
\color[rgb]{0,0,0}%
%EndExpansion
\dot{\hat{\bm{\eta}}}& 
%TCIMACRO{\TeXButton{red}{\color[rgb]{0,0,0}}}%
%BeginExpansion
\color[rgb]{0,0,0}%
%EndExpansion
=\mathrm{dexp}_{\hat{\bm{\eta}}}^{-1}(\hat{\bm{\xi}})=-\mathrm{dexp}_{\hat{%
\bm{\eta}}}^{-1}(\mathcal{A}\dot{\mathbf{r}}),\ \mathrm{with}\ \ g=\exp (%
\hat{\bm{\eta}})g_{0}\mathrm{\ }\ \ \mathrm{(right\ trivialization)}  \notag
\end{align}%
with initial value $g_{0}\in G$, where $\hat{\bm{\eta}}\left( t\right) \in 
\mathfrak{g}$ represents a local canonical parameterization of $G$, and $%
\mathrm{dexp}$ is the right-trivialized differential of the $\exp $ map on $%
G $. 
%TCIMACRO{\TeXButton{red}{\color[rgb]{0,0,0}}}%
%BeginExpansion
\color[rgb]{0,0,0}%
%EndExpansion
The latter is defined by $\dot{g}g^{-1}=\mathrm{dexp}_{\hat{\bm{\eta}}}(\hat{%
\bm{\eta}})$ assuming $g=\exp (\hat{\bm{\eta}})g_{0}$. 
%TCIMACRO{\TeXButton{black}{\color{black}}}%
%BeginExpansion
\color{black}%
%EndExpansion
This replacement is a key step in Lie group integration schemes \cite%
{MuntheKaas-BIT1998,MuntheKaas1999,IserlesMuntheKaasNrsettZanna2000}. For
many Lie groups relevant to solid mechanics, there are closed from
expressions for the $\mathrm{dexp}$ map, in particular for $SO\left(
3\right) $ and $SE\left( 3\right) $ \cite{RSPA2021}. 
%TCIMACRO{\TeXButton{red}{\color[rgb]{0,0,0}}}%
%BeginExpansion
\color[rgb]{0,0,0}%
%EndExpansion
The maps $\mathrm{dexp}$ and $\mathrm{dexp}^{-1}$ can be evaluated using
truncated series expansions. 
%TCIMACRO{\TeXButton{black}{\color{black}}}%
%BeginExpansion
\color{black}%
%EndExpansion
\end{remark}

\subsection{Example: Homogenous ball rolling without slipping or spinning}

As a simple example, consider a ball on a horizontal plane that is subjected
to pure rolling, i.e. in addition to the rolling-without-slipping constraint
it is further constrained so that it cannot spin about its instantaneous
vertical axis, which is parallel to the plane normal. Its configuration is
described by its orientation and the location of its point of contact with
the plane. The configuration space is the bundle $SO\left( 3\right) \times {%
\mathbb{R}}^{2}$ over $M={\mathbb{R}}^{2}$. The motion in the fiber $%
G=SO\left( 3\right) $ (the orientation) is completely determined by the
motion of the contact point so that this example is a 'principal kinematic
case'. Since the constraints and the potential energy (gravity) are right $G$%
-invariant, the kinetic energy is bi-invariant, the configuration space is
regarded as a right-trivialized trivial principal bundle. 
\begin{figure}[b]
\begin{center}
\includegraphics[draft=false,height=4.0cm]{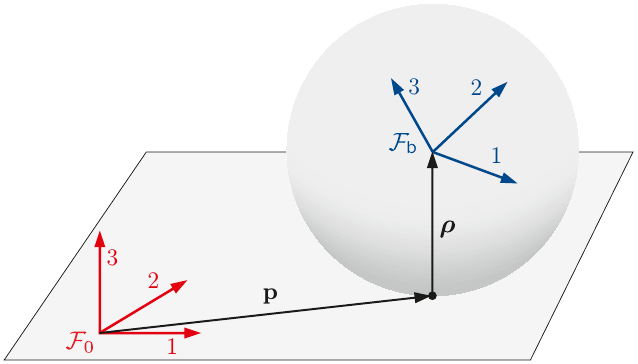}
\end{center}
\par
\vspace{-2ex}
\caption{Ball rolling without spinning about the vertical axis.}
\label{figBall}
\end{figure}
It is assumed that the center of mass (COM) of the ball is at its geometric
center. Denote with $%
%TCIMACRO{\TeXButton{rho}{\bm{\rho}}}%
%BeginExpansion
\bm{\rho}%
%EndExpansion
\in {\mathbb{R}}^{3}$ the vector from contact point to the COM of the ball,
with $\mathbf{p}\in {\mathbb{R}}^{3}$ the position of contact point, and $%
\mathbf{R}\in SO\left( 3\right) $ describes the rotation of the body-fixed
frame $\mathcal{F}_{\mathrm{b}}$ relative to the inertial frame $\mathcal{F}%
_{0}$. The angular velocity $%
%TCIMACRO{\TeXButton{w}{\bm{\omega}}}%
%BeginExpansion
\bm{\omega}%
%EndExpansion
\in {\mathbb{R}}^{3}\cong so\left( 3\right) $ of the ball in spatial
representation is defined by $\hat{%
%TCIMACRO{\TeXButton{w}{\bm{\omega}}}%
%BeginExpansion
\bm{\omega}%
%EndExpansion
}=\dot{\mathbf{R}}\mathbf{R}^{T}\in so\left( 3\right) $. All vectors are
expressed in inertial frame $\mathcal{F}_{0}$. This frame is introduced such
that its 3-axis is parallel and directed along the plane normal, so that $%
%TCIMACRO{\TeXButton{rho}{\bm{\rho}}}%
%BeginExpansion
\bm{\rho}%
%EndExpansion
=\left( 0,0,R\right) ^{T}$, where $R$ is the ball radius.

The rolling condition is $\mathbf{0}=\dot{\mathbf{p}}+\widetilde{%
%TCIMACRO{\TeXButton{rho}{\bm{\rho}}}%
%BeginExpansion
\bm{\rho}%
%EndExpansion
}%
%TCIMACRO{\TeXButton{w}{\bm{\omega}}}%
%BeginExpansion
\bm{\omega}%
%EndExpansion
$, which can be transformed to $\mathbf{0}=\widetilde{%
%TCIMACRO{\TeXButton{rho}{\bm{\rho}}}%
%BeginExpansion
\bm{\rho}%
%EndExpansion
}\dot{\mathbf{p}}+R^{2}%
%TCIMACRO{\TeXButton{w}{\bm{\omega}}}%
%BeginExpansion
\bm{\omega}%
%EndExpansion
-(%
%TCIMACRO{\TeXButton{rho}{\bm{\rho}}}%
%BeginExpansion
\bm{\rho}%
%EndExpansion
%TCIMACRO{\TeXButton{rho}{\bm{\rho}}}%
%BeginExpansion
\bm{\rho}%
%EndExpansion
^{T})%
%TCIMACRO{\TeXButton{w}{\bm{\omega}}}%
%BeginExpansion
\bm{\omega}%
%EndExpansion
$. 
%TCIMACRO{\TeXButton{red}{\color[rgb]{0,0,0}}}%
%BeginExpansion
\color[rgb]{0,0,0}%
%EndExpansion
Here $\widetilde{%
%TCIMACRO{\TeXButton{rho}{\bm{\rho}}}%
%BeginExpansion
\bm{\rho}%
%EndExpansion
}$ is the skew symmetric matrix so that the cross product of vector $\mathbf{%
x}$ and $\mathbf{y}$ is $\widetilde{\mathbf{x}}\mathbf{y}=\mathbf{x}\times 
\mathbf{y}$ that is written in components as $\varepsilon _{ijk}x^{j}y^{k}$,
with the Levi-Civita symbol $\varepsilon _{ijk}$. 
%TCIMACRO{\TeXButton{black}{\color{black}}}%
%BeginExpansion
\color{black}%
%EndExpansion
The non-spinning condition is $0=%
%TCIMACRO{\TeXButton{rho}{\bm{\rho}}}%
%BeginExpansion
\bm{\rho}%
%EndExpansion
^{T}%
%TCIMACRO{\TeXButton{w}{\bm{\omega}}}%
%BeginExpansion
\bm{\omega}%
%EndExpansion
$, which amounts to the constraint $\omega ^{3}=0$, and the rolling
condition simplifies to $\mathbf{0}=%
%TCIMACRO{\TeXButton{w}{\bm{\omega}}}%
%BeginExpansion
\bm{\omega}%
%EndExpansion
+\frac{1}{R^{2}}\widetilde{%
%TCIMACRO{\TeXButton{rho}{\bm{\rho}}}%
%BeginExpansion
\bm{\rho}%
%EndExpansion
}\dot{\mathbf{p}}$. Since $%
%TCIMACRO{\TeXButton{rho}{\bm{\rho}}}%
%BeginExpansion
\bm{\rho}%
%EndExpansion
$ is along the plane normal, the rolling condition is $0=\frac{1}{R^{2}}%
\varepsilon _{\alpha JI}\rho ^{J}=\frac{1}{R}\varepsilon _{\alpha 3I}$. The
'shape coordinates' $\mathbf{r}=(r^{4},r^{5})=(p^{1},p^{2})\in {\mathbb{R}}%
^{2}$ are the coordinates of the contact location in the plane. For sake of
compactness, denote $\bar{I}=I-3$. Local bundle coordinates are $\omega
^{\alpha },\alpha =1,2,3$ and $\dot{r}^{I}=\dot{p}^{\bar{I}},I=4,5$. The $%
\bar{m}=3$ kinematic constraints are $0=u^{\alpha },\alpha =1,2,3$, 
\begin{equation}
u^{\alpha }:=\omega ^{\alpha }+\mathcal{A}_{I}^{\alpha }\dot{r}^{I},\ \ \ 
\mathcal{A}_{I}^{\alpha }:=\left\{ 
\begin{array}{cll}
\frac{1}{R}\varepsilon _{\alpha 3\bar{I}}, & \alpha =1,2 &  \\ 
0, & \alpha =3 & .%
\end{array}%
\right.  \label{uaBall}
\end{equation}%
The Hamel coefficients are $\gamma _{IJ}^{\alpha }=-c_{\beta \mu }^{\alpha }%
\mathcal{A}_{I}^{\beta }\mathcal{A}_{I}^{\mu },\ \ I,J=4,5$ (negative sign
is due to right-trivialization), which is only non-zero for $\alpha =3$ as $%
\beta ,\mu =1,2$. Evaluation yields $c_{\beta \mu }^{\alpha }\mathcal{A}%
_{I}^{\beta }\mathcal{A}_{J}^{\mu }=\frac{1}{R^{4}}\varepsilon _{\alpha
\beta \mu }\varepsilon _{\beta r\bar{I}}\varepsilon _{\mu s\bar{J}}\rho
^{r}\rho ^{s}=\frac{1}{R^{4}}(\delta _{\mu r}\delta _{\alpha \bar{I}}-\delta
_{\mu \bar{I}}\delta _{\alpha r})\varepsilon _{\mu s\bar{J}}\rho ^{r}\rho
^{s}=-\frac{1}{R^{4}}\rho ^{\alpha }\varepsilon _{\bar{I}s\bar{J}}\rho ^{s}$ 
$=-\frac{1}{R^{2}}\rho ^{\alpha }\mathcal{A}_{J}^{\beta }\delta _{\beta }^{%
\bar{I}}$, thus%
\begin{equation}
\gamma _{IJ}^{\alpha }=\left\{ 
\begin{array}{cll}
0, & \alpha =1,2 &  \\ 
\frac{1}{R}\mathcal{A}_{J}^{\bar{I}}, & \alpha =3 & .%
\end{array}%
\right.
\end{equation}%
%TCIMACRO{\TeXButton{red}{\color[rgb]{0,0,0}}}%
%BeginExpansion
\color[rgb]{0,0,0}%
%EndExpansion
In summary, the non-zero coefficients are $\mathcal{A}_{4}^{2}=-\mathcal{A}%
_{5}^{1}=\frac{1}{R}$ and $\mathcal{B}_{54}^{3}=-\mathcal{B}_{45}^{3}=\gamma
_{54}^{3}=\frac{1}{R^{2}}$.%
%TCIMACRO{\TeXButton{black}{\color{black}}}%
%BeginExpansion
\color{black}%
%EndExpansion

The Lagrangian is identical to the kinetic energy (potential energy does not
affect the motion of the ball). The right-reduced Lagrangian is $l(\omega
^{\alpha },\dot{r}^{I})=\frac{m}{2}\delta _{IJ}\dot{r}^{I}\dot{r}^{J}+\frac{1%
}{2}\Theta _{\alpha \beta }\omega ^{\alpha }\omega ^{\beta }$, where $\Theta
_{\alpha \beta }$ is the inertia tensor of the homogenous ball w.r.t. its
COM expressed in the spatial inertial frame $\mathcal{F}_{0}$, and $m$ is
its mass. The Lagrangian in terms of $u^{I}=\dot{r}^{I}$ and $u^{\alpha }$
in (\ref{uaBall}) is%
\begin{equation}
\ell (u^{\alpha },u^{I})=l(\omega ^{\alpha }:=u^{\alpha }-\mathcal{A}%
_{J}^{\alpha }u^{I},\dot{r}^{I}:=u^{I})=\frac{m}{2}\delta _{IJ}u^{I}u^{J}+%
\frac{1}{2}\Theta _{\alpha \beta }(u^{\alpha }-\mathcal{A}_{I}^{\alpha
}u^{I})(u^{\beta }-\mathcal{A}_{J}^{\beta }u^{J}).
\end{equation}%
The Hamel equations (\ref{HE-kinBundle}), with $Q_{I}=0$, are found as (with 
$r^{I}=p^{\bar{I}}$)%
\begin{equation}
(m\delta _{IJ}-\Theta _{\alpha \beta }\mathcal{A}_{I}^{\alpha }\mathcal{A}%
_{J}^{\beta })\ddot{r}^{J}-\frac{1}{R}\Theta _{3\beta }\mathcal{A}%
_{K}^{\beta }\dot{r}^{K}\mathcal{A}_{J}^{\bar{I}}\dot{r}^{J}=0,\ \ I=4,5
\label{HEBall}
\end{equation}%
where $\alpha ,\beta \in \{1,2\}$. 
%TCIMACRO{\TeXButton{red}{\color[rgb]{0,0,0}}}%
%BeginExpansion
\color[rgb]{0,0,0}%
%EndExpansion
The second term in (\ref{HEBall}) is the force due to the gyroscopic torque
caused by the pure rolling. To see this, notice that $\pi _{3}=\Theta
_{3\beta }\mathcal{A}_{K}^{\beta }\dot{r}^{K}$ is the vertical component of
the angular momentum $\pi _{\alpha }=\Theta _{\alpha \beta }\omega ^{\beta }$
due to the rolling (note $\omega ^{3}=0$), and the gyroscopic torque is $%
\varepsilon _{\alpha \beta \lambda }\omega ^{\beta }\pi _{\lambda }$. The
force is obtained as cross product of this torque with $%
%TCIMACRO{\TeXButton{rho}{\bm{\rho}}}%
%BeginExpansion
\bm{\rho}%
%EndExpansion
$, noting that $\pi _{1}=\pi _{3}=0$ and $\omega _{3}=0$.%
%TCIMACRO{\TeXButton{black}{\color{black}}}%
%BeginExpansion
\color{black}%
%EndExpansion

\begin{remark}
The principal kinematic case can be extended by allowing the Lagrangian to
be invariant under actions of a subgroup of the symmetry group $G$ of the
constraints, i.e. infinitesimal generators of this subgroup lie in the
constraint distribution. Then the system is said to possess horizontal
symmetries (relative to the constraints) \cite%
{BlochKrishnaprasadMarsdenMurray1996,BlochBook2003,KoonMarsden1997}. In this
case, the kinematic connection along with a connection accounting for the
horizontal symmetry can be introduced. As an example for mechanical systems
with horizontal symmetries, a ball moving on a plane, i.e. a rolling ball
that is free to spin about the plane normal, was used in \cite%
{BlochKrishnaprasadMarsdenMurray1996,KoonMarsden1997}. The horizontal
symmetry is then due to the momentum being invariant under the subgroup $%
SO\left( 2\right) \subset SO\left( 3\right) $ of rotations about the plane
normal. The $\bar{m}$ rolling constraints are the 1- and 2-component of the
condition $\mathbf{0}=\dot{\mathbf{r}}+\widetilde{%
%TCIMACRO{\TeXButton{rho}{\bm{\rho}}}%
%BeginExpansion
\bm{\rho}%
%EndExpansion
}%
%TCIMACRO{\TeXButton{w}{\bm{\omega}}}%
%BeginExpansion
\bm{\omega}%
%EndExpansion
$ above. The invariant momentum is the 3-component of $\bm{\Theta}%
%TCIMACRO{\TeXButton{w}{\bm{\omega}}}%
%BeginExpansion
\bm{\omega}%
%EndExpansion
$.
\end{remark}

\section{%
%TCIMACRO{\TeXButton{red}{\color[rgb]{0,0,0}}}%
%BeginExpansion
\color[rgb]{0,0,0}%
%EndExpansion
Floating-Base Mechanical Systems with Symmetries\label{secUnConBundle}}

If the Lagrangian $L:TQ\rightarrow {\mathbb{R}}$ is $G$-invariant and can be
written as $L(g,r^{I},\dot{g},\dot{r}^{I})$, the $n$-dimensional
configuration space is regarded as a principal bundle, which admits the
local trivialization $Q=G\times Q/G$. This parallels the formulation of
constrained systems in Sec. \ref{secSymConstr}, where $G$ was the symmetry
group of the constraints. The quotient space $Q/G$ is again the \emph{shape
space}, which can locally be identified with ${\mathbb{V}}^{\bar{\delta}}$.
This setting applies to many mechanical systems possessing symmetry
invariants. Typical examples are floating multibody systems (humanoids,
areal vehicles, space robots), where the fiber $G$ is a subgroup of $%
SE\left( 3\right) $ representing the overall spatial motion, and $Q/G\cong {%
\mathbb{V}}^{\bar{\delta}}$ represents the internal shape. An early problem
that was considered as a controlled floating-base system within this
framework is the falling cat \cite%
{MacDonald1961,KaneScher1969,Essen1981,Enos1993,Montgomery1993,GerritsenKuipers1979}%
.

As an example, consider a human body (or humanoid robot) model. The spatial
hip motion is described by $g\in G=SE\left( 3\right) $, and the motion of
body segments 
%TCIMACRO{\TeXButton{red}{\color[rgb]{0,0,0}}}%
%BeginExpansion
\color[rgb]{0,0,0}%
%EndExpansion
(body shape) relative to the hip 
%TCIMACRO{\TeXButton{black}{\color{black}}}%
%BeginExpansion
\color{black}%
%EndExpansion
is described by the joint variables (angles) $\mathbf{r}\in {\mathbb{V}}^{%
\bar{\delta}}$. The kinetic energy is invariant under left $G$-actions,
while the potential energy is invariant under rotations about the vector of
gravity. Thus, with the obvious choice $G=SE\left( 3\right) $, and with the
Lagrangian defined as kinetic minus potential energy, $Q$ is not a principal
bundle. This situation is referred to as symmetry breaking. However, for
most mechanical systems, a $G$-invariant Lagrangian can be defined, e.g.
restricting it to the kinetic energy and including potential forces
separately. Moreover, for discrete mechanical systems, this bundle is often
trivial, so that the configuration space is $Q=G\times M$, with shape space $%
M$, which will be identified with ${\mathbb{V}}^{\bar{\delta}}$ in the
following, as in Sec. \ref{secSymConstr}. This is naturally so for floating
base multibody systems. It is assumed in the following, that $Q$ is a
principal bundle 
%TCIMACRO{\TeXButton{red}{\color[rgb]{0,0,0}}}%
%BeginExpansion
\color[rgb]{0,0,0}%
%EndExpansion
, and a local trivialization is assumed%
%TCIMACRO{\TeXButton{black}{\color{black}}}%
%BeginExpansion
\color{black}%
%EndExpansion
. In a local bundle trivialization, the Lagrangian is then written as $%
l(r^{I},\xi ^{\alpha },\dot{r}^{I})=T(r^{I},\xi ^{\alpha },\dot{r}%
^{I})-V(r^{I})$ (Note that it is still possible to introduce a potential
depending on shape variables $r^{I}$, e.g. accounting for compliant
elements).

\subsection{Boltzmann-Hamel Equations 
%TCIMACRO{\TeXButton{red}{\color[rgb]{0,0,0}}}%
%BeginExpansion
\color[rgb]{0,0,0}%
%EndExpansion
Reduced Euler--Lagrange Equations}

The Boltzmann-Hamel equations are expressed in terms of the velocity
coordinates $\left( u^{a}\right) =(\xi ^{\alpha },\dot{r}^{I})$, where $\xi
^{\alpha }=A_{\beta }^{\alpha }\left( s^{\alpha }\right) \dot{s}^{\beta }$
are expressed in terms of canonical coordinates $s^{\beta }$ on $G$ (as in
Sec. 
%TCIMACRO{\TeXButton{red}{\color[rgb]{0,0,0}}}%
%BeginExpansion
\color[rgb]{0,0,0}%
%EndExpansion
\ref{secSymConstr}\ref{secHamelCoeff})%
%TCIMACRO{\TeXButton{black}{\color{black}}}%
%BeginExpansion
\color{black}%
%EndExpansion
. The Hamel coefficients $\gamma _{bc}^{I}\equiv 0$ follow immediately
noting that $\dot{r}^{I}$ are integrable, and $\gamma _{bc}^{\alpha }\equiv
0 $ follows from (\ref{gamma}) noting that $\xi ^{\alpha }$ are independent
of $r^{I}$. The remaining coefficients $\gamma _{\beta \delta }^{\alpha
}=\pm c_{\beta \delta }^{\alpha }$ are the structure constants of $G$ (lemma %
\ref{lemKiCon}). The Boltzmann-Hamel equations (\ref{BHEqu2}) to the
Lagrangian $l(r^{I},\xi ^{\alpha },\dot{r}^{I})$ are thus the well-known
equations \cite{MarsdenScheurle1993,MarsdenBook1995}%
\begin{align}
\frac{d}{dt}\frac{\partial l}{\partial \xi ^{\alpha }}\pm \frac{\partial l}{%
\partial \xi ^{\beta }}c_{\alpha \lambda }^{\beta }\xi ^{\lambda }& =0
\label{HEorig1} \\
\frac{d}{dt}\frac{\partial l}{\partial \dot{r}^{I}}-\frac{\partial l}{%
\partial r^{I}}& =0  \label{HEorig2}
\end{align}%
where in (\ref{HEorig1}) the positive sign holds for a left-%
%TCIMACRO{\TeXButton{red}{\color[rgb]{0,0,0}}}%
%BeginExpansion
\color[rgb]{0,0,0}%
%EndExpansion
trivialization 
%TCIMACRO{\TeXButton{black}{\color{black}}}%
%BeginExpansion
\color{black}%
%EndExpansion
($\hat{\xi}=g^{-1}\dot{g}$ body velocities), and the negative sign when
right-trivialization is used ($\hat{\xi}=\dot{g}g^{-1}$ spatial velocities).
The base motion is reconstructed with the respective kinematic equation in (%
\ref{KinRec}).

The above equations obviously split into the Euler-Poincar\'{e} equations (%
\ref{HEorig1}) and Euler-Lagrange equations (\ref{HEorig2}), where the first
equation (\ref{HEorig1}) can also be written as Lie-Poisson equations on $G$
when expressed with momentum $\Pi _{\alpha }=\frac{\partial l}{\partial \xi
^{\alpha }}$. While here they have been solely derived with the Hamel
formalism, they were derived as reduced Euler-Lagrange equations \cite%
{MarsdenScheurle1993,MarsdenBook1995} from a variational principle on $Q$.
The involved variations are not intrinsic in the sense that they are not
split into variations in the fiber and the base manifold, respectively. Such
a splitting leads to a decoupling of the equations and to a
block-diagonalization of the mass matrix defining the kinetic energy. This
is achieved by variations of $\xi ^{\alpha }$ with zero variations of $\dot{r%
}^{I}$, i.e. setting the shape velocity to zero, which in geometric terms is
equivalent to variations in the vertical space of the principle bundle. This
is formalized using a connection, as described next.

\subsection{The Mechanical Connection 
%TCIMACRO{\TeXButton{red}{\color[rgb]{0,0,0}}}%
%BeginExpansion
\color[rgb]{0,0,0}%
%EndExpansion
and Locked Velocity%
%TCIMACRO{\TeXButton{black}{\color{black}}}%
%BeginExpansion
\color{black}%
%EndExpansion
%TCIMACRO{\TeXButton{secMechConn}{\label{secMechConn}}}%
%BeginExpansion
\label{secMechConn}%
%EndExpansion
}

Introducing a connection is not as obvious as in case of kinematic
constraint. The starting point is the reduced Lagrangian $l(r^{I},\xi
^{\alpha },\dot{r}^{I})=T(r^{I},\xi ^{\alpha },\dot{r}^{I})-V(r^{I})$ with
potential energy $V(r^{I})$, and kinetic energy%
\begin{equation}
T(r^{I},\xi ^{\alpha },\dot{r}^{I})=\frac{1}{2}\left( 
%TCIMACRO{\TeXButton{xi}{\bm{\xi}}}%
%BeginExpansion
\bm{\xi}%
%EndExpansion
^{T}\ \ \dot{\mathbf{r}}^{T}\right) \left( 
\begin{array}{cc}
\mathbf{L}\left( \mathbf{r}\right) & \mathbf{K}\left( \mathbf{r}\right) \\ 
\mathbf{K}^{T}\left( \mathbf{r}\right) & \mathbf{S}\left( \mathbf{r}\right)%
\end{array}%
\right) \left( 
\begin{array}{c}
%TCIMACRO{\TeXButton{xi}{\bm{\xi}} }%
%BeginExpansion
\bm{\xi}
%EndExpansion
\\ 
\dot{\mathbf{r}}%
\end{array}%
\right)  \label{T}
\end{equation}%
expressed in terms of the mass matrix%
\begin{equation}
\mathbf{M}\left( \mathbf{r}\right) =\left( 
\begin{array}{cc}
\mathbf{L}\left( \mathbf{r}\right) & \mathbf{K}\left( \mathbf{r}\right) \\ 
\mathbf{K}^{T}\left( \mathbf{r}\right) & \mathbf{S}\left( \mathbf{r}\right)%
\end{array}%
\right)  \label{M}
\end{equation}%
which defines a $G$-invariant metric on $Q$. The momentum $\Pi _{\alpha }=%
\frac{\partial l}{\partial \xi ^{\alpha }}$ associated to the fiber is thus $%
%TCIMACRO{\TeXButton{Pi}{\bm{\Pi}}}%
%BeginExpansion
\bm{\Pi}%
%EndExpansion
\left( \mathbf{r},%
%TCIMACRO{\TeXButton{xi}{\bm{\xi}}}%
%BeginExpansion
\bm{\xi}%
%EndExpansion
,\dot{\mathbf{r}}\right) =\mathbf{L}\left( \mathbf{r}\right) 
%TCIMACRO{\TeXButton{xi}{\bm{\xi}}}%
%BeginExpansion
\bm{\xi}%
%EndExpansion
+\mathbf{K}\left( \mathbf{r}\right) \dot{\mathbf{r}}$. 
%TCIMACRO{\TeXButton{red}{\color[rgb]{0,0,0}}}%
%BeginExpansion
\color[rgb]{0,0,0}%
%EndExpansion
If this relation of $%
%TCIMACRO{\TeXButton{xi}{\bm{\xi}}}%
%BeginExpansion
\bm{\xi}%
%EndExpansion
$ and $\dot{\mathbf{r}}$ cannot be integrated to define a relation of $g$
and $\mathbf{r}$, the momentum is said to be non-holonomic.%
%TCIMACRO{\TeXButton{black}{\color{black}}}%
%BeginExpansion
\color{black}%
%EndExpansion

For a floating multibody system, for example, where $\xi ^{\alpha }$ are
twist coordinates of the floating base body, the momentum co-screw $%
%TCIMACRO{\TeXButton{Pi}{\bm{\Pi}}}%
%BeginExpansion
\bm{\Pi}%
%EndExpansion
\in \mathfrak{g}^{\ast }=se^{\ast }\left( 3\right) $ comprises the angular
and linear momentum. One can then introduce the \emph{locked velocity} $%
%TCIMACRO{\TeXButton{Omega}{\bm{\Omega}}}%
%BeginExpansion
\bm{\Omega}%
%EndExpansion
=\left( \Omega ^{\alpha }\right) \in \mathfrak{g}$ such that $\mathbf{L}%
\left( \mathbf{r}\right) 
%TCIMACRO{\TeXButton{Omega}{\bm{\Omega}}}%
%BeginExpansion
\bm{\Omega}%
%EndExpansion
=\mathbf{L}\left( \mathbf{r}\right) 
%TCIMACRO{\TeXButton{xi}{\bm{\xi}}}%
%BeginExpansion
\bm{\xi}%
%EndExpansion
+\mathbf{K}\left( \mathbf{r}\right) \dot{\mathbf{r}}$. The name stems from
the observation that $%
%TCIMACRO{\TeXButton{Omega}{\bm{\Omega}}}%
%BeginExpansion
\bm{\Omega}%
%EndExpansion
$ is the base body velocity which generates the same momentum $%
%TCIMACRO{\TeXButton{Pi}{\bm{\Pi}}}%
%BeginExpansion
\bm{\Pi}%
%EndExpansion
$ when the system is regard as a rigid body, i.e. when $\dot{\mathbf{r}}=0$.
In the humanoid example, this is the velocity the hip would attain, when all
joints are locked instantaneously. Accordingly, $\mathbf{L}\left( \mathbf{r}%
\right) :\mathfrak{g\rightarrow g}^{\ast }$ is called the \emph{locked
inertia tensor}, while $\mathbf{S}$ is the inertia related to the shape
coordinates, and $\mathbf{K}$ is the cross coupling inertia. If $\xi
^{\alpha }$ are coordinates of the base twist in body-fixed representation,
i.e. left-invariant, then $%
%TCIMACRO{\TeXButton{Omega}{\bm{\Omega}}}%
%BeginExpansion
\bm{\Omega}%
%EndExpansion
\in \mathfrak{g}$ is usually called the locked body velocity.

This change of coordinates on $\mathfrak{g}$ is formalized by means of a
connection $\mathcal{A}^{\mathrm{mech}}=Ad_{g}(g^{-1}dg+\mathcal{A}d\mathbf{r%
})$ on the principal bundle $Q$. The local connection one-form is defined
with $\left( \mathcal{A}_{I}^{\alpha }\right) :=\mathbf{L}^{-1}\mathbf{K}$
so that the velocity shift is%
\begin{equation}
\Omega ^{\alpha }=\xi ^{\alpha }+\mathcal{A}_{I}^{\alpha }(r^{I})\dot{r}^{I}.
\label{mechConn}
\end{equation}%
The so-defined connection, is called the \emph{mechanical connection} \cite%
{Marsden1992,MarsdenScheurle1993} building upon a concept discussed in \cite%
{Smale1970}. In contrast to the kinematic connection, it is defined via the
momentum. The connection may be considered to be in Chaplygin form since it
is independent of group variables. Notice that this strictly relies on a
local bundle trivialization since $Q$ may not be a trivial principal bundle.
The locked velocity $%
%TCIMACRO{\TeXButton{Omega}{\bm{\Omega}}}%
%BeginExpansion
\bm{\Omega}%
%EndExpansion
$ is the vertical part relative to the mechanical connection.

\begin{remark}
%TCIMACRO{\TeXButton{remMomentumInt2}{\label{remMomentumInt2}}}%
%BeginExpansion
\label{remMomentumInt2}%
%EndExpansion
%TCIMACRO{\TeXButton{red}{\color[rgb]{0,0,0}}}%
%BeginExpansion
\color[rgb]{0,0,0}%
%EndExpansion
An important aspect of the locked velocity is that it cannot be associated
to a frame whose motion is described by $h\left( g,\mathbf{r}\right) \in G$
depending on some $g$ and $\mathbf{r}$, so that $\hat{%
%TCIMACRO{\TeXButton{Omega}{\bm{\Omega}}}%
%BeginExpansion
\bm{\Omega}%
%EndExpansion
}=h^{-1}\dot{h}$. This is an immediate consequence of the fact that $\Omega $
is defined by the non-holonomic momentum $%
%TCIMACRO{\TeXButton{Pi}{\bm{\Pi}}}%
%BeginExpansion
\bm{\Pi}%
%EndExpansion
\left( \mathbf{r},%
%TCIMACRO{\TeXButton{xi}{\bm{\xi}}}%
%BeginExpansion
\bm{\xi}%
%EndExpansion
,\dot{\mathbf{r}}\right) $ (see Sec. \ref{secUnConBundle}\ref{secDecoupling}%
).%
%TCIMACRO{\TeXButton{black}{\color{black}}}%
%BeginExpansion
\color{black}%
%EndExpansion
\end{remark}

\begin{remark}
The principal bundle view on floating-base systems has an interesting
connection to gauge theory. In gauge theory, the symmetry is related to
gauge invariance, $G$ is called the 'gauge group', and the connection
one-form $\mathcal{A}$ in (\ref{mechConn}) is the 'gauge potential' \cite%
{Choquet1989}. This was discussed for the falling cat and similar
non-holonomic control systems in \cite{Montgomery1993}, and more generally,
for 'deformable bodies' (mechanical structures that can change their shape)
in \cite{WilczekShapere1987,ShapereWilczek1989}, and for Maxwell or
Yang-Mills fields in \cite{MontgomeryMarsdenRatiu1984}. In gauge theory, the
equivariance condition on the connection describes a gauge transformation
from a local gauge $\mathcal{A}$ to a new $\mathcal{A}^{\prime
}=Ad_{g}(g^{-1}dg+\mathcal{A}d\mathbf{r})$ \cite{BaezMuniain1994,Isham1999}.
In case of mechanical systems, it describes a transformation from one
body-fixed frame to another in which velocities are measured. In \cite%
{WilczekShapere1987,ShapereWilczek1989}, $\mathcal{A}_{I}^{\beta }$ was
called the master gauge, while its curvature is considered as field strength.
\end{remark}

\subsection{Hamel Coefficients and the 
%TCIMACRO{\TeXButton{red}{\color[rgb]{0,0,0}}}%
%BeginExpansion
\color[rgb]{0,0,0}%
%EndExpansion
Mechanical Connection%
%TCIMACRO{\TeXButton{black}{\color{black}}}%
%BeginExpansion
\color{black}%
%EndExpansion
}

Here again, the original definition (\ref{gamma}) of the Hamel coefficients
is employed to derive relations that are 
%TCIMACRO{\TeXButton{red}{\color[rgb]{0,0,0}}}%
%BeginExpansion
\color[rgb]{0,0,0}%
%EndExpansion
today 
%TCIMACRO{\TeXButton{black}{\color{black}}}%
%BeginExpansion
\color{black}%
%EndExpansion
obtained using the machinery of geometric mechanics.

\begin{lemma}
The Hamel coefficients in bundle coordinates $(\Omega ^{\alpha },\dot{r}^{I})
$ are $\gamma _{ab}^{I}\equiv 0$, and%
\begin{align}
\gamma _{\beta \delta }^{\alpha }& =\pm c_{\beta \delta }^{\alpha }
\label{gamma41} \\
\gamma _{J\beta }^{\alpha }& =-\gamma _{\beta J}^{\alpha }=\pm c_{\beta
\delta }^{\alpha }\mathcal{A}_{J}^{\delta }  \label{gamma42} \\
\gamma _{IJ}^{\alpha }& =\frac{\partial \mathcal{A}_{I}^{\alpha }}{\partial
r^{J}}-\frac{\partial \mathcal{A}_{J}^{\alpha }}{\partial r^{I}}+\gamma
_{J\beta }^{\alpha }\mathcal{A}_{I}^{\beta }=\frac{\partial \mathcal{A}%
_{I}^{\alpha }}{\partial r^{J}}-\frac{\partial \mathcal{A}_{J}^{\alpha }}{%
\partial r^{I}}\pm c_{\beta \lambda }^{\alpha }\mathcal{A}_{I}^{\beta }%
\mathcal{A}_{J}^{\lambda }  \label{gamma43}
\end{align}%
where the positive sign of $\pm $ applies to left-, and the negative sign to
a right-trivialization of the $G$-bundle.
\end{lemma}

\begin{proof}
Canonical coordinates $s^{\alpha }$ are introduced on $G$, which are related
to the fiber coordinates via $\xi ^{\alpha }=A_{\beta }^{\alpha }\left(
s^{\alpha }\right) \dot{s}^{\beta }$ and $\dot{s}^{\alpha }=B_{\beta
}^{\alpha }\left( s^{\alpha }\right) \xi ^{\beta }$, respectively (see Sec. %
\ref{secSymConstr}\ref{secHamelCoeff}). Relation (\ref{mechConn}) and its
inverse are then written as%
\begin{align}
\Omega ^{\alpha }& =A_{\beta }^{\alpha }\left( s^{\alpha }\right) \dot{s}%
^{\beta }+\mathcal{A}_{I}^{\alpha }(r^{J})\dot{r}^{I} \\
\dot{s}^{\alpha }& =B_{\beta }^{\alpha }\left( s^{a}\right) \Omega ^{\beta
}-B_{\beta }^{\alpha }\left( s^{a}\right) \mathcal{A}_{I}^{\beta
}(r^{J})u^{I}
\end{align}%
With the locked velocity, the quasi-velocities are $\left( u^{a}\right)
=(u^{\alpha },u^{I})=(\Omega ^{\alpha },\dot{r}^{I})$, thus expressed in the
form (\ref{qdu}), with $B_{J}^{I}=\delta _{J}^{I}$ and $B_{\beta }^{I}\equiv
0$, or in matrix form%
\begin{equation}
\mathbf{u}=\left( 
\begin{array}{cc}
\left( A_{\beta }^{\alpha }\right) & \left( \mathcal{A}_{I}^{\alpha }\right)
\\ 
\mathbf{0} & \mathbf{I}%
\end{array}%
\right) \left( 
\begin{array}{c}
\dot{\mathbf{s}} \\ 
\dot{\mathbf{r}}%
\end{array}%
\right) ,\ \ \ \left( 
\begin{array}{c}
\dot{\mathbf{s}} \\ 
\dot{\mathbf{r}}%
\end{array}%
\right) =\left( 
\begin{array}{cc}
\left( B_{\beta }^{\alpha }\right) & \ \ -\left( B_{\beta }^{\alpha }%
\mathcal{A}_{I}^{\beta }\right) \\ 
\mathbf{0} & \mathbf{I}%
\end{array}%
\right) \mathbf{u}.
\end{equation}%
The relation (\ref{gamma}) is separated for the coordinates $\Omega ^{\alpha
}$ and $\dot{r}^{I}$. The coefficients $\gamma _{\beta \delta }^{\alpha
}=\left( \frac{\partial A_{\mu }^{\alpha }}{\partial s^{\lambda }}-\frac{%
\partial A_{\lambda }^{\alpha }}{\partial s^{\mu }}\right) B_{\beta }^{\mu
}B_{\delta }^{\lambda }=\pm c_{\beta \delta }^{\alpha }$ are again
determined by the structure constants of $G$, and thus%
\begin{equation*}
\gamma _{\beta J}^{\alpha }=\left( \frac{\partial A_{\delta }^{\alpha }}{%
\partial s^{\lambda }}-\frac{\partial A_{\lambda }^{\alpha }}{\partial
s^{\delta }}\right) B_{\beta }^{\lambda }B_{\gamma }^{\delta }\mathcal{A}%
_{J}^{\gamma }=\pm c_{\beta \delta }^{\alpha }\mathcal{A}_{J}^{\delta }.
\end{equation*}%
The remaining coefficients $\gamma _{IJ}^{\alpha }$ are given in (\ref%
{gamma3}).
\end{proof}

%TCIMACRO{\TeXButton{red}{\color[rgb]{0,0,0}}}%
%BeginExpansion
\color[rgb]{0,0,0}%
%EndExpansion
The Hamel coefficients (\ref{gamma43}) are identical to the components of
the curvature (\ref{curvature}), $\mathcal{B}_{IJ}^{\alpha }=\gamma
_{IJ}^{\alpha }$, of the mechanical connection in bundle coordinates. 
%TCIMACRO{\TeXButton{black}{\color{black}}}%
%BeginExpansion
\color{black}%
%EndExpansion
They are indeed formally identical to the curvature coefficients (\ref%
{gamma3}) of the kinematic connection.

\subsection{Boltzmann-Hamel Equations as reduced Euler-Lagrange Equations in
Bundle Coordinates%
%TCIMACRO{\TeXButton{secBHBudle}{\label{secBHBudle}}}%
%BeginExpansion
\label{secBHBudle}%
%EndExpansion
}

Denote with $\ell (r^{I},\Omega ^{\alpha },\dot{r}^{I}):=l(r^{I},\xi
^{\alpha }:=\Omega ^{\alpha }-\mathcal{A}_{I}^{\alpha }\dot{r}^{I},\dot{r}%
^{I})$ the reduced Lagrangian in terms of the locked velocity $\Omega
^{\alpha }$. The Hamel equations (\ref{BHEqu2}) are%
\begin{eqnarray}
\frac{d}{dt}\frac{\partial \ell }{\partial \Omega ^{\alpha }}+\frac{\partial
\ell }{\partial \Omega ^{\beta }}\left( \gamma _{\alpha I}^{\beta }\dot{r}%
^{I}+\gamma _{\alpha \lambda }^{\beta }\Omega ^{\lambda }\right)
&=&Q_{\alpha }  \label{BHESym1} \\
\frac{d}{dt}\frac{\partial \ell }{\partial \dot{r}^{I}}-\frac{\partial \ell 
}{\partial r^{I}}+\frac{\partial \ell }{\partial \Omega ^{\beta }}\left(
\gamma _{IJ}^{\beta }\dot{r}^{J}+\gamma _{I\alpha }^{\beta }\Omega ^{\alpha
}\right) &=&Q_{I}.  \label{BHESym2}
\end{eqnarray}%
The motion in $G$ is obtained from the reconstruction equations (\ref{KinRec}%
), now with $\xi ^{\alpha }=\Omega ^{\alpha }-\mathcal{A}_{I}^{\alpha
}\left( \mathbf{r}\right) \dot{r}^{I}$ defined by the mechanical connection.
It is important that the motion is deduced from $\xi ^{\alpha }$, and not
from $\Omega ^{\alpha }$, as the latter can in general not be attributed to
a frame motion (Rem. \ref{remMomentumInt2} and Sec. \ref{secUnConBundle}\ref%
{secDecoupling}). Equations (\ref{BHESym1},\ref{BHESym2}), admit the
following geometric interpretation.

\begin{proposition}
The Hamel equations (\ref{BHESym1},\ref{BHESym2}) for a left $G$-invariant
Lagrangian $\ell (r^{I},\Omega ^{\alpha },\dot{r}^{I})$ are the reduced
Euler-Lagrange equations (\ref{LP1},\ref{LP2}) 
%TCIMACRO{\TeXButton{red}{\color[rgb]{0,0,0}}}%
%BeginExpansion
\color[rgb]{0,0,0}%
%EndExpansion
in terms of bundle coordinates $(\Omega ^{\alpha },\dot{r}^{I})$, with the
coefficients $\mathcal{E}_{\beta I}^{\alpha }=\gamma _{I\beta }^{\alpha
}=c_{\beta \lambda }^{\alpha }\mathcal{A}_{I}^{\lambda },c_{\beta \lambda
}^{\alpha }=\gamma _{\beta \lambda }^{\alpha }$, and curvature $\mathcal{B}%
_{IJ}^{\alpha }=\gamma _{IJ}^{\alpha }$ determined by the Hamel-coefficients.
\end{proposition}

The equations (\ref{LP1},\ref{LP2}) have been derived in \cite%
{MarsdenScheurle1993}, and presented in \cite[p. 397]{MarsdenBook1995}, by
introducing the locked velocity (\ref{mechConn}) into the Lagrangian $\ell
(r^{I},\Omega ^{\alpha },\dot{r}^{I})$ before taking the Euler-Lagrange
derivative. Their derivation as Hamel equation in terms of the locked
velocity has not been reported in the literature. The first equation (\ref%
{LP1}), respectively (\ref{BHESym1}), is indeed the Euler-Poincar\'{e}
equation (\ref{HEorig1}) with $\xi ^{\alpha }$ replaced by the velocity $%
\Omega ^{\alpha }$ of the locked system, which is why (\ref{LP1},\ref{LP2})
are also called Lagrange-Poincar\'{e} equations \cite[p. 3395]%
{MarsdenRatiuScheurle2000},\cite[p. 146]{BlochBook2003}. The terms with $%
\gamma _{\alpha I}^{\beta }$ in (\ref{BHESym1}) and (\ref{BHESym2}) can be
regarded as interaction (or coupling) terms. Clearly, as remarked in \cite[%
pp. 912, 913]{KoonMarsden1997}, the equations (\ref{BHESym1},\ref{BHESym2})
reduce to the Hamel equations (\ref{HEorig1},\ref{HEorig2}) if the
coefficients of connection and curvature vanish, i.e. when expressed in
local coordinates $(\xi ^{\alpha },\dot{r}^{I})$. However, since the Hamel
formalism applies to any choice of local coordinates, as shown in this
paper, it should not be said that (\ref{BHESym1},\ref{BHESym2}) reduce to 
\emph{the} Hamel equations when using local coordinates $(\xi ^{\alpha },%
\dot{r}^{I})$, as occasionally stated, e.g. \cite[p. 226]%
{CendraMarsdenRatiu2001}.

It follows from the definition of the Hamel coefficients that the mechanical
connection is flat, i.e. the curvature $\mathcal{B}_{IJ}^{\alpha }$
vanishes, if and only if the momentum $\Pi _{\alpha }=\frac{\partial \ell }{%
\partial \Omega ^{\alpha }}$ (equivalently $\Pi _{\alpha }=\frac{\partial l}{%
\partial \xi ^{\alpha }}$) 
%TCIMACRO{\TeXButton{red}{\color[rgb]{0,0,0}}}%
%BeginExpansion
\color[rgb]{0,0,0}%
%EndExpansion
defines a non-integrable relation of $\xi ^{\alpha }$ and $\dot{r}^{I}$.

\begin{remark}
%TCIMACRO{\TeXButton{remSign}{\label{remSign}}}%
%BeginExpansion
\label{remSign}%
%EndExpansion
A note on the sign convention for the curvature coefficients $\mathcal{B}%
_{IJ}^{\alpha }$ is in order. Given a connection, the local curvature is
usually defined as $\mathcal{B}_{IJ}^{\alpha }=-\gamma _{IJ}^{\alpha }$ with 
$\gamma _{IJ}^{\alpha }$ in (\ref{gamma3}) (see \cite{Choquet1989} for right
bundles), which agrees with the definition of curvature used in gauge theory 
\cite[p. 247]{BaezMuniain1994}. This convention is used in \cite[p. 157]%
{MarsdenScheurle1993} and \cite[p. 3395]{MarsdenRatiuScheurle2000}, and thus 
$\mathcal{B}_{IJ}^{\alpha }$ appears with a positive sign in the reduced
Euler-Lagrange equations (\ref{LP2}). In \cite[p. 910]{KoonMarsden1997} and 
\cite[notice the correction on p. 44]{BlochKrishnaprasadMarsdenMurray1996},
the local curvature is introduced as $\mathcal{B}_{IJ}^{\alpha }=\gamma
_{IJ}^{\alpha }$, as in this paper, along with the reduced Euler-Lagrange
equations (\ref{LP1},\ref{LP2}). In \cite[pp. 117,146]{BlochBook2003}, the
curvature is derived as $\mathcal{B}_{IJ}^{\alpha }=-\gamma _{IJ}^{\alpha }$%
, but is then used with a negative sign in the reduced Euler-Lagrange
equations. Similarly in \cite{CendraMarsdenRatiu2001}, the curvature is
introduced as $\mathcal{B}_{IJ}^{\alpha }=\gamma _{IJ}^{\alpha }$ but used
with positive sign in (\ref{LP2}). These inconsistencies deserve particular
attention when applying equations (\ref{HE-kinBundle}) and (\ref{LP1},\ref%
{LP2}).
\end{remark}

\subsection{Inertial Decoupling of Poincar\'{e} and Lagrange Equations\label%
{secDecoupling}}

The mechanical connection allows to intrinsically split variation into the
vertical and horizontal parts\footnote{%
%TCIMACRO{\TeXButton{red}{\color[rgb]{0,0,0}}}%
%BeginExpansion
\color[rgb]{0,0,0}%
%EndExpansion
The vertical space is the tangent space $\ker T_{q}\pi $ to the group
orbits, i.e. possible velocities of the base body for locked shape
coordinates. The horizontal space is the space of velocities not producing a
locked velocity.}. As a consequence, the mass matrix $\mathbf{M}$ of the
equations in terms of the locked velocity, is diagonal \cite[p. 147]%
{BlochBook2003}. Inserting (\ref{mechConn}) into (\ref{T}) yields the mass
matrix%
\begin{equation}
\mathbf{M}^{\Omega }\left( \mathbf{r}\right) =\left( 
\begin{array}{cc}
\mathbf{L}\left( \mathbf{r}\right) & \mathbf{0} \\ 
\mathbf{0} & \mathbf{S}\left( \mathbf{r}\right) -\mathcal{A}^{T}\left( 
\mathbf{r}\right) \mathbf{L}\left( \mathbf{r}\right) \mathcal{A}\left( 
\mathbf{r}\right)%
\end{array}%
\right)  \label{MOmega}
\end{equation}%
and the kinetic energy $\ell (r^{I},\Omega ^{\alpha },\dot{r}^{I})=\frac{1}{2%
}%
%TCIMACRO{\TeXButton{Omega}{\bm{\Omega}}}%
%BeginExpansion
\bm{\Omega}%
%EndExpansion
^{T}\mathbf{L}\left( \mathbf{r}\right) 
%TCIMACRO{\TeXButton{Omega}{\bm{\Omega}}}%
%BeginExpansion
\bm{\Omega}%
%EndExpansion
+\frac{1}{2}\dot{\mathbf{r}}^{T}\left( \mathbf{S}\left( \mathbf{r}\right) -%
\mathcal{A}^{T}\left( \mathbf{r}\right) \mathbf{L}\left( \mathbf{r}\right) 
\mathcal{A}\left( \mathbf{r}\right) \right) \dot{\mathbf{r}}$. Indeed, the
momentum $\Pi _{\alpha }=\frac{\partial \ell }{\partial \Omega ^{\alpha }}$
only depends on the fiber coordinates $\Omega ^{\alpha }$, so that the
equations (\ref{BHESym1}) and (\ref{BHESym2}) are inertially decoupled (not
coupled on accelerations level). Coupling of the equations is via the
velocity terms involving $\gamma _{\alpha I}^{\beta }$. Inertial decoupling
using the locked velocity has been addressed for modeling of floating base
robots \cite{Garofalo2015,Mishra2020} and space robots in \cite%
{Giordano2018,GiordanoRAL2019}.

A closely related concept for decoupling the equations is that of the \emph{%
centroidal momentum} as introduced in \cite{Essen1993}, which is widely used
for whole-body control of humanoid robots \cite%
{OrinGoswamiLee2013,Nava2018,LeeGoswami2007} for instance. In this context $%
\mathfrak{g}=se\left( 3\right) $, and $\mathbf{V}\in se\left( 3\right) $ is
the velocity (twist) of the base body (using symbol $\mathbf{V}$ instead of $%
\xi $), and $\dot{\mathbf{r}}$ are the joint velocities. A frame $\mathcal{F}%
_{\mathrm{G}}$ is introduced that is located at the total COM of the system
and aligned with the inertial frame $\mathcal{F}_{0}$. The configuration of $%
\mathcal{F}_{\mathrm{G}}$ relative to the frame $\mathcal{F}_{\mathrm{b}}$
attached at the base body is described by $g_{\mathrm{bG}}\in SE\left(
3\right) =G$. The centroidal momentum is defined as $%
%TCIMACRO{\TeXButton{Pi}{\bm{\Pi}}}%
%BeginExpansion
\bm{\Pi}%
%EndExpansion
_{\mathrm{G}}=\mathbf{Ad}_{g_{\mathrm{bG}}}^{T}%
%TCIMACRO{\TeXButton{Pi}{\bm{\Pi}}}%
%BeginExpansion
\bm{\Pi}%
%EndExpansion
$, with momentum co-screw $%
%TCIMACRO{\TeXButton{Pi}{\bm{\Pi}}}%
%BeginExpansion
\bm{\Pi}%
%EndExpansion
\in se^{\ast }\left( 3\right) $ defined in Sec. \ref{secUnConBundle}\ref%
{secMechConn}. This is also expressed as $%
%TCIMACRO{\TeXButton{Pi}{\bm{\Pi}}}%
%BeginExpansion
\bm{\Pi}%
%EndExpansion
_{\mathrm{G}}=\mathbf{M}_{\mathrm{G}}\mathbf{V}_{\mathrm{G}}$, where $%
\mathbf{V}_{\mathrm{G}}$ is referred to as the \emph{average velocity}, and $%
\mathbf{M}_{\mathrm{G}}=\mathbf{Ad}_{g_{\mathrm{bG}}}^{T}\mathbf{M}_{\mathrm{%
bb}}\mathbf{Ad}_{g_{\mathrm{bG}}}$ is called the \emph{centroidal composite
inertia matrix}, and $\mathbf{Ad}_{g_{\mathrm{bG}}}^{T}(\mathbf{L},\mathbf{K}%
)$ the \emph{centroidal momentum matrix} \cite%
{OrinGoswami2008,OrinGoswamiLee2013}. Comparing this with the definition of
the locked velocity $\mathbf{L}\left( \mathbf{r}\right) \mathbf{V}_{\mathrm{%
loc}}=%
%TCIMACRO{\TeXButton{Pi}{\bm{\Pi}}}%
%BeginExpansion
\bm{\Pi}%
%EndExpansion
$ shows that $\mathbf{V}_{\mathrm{G}}=\mathbf{Ad}_{g_{\mathrm{bG}}}^{-1}%
\mathbf{V}_{\mathrm{loc}}$. Moreover, $\mathbf{V}_{\mathrm{G}}=\left( 
%TCIMACRO{\TeXButton{w}{\mathbold{\omega}}}%
%BeginExpansion
\mathbold{\omega}%
%EndExpansion
_{\mathrm{ave}},\dot{\mathbf{p}}_{\mathrm{com}}\right) $, where $%
%TCIMACRO{\TeXButton{w}{\mathbold{\omega}}}%
%BeginExpansion
\mathbold{\omega}%
%EndExpansion
_{\mathrm{ave}}=%
%TCIMACRO{\TeXButton{w}{\mathbold{\omega}}}%
%BeginExpansion
\mathbold{\omega}%
%EndExpansion
_{\mathrm{loc}}$ is called the \emph{average angular velocity} \cite%
{Essen1993}, and $\dot{\mathbf{p}}_{\mathrm{com}}$ is the velocity of the
total COM of the system, both expressed in $\mathcal{F}_{0}$. The important
point is that $\mathbf{M}_{\mathrm{G}}=\left( 
\begin{array}{cc}
%TCIMACRO{\TeXButton{Theta}{{\bm{\Theta}}}}%
%BeginExpansion
{\bm{\Theta}}%
%EndExpansion
_{\mathrm{G}} & \mathbf{0} \\ 
\mathbf{0} & \bar{m}\mathbf{I}%
\end{array}%
\right) $ is a block diagonal matrix, where $\bar{m}$ is the total mass, and 
$%
%TCIMACRO{\TeXButton{Theta}{{\bm{\Theta}}}}%
%BeginExpansion
{\bm{\Theta}}%
%EndExpansion
_{\mathrm{G}}$ is the total inertia tensor w.r.t. to the total COM expressed
in $\mathcal{F}_{0}$. This would replace the locked inertia $\mathbf{L}$ in (%
\ref{MOmega}) when the EOM are expressed with $\mathbf{V}_{\mathrm{G}}$.

%TCIMACRO{\TeXButton{red}{\color[rgb]{0,0,0}}}%
%BeginExpansion
\color[rgb]{0,0,0}%
%EndExpansion
The motivation for using the centroidal momentum $%
%TCIMACRO{\TeXButton{Pi}{\bm{\Pi}}}%
%BeginExpansion
\bm{\Pi}%
%EndExpansion
_{\mathrm{G}}=\left( 
%TCIMACRO{\TeXButton{Theta}{{\bm{\Theta}}}}%
%BeginExpansion
{\bm{\Theta}}%
%EndExpansion
_{\mathrm{G}}%
%TCIMACRO{\TeXButton{w}{\mathbold{\omega}}}%
%BeginExpansion
\mathbold{\omega}%
%EndExpansion
_{\mathrm{ave}},\bar{m}\dot{\mathbf{p}}_{\mathrm{com}}\right) $ is that the
linear and angular momentum are decoupled (in addition to the inertial
decoupling of (\ref{BHESym1}) and (\ref{BHESym2})), and can be controlled
independently. However, there is generally no frame associated to $\mathbf{V}%
_{\mathrm{G}}$ (Rem. \ref{remMomentumInt2}) that could serve to represent
the system orientation. This would imply that the motion of this frame is
represented by a $g_{\mathrm{bG}}\left( g,\mathbf{r}\right) $ such that $%
\hat{\mathbf{V}}_{\mathrm{G}}=\dot{g}_{\mathrm{bG}}g_{\mathrm{bG}}^{-1}$. It
is clear from the definition of $\mathbf{V}_{\mathrm{G}}$ (and $\mathbf{V}_{%
\mathrm{loc}}$) by means of the momentum that such a frame exists if and
only if the momentum defines integrable relation, i.e. if the curvature
vanishes, which is generally not the case. 
%TCIMACRO{\TeXButton{black}{\color{black}}}%
%BeginExpansion
\color{black}%
%EndExpansion
This seemingly obvious fact was proven in \cite{Saccon2017}. In order to
determine the base configuration $g\in SE\left( 3\right) $ w.r.t. $\mathcal{F%
}_{0}$, the reconstruction equations (\ref{KinRec}), which are now $\hat{%
\mathbf{V}}_{\mathrm{b}}=g^{-1}\dot{g}$, must be solved with $\mathbf{V}_{%
\mathrm{b}}=\mathbf{V}_{\mathrm{loc}}-\mathbf{L}^{-1}\mathbf{K}\dot{\mathbf{r%
}}$, as proposed in \cite{Nava2018}. Finally it should be remarked that the
centroidal kinematics and dynamics can be expressed in terms of barycentric
vectors \cite{PapadopoulosDubowskyTRO1993,Papadopoulos1993}.

\subsection{Example: Satellite with three symmetric reaction wheels\label%
{secSatellite}}

\begin{figure}[bh]
\begin{center}
\includegraphics[draft=false,height=6.5cm]{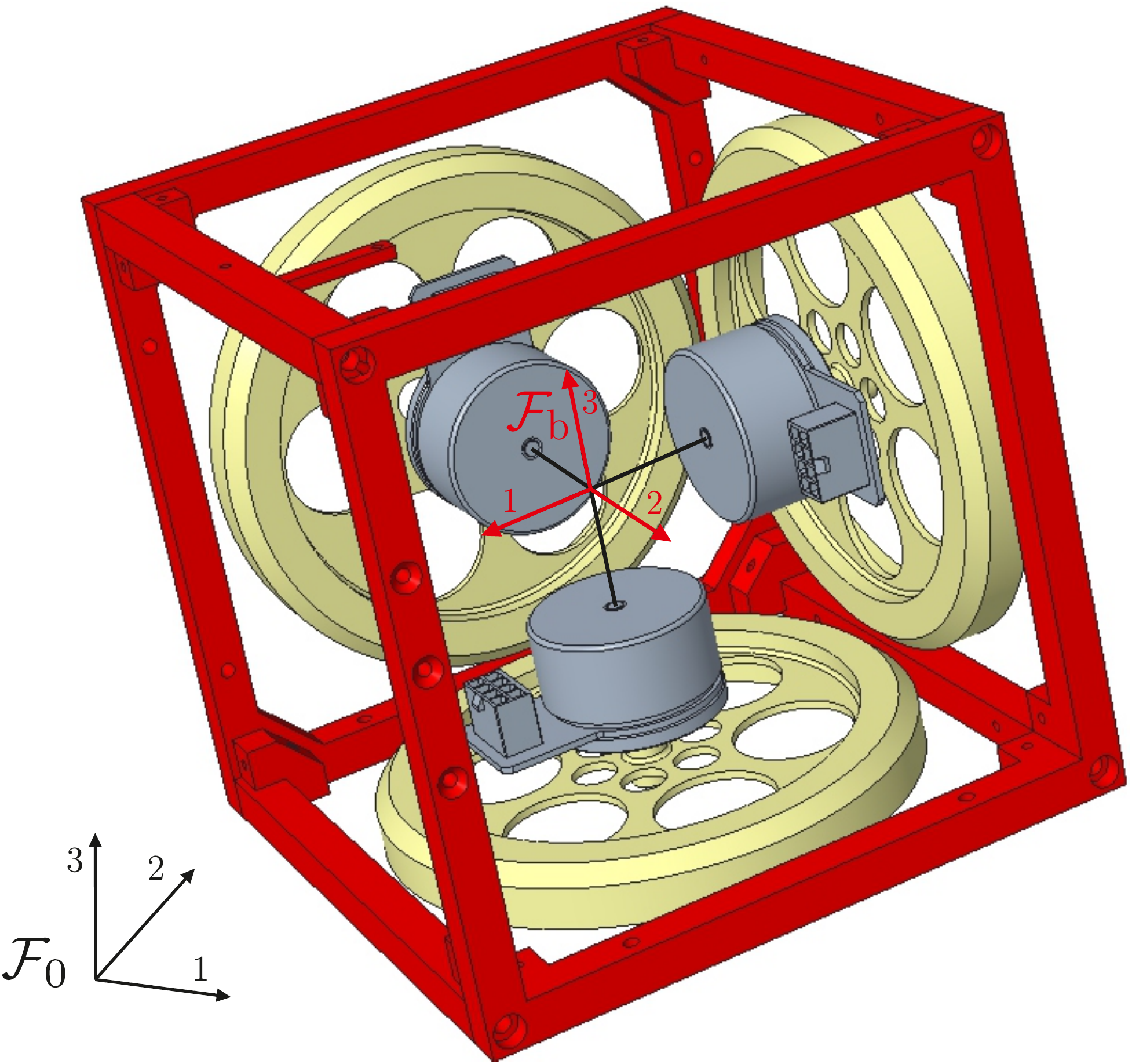}
\end{center}
\caption{Schematic drawing of a satellite model. The figure shows a
self-stabilizing cube reported in \protect\cite{CubeRAAD2016}.}
\label{figCube}
\end{figure}

The simplified model of a satellite equipped with three reaction wheels is
considered. Fig. \ref{figCube} shows a schematic drawing of the principle
mechanical setup. In the following, the reaction wheels are called rotors,
for simplicity. The axes of the three rotors are mutually orthogonal, the
rotors are located arbitrarily at the satellite, and are assumed to be
symmetric (so that the total COM of the satellite is constant). The
satellite's main body is modeled as a rigid body. It is assumed that there
are no gravity or other potential forces acting on the satellite, thus the
Lagrangian is the kinetic energy. The motion of the main body is a rigid
body motion evolving in a Lie group $G$, and represented by $g\in G$. The
kinetic energy is invariant w.r.t. $G$-actions. The rotations of the three
rotors are described by the rotation angles $\varphi ^{i},i=1,2,3$. The
configuration space of the satellite model is thus $Q=G\times T^{3}$, with
configuration $q=\left( g,%
%TCIMACRO{\TeXButton{phi}{\mathbold{\varphi}}}%
%BeginExpansion
\mathbold{\varphi}%
%EndExpansion
\right) $, with $%
%TCIMACRO{\TeXButton{phi}{\mathbold{\varphi}}}%
%BeginExpansion
\mathbold{\varphi}%
%EndExpansion
=(\varphi ^{1},\varphi ^{2},\varphi ^{3})$. The latter serve as coordinates
on the shape space, $r^{I}:=\varphi ^{\bar{I}}$. Different choices for $G$
are used in the literature. Most of the original formulations in multibody
system dynamics use the direct product group $G=SO\left( 3\right) \times {%
\mathbb{R}}^{3}$, while recent research uses the proper rigid body motion
group $SE\left( 3\right) $ (mainly triggered by development of Lie group
integration methods \cite{CelledoniOwren2003}, and geometrically exact
modeling of Cosserat continua \cite%
{BorriBottasso1994a,SonnevilleCardonaBruls2014}). The particular choice of
symmetry group, but also whether left- or right trivialization is used,
leads to different definition of rigid body velocities and equations of
motion. In the following, notation from multibody dynamics and robotics is
adopted, where $\mathbf{V}\in {\mathbb{R}}^{6}\cong \mathfrak{g}$ denotes
velocity of a frame (rigid body), and $\mathfrak{g}$ is either $so\left(
3\right) \times {\mathbb{R}}^{3}$ or $se\left( 3\right) $. 
%TCIMACRO{\TeXButton{red}{\color[rgb]{0,0,0}}}%
%BeginExpansion
\color[rgb]{0,0,0}%
%EndExpansion
A detailed description and numerical results can be found in the supplement 
\cite{Supplement}.%
%TCIMACRO{\TeXButton{black}{\color{black}}}%
%BeginExpansion
\color{black}%
%EndExpansion

\subsubsection{Mixed Representation of Rigid Body Velocity --Symmetry Group $%
G=SO\left( 3\right) \times {\mathbb{R}}^{3}$}

The configuration (pose) of the main body is represented as $\left( \mathbf{R%
},\mathbf{p}\right) \in SO\left( 3\right) \times {\mathbb{R}}^{3}$, where $%
\mathbf{R}\in SO\left( 3\right) $ and $\mathbf{p}\in {\mathbb{R}}^{3}$
describes the rotation and translation of a body-fixed reference frame (RFR) 
$\mathcal{F}_{\mathrm{b}}$ relative to an inertial frame (IFR) $\mathcal{F}%
_{0}$. The group multiplication on the direct product group is $\left( 
\mathbf{R}_{1},\mathbf{p}_{1}\right) \cdot \left( \mathbf{R}_{2},\mathbf{p}%
_{2}\right) =\left( \mathbf{R}_{1}\mathbf{R}_{2},\mathbf{p}_{1}+\mathbf{p}%
_{2}\right) $. Since rotations and translations are decoupled, this is
clearly not a frame transformation (i.e. a rigid body motion). The
corresponding velocity defined via left-trivialization is $\hat{\mathbf{V}}_{%
\mathrm{b}}=g^{-1}\dot{g}=(\mathbf{R}^{-1}\dot{\mathbf{R}},\dot{\mathbf{p}}%
)=\left( \hat{%
%TCIMACRO{\TeXButton{w}{\mathbold{\omega}}}%
%BeginExpansion
\mathbold{\omega}%
%EndExpansion
},\dot{\mathbf{p}}\right) \in \mathfrak{g}=so\left( 3\right) \times {\mathbb{%
R}}^{3}$, and in vector representation $\mathbf{V}_{\mathrm{b}}=\left( 
%TCIMACRO{\TeXButton{w}{\mathbold{\omega}}}%
%BeginExpansion
\mathbold{\omega}%
%EndExpansion
,\dot{\mathbf{p}}\right) \in {\mathbb{R}}^{6}$, where $%
%TCIMACRO{\TeXButton{w}{\mathbold{\omega}}}%
%BeginExpansion
\mathbold{\omega}%
%EndExpansion
\in {\mathbb{R}}^{3}$ is the angular velocity of the main body relative to $%
\mathcal{F}_{0}$ resolved in $\mathcal{F}_{\mathrm{b}}$. This is referred to
as the mixed representation of rigid body velocities since $%
%TCIMACRO{\TeXButton{w}{\mathbold{\omega}}}%
%BeginExpansion
\mathbold{\omega}%
%EndExpansion
$ is resolved in the body frame, and $\dot{\mathbf{p}}$ in the inertia frame 
\cite{Murray,MUBOScrews1}. Regarding the dynamics, decoupling of rotation
and translation is valid only if the body-fixed RFR is located at the COM,
which is the main premise when using the direct product group $G$, since
then the angular and linear momenta are decoupled. Therefore, the velocities
of main body and rotors will be measured at the total COM of the satellite
(main body including the rotors), thus $\mathbf{p}$ is the position vector
of the total COM resolved in the IFR, and $\mathcal{F}_{\mathrm{b}}$ is
located at the total COM. Denote with $\mathbf{V}_{i}=\left( 
%TCIMACRO{\TeXButton{w}{\mathbold{\omega}}}%
%BeginExpansion
\mathbold{\omega}%
%EndExpansion
_{i},\dot{\mathbf{p}}\right) $ the hybrid velocity of rotor $i=1,2,3$.
W.l.o.g. the RFR is aligned with the rotor axes. Then $%
%TCIMACRO{\TeXButton{w}{\mathbold{\omega}}}%
%BeginExpansion
\mathbold{\omega}%
%EndExpansion
_{i}=%
%TCIMACRO{\TeXButton{w}{\mathbold{\omega}}}%
%BeginExpansion
\mathbold{\omega}%
%EndExpansion
+\mathbf{e}_{i}\dot{\varphi}^{i}$, where $\mathbf{e}_{i}\in {\mathbb{R}}^{3}$
is the $i$-th unit vector (e.g. $\mathbf{e}_{1}=\left( 1,0,0\right) $), and
thus $\mathbf{V}_{i}=\mathbf{V}_{\mathrm{b}}+\overline{\mathbf{V}}_{i}$,
with $\overline{\mathbf{V}}_{i}=(\mathbf{e}_{i}\dot{\varphi}^{i},\mathbf{0})$%
.

The momentum of the main body in mixed representation is $%
%TCIMACRO{\TeXButton{Pi}{\bm{\Pi}}}%
%BeginExpansion
\bm{\Pi}%
%EndExpansion
^{\mathrm{b}}=\mathbf{M}^{\mathrm{b}}\mathbf{V}_{\mathrm{b}}\in {\mathbb{R}}%
^{6}\cong \mathfrak{g}^{\ast }=so^{\ast }\left( 3\right) \times {\mathbb{R}}%
^{3}$, and of the $i$-th rotor $%
%TCIMACRO{\TeXButton{Pi}{\bm{\Pi}}}%
%BeginExpansion
\bm{\Pi}%
%EndExpansion
^{i}=\mathbf{M}^{i}\mathbf{V}_{i}$, with the inertia matrix of the main body
and the $i$-th rotor, respectively,%
\begin{equation}
\mathbf{M}^{\mathrm{b}}=\left( 
\begin{array}{cc}
%TCIMACRO{\TeXButton{Theta}{{\bm{\Theta}}}}%
%BeginExpansion
{\bm{\Theta}}%
%EndExpansion
^{\mathrm{b}} & \mathbf{0} \\ 
\mathbf{0} & m_{\mathrm{b}}\mathbf{I}%
\end{array}%
\right) ,\ \ \mathbf{M}^{i}=\left( 
\begin{array}{cc}
%TCIMACRO{\TeXButton{Theta}{{\bm{\Theta}}}}%
%BeginExpansion
{\bm{\Theta}}%
%EndExpansion
^{i} & \mathbf{0} \\ 
\mathbf{0} & m_{i}\mathbf{I}%
\end{array}%
\right)
\end{equation}%
where $%
%TCIMACRO{\TeXButton{Theta}{{\bm{\Theta}}}}%
%BeginExpansion
{\bm{\Theta}}%
%EndExpansion
^{\mathrm{b}}$ and $%
%TCIMACRO{\TeXButton{Theta}{{\bm{\Theta}}}}%
%BeginExpansion
{\bm{\Theta}}%
%EndExpansion
^{i}$ is the inertia tensor of the main body and the $i$-th rotor w.r.t. the
total COM, and $m_{\mathrm{b}}$ and $m_{i}$ is the mass of the main body and 
$i$-th rotor.

The velocity coordinates are $(\xi ^{\alpha },\dot{r}^{I})=(\omega
^{1},\omega ^{2},\omega ^{3},\dot{p}^{1},\dot{p}^{2},\dot{p}^{3},\dot{\varphi%
}^{1},\dot{\varphi}^{2},\dot{\varphi}^{3})=(V_{\mathrm{b}}^{\alpha },\dot{%
\varphi}^{i})=(\mathbf{V}_{\mathrm{b}},\dot{%
%TCIMACRO{\TeXButton{phi}{\mathbold{\varphi}}}%
%BeginExpansion
\mathbold{\varphi}%
%EndExpansion
})$, with fiber coordinates $\left( \xi ^{\alpha }\right) =\left( V_{\mathrm{%
b}}^{\alpha }\right) =(\omega ^{1},\omega ^{2},\omega ^{3},\dot{p}^{1},\dot{p%
}^{2},\dot{p}^{3}),\alpha =1,\ldots ,6$ and $(\dot{r}^{I})=(\dot{\varphi}%
^{1},\dot{\varphi}^{2},\dot{\varphi}^{3})$, $I=7,8,9$. In the following,
indexes $i,j,k,l=1,2,3$ and the notation $\bar{I}=I-6$ are used. The kinetic
energy of the satellite is 
\begin{eqnarray}
T\left( \mathbf{V}_{\mathrm{b}},\dot{%
%TCIMACRO{\TeXButton{phi}{\mathbold{\varphi}}}%
%BeginExpansion
\mathbold{\varphi}%
%EndExpansion
}\right) &=&\frac{1}{2}\mathbf{V}_{\mathrm{b}}^{T}%
%TCIMACRO{\TeXButton{Pi}{\bm{\Pi}}}%
%BeginExpansion
\bm{\Pi}%
%EndExpansion
^{\mathrm{b}}+\frac{1}{2}\sum_{i=1}^{3}\mathbf{V}_{i}^{T}%
%TCIMACRO{\TeXButton{Pi}{\bm{\Pi}}}%
%BeginExpansion
\bm{\Pi}%
%EndExpansion
^{i}  \notag \\
&=&\frac{1}{2}\Theta _{\alpha \beta }^{\mathrm{b}}\omega ^{\alpha }\omega
^{\beta }+\frac{1}{2}\sum_{i=1}^{3}\left[ \Theta _{jk}^{i}(\omega
^{j}+\delta _{i}^{j}\dot{\varphi}^{i})(\omega ^{k}+\delta _{i}^{k}\dot{%
\varphi}^{i})+\bar{m}\dot{p}^{i}\dot{p}^{i}\right]
\end{eqnarray}%
where $\bar{m}:=m_{\mathrm{b}}+\sum_{i=1}^{3}m_{i}$ is the total mass of the
satellite. The kinetic energy is invariant under left-action of $G$ (due to
the body-fixed angular velocity).

\paragraph{Hamel Equations}

The structure coefficients on the direct product group $SO\left( 3\right)
\times {\mathbb{R}}^{3}$ are $c_{\alpha \lambda }^{\beta }=\varepsilon
_{\alpha \lambda \beta }$, for $\alpha ,\beta ,\lambda =1,2,3$, and $%
c_{\alpha \lambda }^{\beta }=0$ otherwise. The Euler-Poincar\'{e} equations (%
\ref{HEorig1}) are found as%
\begin{align}
\frac{d}{dt}\frac{\partial T}{\partial \omega ^{i}}+\frac{\partial T}{%
\partial \omega ^{j}}c_{ik}^{j}\omega ^{k} &=\bar{\Theta}_{ij}\dot{\omega}%
^{j}+\varepsilon _{ikj}\omega ^{k}\bar{\Theta}_{jl}\omega
^{l}+\sum_{j=1}^{3}(\Theta _{ij}^{j}\ddot{\varphi}^{j}+\varepsilon
_{ikl}\omega ^{k}\Theta _{lj}\dot{\varphi}^{j})  \notag \\
\frac{d}{dt}\frac{\partial T}{\partial \dot{p}^{j}} &=\bar{m}\ddot{p}^{j}
\end{align}%
where $\bar{%
%TCIMACRO{\TeXButton{Theta}{{\bm{\Theta}}}}%
%BeginExpansion
{\bm{\Theta}}%
%EndExpansion
}:=%
%TCIMACRO{\TeXButton{Theta}{{\bm{\Theta}}}}%
%BeginExpansion
{\bm{\Theta}}%
%EndExpansion
^{\mathrm{b}}+\sum_{i=1}^{3}%
%TCIMACRO{\TeXButton{Theta}{{\bm{\Theta}}}}%
%BeginExpansion
{\bm{\Theta}}%
%EndExpansion
^{i}$ is the composite inertia tensor of the satellite including main body
and rotors. They can also be written in vector form, with matrix $\mathbf{ad}%
_{%
%TCIMACRO{\TeXButton{w}{\mathbold{\omega}}}%
%BeginExpansion
\mathbold{\omega}%
%EndExpansion
}=\widetilde{%
%TCIMACRO{\TeXButton{w}{\mathbold{\omega}}}%
%BeginExpansion
\mathbold{\omega}%
%EndExpansion
}$,%
\begin{align}
\frac{d}{dt}\frac{\partial T}{\partial 
%TCIMACRO{\TeXButton{w}{\mathbold{\omega}}}%
%BeginExpansion
\mathbold{\omega}%
%EndExpansion
}-\mathbf{ad}_{%
%TCIMACRO{\TeXButton{w}{\mathbold{\omega}}}%
%BeginExpansion
\mathbold{\omega}%
%EndExpansion
}^{T}\frac{\partial T}{\partial 
%TCIMACRO{\TeXButton{w}{\mathbold{\omega}}}%
%BeginExpansion
\mathbold{\omega}%
%EndExpansion
} &=\bar{%
%TCIMACRO{\TeXButton{Theta}{{\bm{\Theta}}}}%
%BeginExpansion
{\bm{\Theta}}%
%EndExpansion
}\dot{%
%TCIMACRO{\TeXButton{w}{\mathbold{\omega}}}%
%BeginExpansion
\mathbold{\omega}%
%EndExpansion
}+\widetilde{%
%TCIMACRO{\TeXButton{w}{\mathbold{\omega}}}%
%BeginExpansion
\mathbold{\omega}%
%EndExpansion
}\bar{%
%TCIMACRO{\TeXButton{Theta}{{\bm{\Theta}}}}%
%BeginExpansion
{\bm{\Theta}}%
%EndExpansion
}%
%TCIMACRO{\TeXButton{w}{\mathbold{\omega}}}%
%BeginExpansion
\mathbold{\omega}%
%EndExpansion
+\sum_{i=1}^{3}\left( 
%TCIMACRO{\TeXButton{theta}{{\mathbold{\theta}}}}%
%BeginExpansion
{\mathbold{\theta}}%
%EndExpansion
^{i}\ddot{\varphi}^{i}+\widetilde{%
%TCIMACRO{\TeXButton{w}{\mathbold{\omega}}}%
%BeginExpansion
\mathbold{\omega}%
%EndExpansion
}%
%TCIMACRO{\TeXButton{theta}{{\mathbold{\theta}}}}%
%BeginExpansion
{\mathbold{\theta}}%
%EndExpansion
^{i}\dot{\varphi}^{i}\right)  \label{EP1SO3} \\
\frac{d}{dt}\frac{\partial T}{\partial \dot{\mathbf{p}}} &=\bar{m}\ddot{%
\mathbf{p}}  \label{EP2SO3}
\end{align}%
where $%
%TCIMACRO{\TeXButton{theta}{{\mathbold{\theta}}}}%
%BeginExpansion
{\mathbold{\theta}}%
%EndExpansion
^{i}:=%
%TCIMACRO{\TeXButton{Theta}{{\bm{\Theta}}}}%
%BeginExpansion
{\bm{\Theta}}%
%EndExpansion
^{i}\mathbf{e}_{i}$ is the $i$-th column of $%
%TCIMACRO{\TeXButton{Theta}{{\bm{\Theta}}}}%
%BeginExpansion
{\bm{\Theta}}%
%EndExpansion
^{i}$.

The Euler-Lagrange equations (\ref{HEorig2}) are, noting that $T$ does not
depend on $\varphi ^{i}(=r^{\bar{I}})$,%
\begin{equation}
\frac{d}{dt}\frac{\partial T}{\partial \dot{\varphi}^{i}}=\Theta _{ij}^{i}%
\dot{\omega}^{j}+\Theta _{ii}^{i}\ddot{\varphi}^{i}  \label{ELSO3}
\end{equation}%
where the diagonal element $\Theta _{ii}^{i}$ is the moment of inertia about
the axis of rotator $i$.

The above equations are summarized to the set of EOM for the satellite%
\begin{equation}
\left( 
\begin{array}{ccccc}
\bar{%
%TCIMACRO{\TeXButton{Theta}{{\bm{\Theta}}}}%
%BeginExpansion
{\bm{\Theta}}%
%EndExpansion
} & \mathbf{0} & 
%TCIMACRO{\TeXButton{theta}{{\mathbold{\theta}}}}%
%BeginExpansion
{\mathbold{\theta}}%
%EndExpansion
^{1} & 
%TCIMACRO{\TeXButton{theta}{{\mathbold{\theta}}}}%
%BeginExpansion
{\mathbold{\theta}}%
%EndExpansion
^{2} & 
%TCIMACRO{\TeXButton{theta}{{\mathbold{\theta}}}}%
%BeginExpansion
{\mathbold{\theta}}%
%EndExpansion
^{3} \\ 
\mathbf{0} & \bar{m}\mathbf{I} & \mathbf{0} & \mathbf{0} & \mathbf{0} \\ 
%TCIMACRO{\TeXButton{theta}{{\mathbold{\theta}}}}%
%BeginExpansion
{\mathbold{\theta}}%
%EndExpansion
^{1^{T}} & \mathbf{0} & \Theta _{11}^{1} & \mathbf{0} & \mathbf{0} \\ 
%TCIMACRO{\TeXButton{theta}{{\mathbold{\theta}}}}%
%BeginExpansion
{\mathbold{\theta}}%
%EndExpansion
^{2^{T}} & \mathbf{0} & \mathbf{0} & \Theta _{22}^{2} & \mathbf{0} \\ 
%TCIMACRO{\TeXButton{theta}{{\mathbold{\theta}}}}%
%BeginExpansion
{\mathbold{\theta}}%
%EndExpansion
^{3^{T}} & \mathbf{0} & \mathbf{0} & \mathbf{0} & \Theta _{33}^{3}%
\end{array}%
\right) \left( 
\begin{array}{c}
\dot{%
%TCIMACRO{\TeXButton{w}{\mathbold{\omega}}}%
%BeginExpansion
\mathbold{\omega}%
%EndExpansion
} \\ 
\ddot{\mathbf{p}} \\ 
\ddot{\varphi}^{1} \\ 
\ddot{\varphi}^{2} \\ 
\ddot{\varphi}^{3}%
\end{array}%
\right) +\left( 
\begin{array}{c}
\widetilde{%
%TCIMACRO{\TeXButton{w}{\mathbold{\omega}}}%
%BeginExpansion
\mathbold{\omega}%
%EndExpansion
}\bar{%
%TCIMACRO{\TeXButton{Theta}{{\bm{\Theta}}}}%
%BeginExpansion
{\bm{\Theta}}%
%EndExpansion
}%
%TCIMACRO{\TeXButton{w}{\mathbold{\omega}}}%
%BeginExpansion
\mathbold{\omega}%
%EndExpansion
+\sum_{i=1}^{3}\widetilde{%
%TCIMACRO{\TeXButton{w}{\mathbold{\omega}}}%
%BeginExpansion
\mathbold{\omega}%
%EndExpansion
}%
%TCIMACRO{\TeXButton{theta}{{\mathbold{\theta}}}}%
%BeginExpansion
{\mathbold{\theta}}%
%EndExpansion
^{i}\dot{\varphi}^{i} \\ 
\mathbf{0} \\ 
\mathbf{0} \\ 
\mathbf{0} \\ 
\mathbf{0}%
\end{array}%
\right) =\mathbf{0.}  \label{EOM1}
\end{equation}%
The mass matrix has the form (\ref{M}) with non-zero submatrix $\mathbf{K}$.
It is constant due to the assumption of symmetric rotors and axes aligned
with the RFR axes.

The pose of the satellite is obtained by solving the kinematic
reconstruction equations (\ref{KinRec}). To this end, the equations $\dot{%
\mathbf{X}}=\mathbf{dexp}_{-\mathbf{X}}^{-1}\mathbf{V}_{\mathrm{b}}=\left( 
\mathbf{dexp}_{-\mathbf{x}}^{-1}%
%TCIMACRO{\TeXButton{w}{\mathbold{\omega}}}%
%BeginExpansion
\mathbold{\omega}%
%EndExpansion
_{\mathrm{b}},\dot{\mathbf{p}}\right) $ are solved for the coordinate vector 
$\mathbf{X}=\left( \mathbf{x},\mathbf{p}\right) \in {\mathbb{R}}^{6}\cong
so\left( 3\right) \times {\mathbb{R}}^{3}$. Here $\mathbf{dexp}_{\mathbf{x}}$
is the matrix form of the right-trivialized differential of the exp map on $%
SO\left( 3\right) $ \cite{RSPA2021}.

\paragraph{Euler-Lagrange Equations on a Trivial Principle Bundle}

With the partitioning (\ref{M}) of the mass matrix, the connection
coefficients $\mathcal{A}_{I}^{\alpha }$ are defined by%
\begin{equation}
\mathcal{A}=\mathbf{L}^{-1}\mathbf{K}=\left( 
\begin{array}{ccc}
\bar{%
%TCIMACRO{\TeXButton{theta}{{\mathbold{\theta}}}}%
%BeginExpansion
{\mathbold{\theta}}%
%EndExpansion
}^{1} & \bar{%
%TCIMACRO{\TeXButton{theta}{{\mathbold{\theta}}}}%
%BeginExpansion
{\mathbold{\theta}}%
%EndExpansion
}^{2} & \bar{%
%TCIMACRO{\TeXButton{theta}{{\mathbold{\theta}}}}%
%BeginExpansion
{\mathbold{\theta}}%
%EndExpansion
}^{3} \\ 
\mathbf{0} & \mathbf{0} & \mathbf{0}%
\end{array}%
\right) .  \label{A1}
\end{equation}%
The locked velocity (\ref{mechConn}) is $\mathbf{V}_{\mathrm{loc}}=\mathbf{V}%
_{\mathrm{b}}+\mathbf{L}^{-1}\mathbf{K}\dot{%
%TCIMACRO{\TeXButton{phi}{\mathbold{\varphi}}}%
%BeginExpansion
\mathbold{\varphi}%
%EndExpansion
}$, and thus with (\ref{A1}), $\omega _{\mathrm{loc}}^{\alpha }=\omega
^{\alpha }-\mathcal{A}_{I}^{\alpha }\dot{\varphi}^{\bar{I}}$ and $\dot{p}_{%
\mathrm{loc}}^{\alpha }=\dot{p}^{\alpha }$. In terms of the locked velocity,
the kinetic energy is 
\begin{align}
T\left( \mathbf{V}_{\mathrm{loc}},\dot{%
%TCIMACRO{\TeXButton{phi}{\mathbold{\varphi}}}%
%BeginExpansion
\mathbold{\varphi}%
%EndExpansion
}\right) & =\frac{1}{2}\Theta _{\alpha \beta }^{\mathrm{b}}\left( \omega _{%
\mathrm{loc}}^{\alpha }-\mathcal{A}_{\bar{I}}^{\alpha }\dot{\varphi}^{\bar{I}%
}\right) \left( \omega _{\mathrm{loc}}^{\beta }-\mathcal{A}_{\bar{J}}^{\beta
}\dot{\varphi}^{\bar{J}}\right) \\
& +\frac{1}{2}\sum_{i=1}^{3}\Theta _{\alpha \beta }^{i}\left( \omega _{%
\mathrm{loc}}^{\alpha }+\left( \delta _{\bar{I}}^{\alpha }-\mathcal{A}_{\bar{%
I}}^{\alpha }\right) \dot{\varphi}^{\bar{I}}\right) \left( \omega _{\mathrm{%
loc}}^{\beta }+(\delta _{\bar{J}}^{\beta }-\mathcal{A}_{\bar{J}}^{\beta })%
\dot{\varphi}^{\bar{J}}\right) +\frac{\bar{m}}{2}\dot{p}^{\alpha }\dot{p}%
^{\alpha }.  \notag
\end{align}%
A straightforward calculation yields, with block matrices $\mathbf{S,K},%
\mathbf{L}$ deduced from (\ref{EOM1}),%
\begin{equation}
\frac{\partial T}{\partial \mathbf{V}_{\mathrm{loc}}}=\left( 
\begin{array}{cc}
\bar{%
%TCIMACRO{\TeXButton{Theta}{{\bm{\Theta}}}}%
%BeginExpansion
{\bm{\Theta}}%
%EndExpansion
} & \mathbf{0} \\ 
\mathbf{0} & \bar{m}\mathbf{I}%
\end{array}%
\right) \mathbf{V}_{\mathrm{loc}},\ \ \frac{\partial T}{\partial \dot{%
%TCIMACRO{\TeXButton{phi}{\mathbold{\varphi}}}%
%BeginExpansion
\mathbold{\varphi}%
%EndExpansion
}}=\left( \mathbf{S}-\mathcal{A}^{T}\mathbf{L}\mathcal{A}\right) \dot{%
%TCIMACRO{\TeXButton{phi}{\mathbold{\varphi}}}%
%BeginExpansion
\mathbold{\varphi}%
%EndExpansion
}=\left( \mathbf{S}-\mathbf{K}^{T}\mathbf{L}^{-1}\mathbf{K}\right) \dot{%
%TCIMACRO{\TeXButton{phi}{\mathbold{\varphi}}}%
%BeginExpansion
\mathbold{\varphi}%
%EndExpansion
}  \label{Der1}
\end{equation}%
with 
\begin{equation}
\mathbf{S}-\mathbf{K}^{T}\mathbf{L}^{-1}\mathbf{K}=\left( 
\begin{array}{ccc}
\Theta _{11}^{1}-%
%TCIMACRO{\TeXButton{theta}{{\mathbold{\theta}}}}%
%BeginExpansion
{\mathbold{\theta}}%
%EndExpansion
^{1^{T}}\bar{%
%TCIMACRO{\TeXButton{theta}{{\mathbold{\theta}}}}%
%BeginExpansion
{\mathbold{\theta}}%
%EndExpansion
}^{1} & -%
%TCIMACRO{\TeXButton{theta}{{\mathbold{\theta}}}}%
%BeginExpansion
{\mathbold{\theta}}%
%EndExpansion
^{1^{T}}\bar{%
%TCIMACRO{\TeXButton{theta}{{\mathbold{\theta}}}}%
%BeginExpansion
{\mathbold{\theta}}%
%EndExpansion
}^{2} & -%
%TCIMACRO{\TeXButton{theta}{{\mathbold{\theta}}}}%
%BeginExpansion
{\mathbold{\theta}}%
%EndExpansion
^{1^{T}}\bar{%
%TCIMACRO{\TeXButton{theta}{{\mathbold{\theta}}}}%
%BeginExpansion
{\mathbold{\theta}}%
%EndExpansion
}^{3} \\ 
-%
%TCIMACRO{\TeXButton{theta}{{\mathbold{\theta}}}}%
%BeginExpansion
{\mathbold{\theta}}%
%EndExpansion
^{2^{T}}\bar{%
%TCIMACRO{\TeXButton{theta}{{\mathbold{\theta}}}}%
%BeginExpansion
{\mathbold{\theta}}%
%EndExpansion
}^{1} & \Theta _{22}^{2}-%
%TCIMACRO{\TeXButton{theta}{{\mathbold{\theta}}}}%
%BeginExpansion
{\mathbold{\theta}}%
%EndExpansion
^{2^{T}}\bar{%
%TCIMACRO{\TeXButton{theta}{{\mathbold{\theta}}}}%
%BeginExpansion
{\mathbold{\theta}}%
%EndExpansion
}^{2} & -%
%TCIMACRO{\TeXButton{theta}{{\mathbold{\theta}}}}%
%BeginExpansion
{\mathbold{\theta}}%
%EndExpansion
^{2^{T}}\bar{%
%TCIMACRO{\TeXButton{theta}{{\mathbold{\theta}}}}%
%BeginExpansion
{\mathbold{\theta}}%
%EndExpansion
}^{3} \\ 
-%
%TCIMACRO{\TeXButton{theta}{{\mathbold{\theta}}}}%
%BeginExpansion
{\mathbold{\theta}}%
%EndExpansion
^{3^{T}}\bar{%
%TCIMACRO{\TeXButton{theta}{{\mathbold{\theta}}}}%
%BeginExpansion
{\mathbold{\theta}}%
%EndExpansion
}^{1} & -%
%TCIMACRO{\TeXButton{theta}{{\mathbold{\theta}}}}%
%BeginExpansion
{\mathbold{\theta}}%
%EndExpansion
^{3^{T}}\bar{%
%TCIMACRO{\TeXButton{theta}{{\mathbold{\theta}}}}%
%BeginExpansion
{\mathbold{\theta}}%
%EndExpansion
}^{2} & \Theta _{33}^{3}-%
%TCIMACRO{\TeXButton{theta}{{\mathbold{\theta}}}}%
%BeginExpansion
{\mathbold{\theta}}%
%EndExpansion
^{3^{T}}\bar{%
%TCIMACRO{\TeXButton{theta}{{\mathbold{\theta}}}}%
%BeginExpansion
{\mathbold{\theta}}%
%EndExpansion
}^{3}%
\end{array}%
\right)
\end{equation}%
Thus the mass matrix in the EOM in terms of the locked velocity has the
block-diagonal form (\ref{MOmega}). The Hamel coefficients (\ref{gamma41})-(%
\ref{gamma43}) are determined by the non-zero structure coefficients $%
c_{\alpha \lambda }^{\beta }=\varepsilon _{\alpha \lambda \beta }$, for $%
\alpha ,\beta ,\lambda =1,2,3$. 
%TCIMACRO{\TeXButton{red}{\color[rgb]{0,0,0}}}%
%BeginExpansion
\color[rgb]{0,0,0}%
%EndExpansion
Since the connection coefficients are constant, the curvature coefficients (%
\ref{gamma43}) are $\mathcal{B}_{IJ}^{\alpha }=[\mathcal{A}_{I},\mathcal{A}%
_{J}]^{\alpha }$. In vector representation, $\mathcal{A}_{I}$ is the $\bar{I}
$-th column in (\ref{A1}), and $[\mathcal{A}_{I},\mathcal{A}_{J}]=(\bar{%
%TCIMACRO{\TeXButton{theta}{{\mathbold{\theta}}}}%
%BeginExpansion
{\mathbold{\theta}}%
%EndExpansion
}^{\bar{I}}\times \bar{%
%TCIMACRO{\TeXButton{theta}{{\mathbold{\theta}}}}%
%BeginExpansion
{\mathbold{\theta}}%
%EndExpansion
}^{\bar{J}},\mathbf{0})$. It is non-zero due the non-commutativity of vector
fields $\mathcal{A}_{I}$ w.r.t. the Lie bracket on $\mathfrak{g}=so\left(
3\right) \times {\mathbb{R}}^{3}$ (non-parallel rotor axes). The non-zero $%
\mathcal{B}_{IJ}^{\alpha }$ (non-flat connection) implies that the momentum
does not define an integrable relation of rotor and base motion.%
%TCIMACRO{\TeXButton{black}{\color{black}}}%
%BeginExpansion
\color{black}%
%EndExpansion

\subsubsection{Body-fixed Representation of Rigid Body Velocity --- Symmetry
Group $G=SE\left( 3\right) $%
%TCIMACRO{\TeXButton{secSatelliteSE3}{\label{secSatelliteSE3}}}%
%BeginExpansion
\label{secSatelliteSE3}%
%EndExpansion
}

The semi-direct product group $SE\left( 3\right) =SO\left( 3\right) \ltimes {%
\mathbb{R}}^{3}$ describes proper rigid body motions. The configuration of
the main body is again represented as $\left( \mathbf{R},\mathbf{p}\right)
\in SE\left( 3\right) $, but with group multiplication $\left( \mathbf{R}%
_{1},\mathbf{p}_{1}\right) \cdot \left( \mathbf{R}_{2},\mathbf{p}_{2}\right)
=\left( \mathbf{R}_{1}\mathbf{R}_{2},\mathbf{p}_{1}+\mathbf{R}_{1}\mathbf{p}%
_{2}\right) $, which correctly accounts for coupling of rotations and
translations. Thus, the body-fixed RFR can be located arbitrarily. The
velocity (also called twists) of the main body, i.e. of $\mathcal{F}_{%
\mathrm{b}}$, in body-fixed representation \cite{Murray,MUBOScrews1} is
defined via left-trivialization as $\hat{\mathbf{V}}_{\mathrm{b}}=g^{-1}\dot{%
g}=(\mathbf{R}^{-1}\dot{\mathbf{R}},\mathbf{R}^{-1}\dot{\mathbf{p}})=(\hat{%
%TCIMACRO{\TeXButton{w}{\mathbold{\omega}}}%
%BeginExpansion
\mathbold{\omega}%
%EndExpansion
},\mathbf{v})\in \mathfrak{g}=se\left( 3\right) $, and in vector
representation $\mathbf{V}_{\mathrm{b}}=\left( 
%TCIMACRO{\TeXButton{w}{\mathbold{\omega}}}%
%BeginExpansion
\mathbold{\omega}%
%EndExpansion
,\mathbf{v}\right) \in {\mathbb{R}}^{6}$, where now $\mathbf{v}\in {\mathbb{R%
}}^{3}$ is the linear velocity of the main body relative to $\mathcal{F}_{0}$
resolved in $\mathcal{F}_{\mathrm{b}}$. To simplify the derivation, the
velocity of main body and rotors are expressed in the body-fixed RFR at the
main body. The velocity of rotor $i=1,2,3$ is $\mathbf{V}_{i}=\left( 
%TCIMACRO{\TeXButton{w}{\mathbold{\omega}}}%
%BeginExpansion
\mathbold{\omega}%
%EndExpansion
_{i},\mathbf{v}_{i}\right) $. Assuming again that the RFR is aligned with
the rotor axes, it holds true that $\mathbf{V}_{i}=\mathbf{V}_{\mathrm{b}}+%
\overline{\mathbf{V}}_{i}$.

The momentum of the main body in body-fixed representation is $%
%TCIMACRO{\TeXButton{Pi}{\bm{\Pi}}}%
%BeginExpansion
\bm{\Pi}%
%EndExpansion
^{\mathrm{b}}=\mathbf{M}^{\mathrm{b}}\mathbf{V}_{\mathrm{b}}\in {\mathbb{R}}%
^{6}\cong \mathfrak{g}^{\ast }=se^{\ast }\left( 3\right) $, and of the $i$%
-th rotor $%
%TCIMACRO{\TeXButton{Pi}{\bm{\Pi}}}%
%BeginExpansion
\bm{\Pi}%
%EndExpansion
^{i}=\mathbf{M}^{i}\mathbf{V}_{i}$, with the inertia matrix of the main body
and of the $i$-th rotor w.r.t. an arbitrary RFR%
\begin{equation}
\mathbf{M}^{\mathrm{b}}=\left( 
\begin{array}{cc}
%TCIMACRO{\TeXButton{Theta}{\boldmath{\Theta}}}%
%BeginExpansion
\boldmath{\Theta}%
%EndExpansion
^{\mathrm{b}} & m_{\mathrm{b}}\widetilde{\mathbf{d}}_{\mathrm{b}} \\ 
-m_{\mathrm{b}}\widetilde{\mathbf{d}}_{\mathrm{b}} & m_{\mathrm{b}}\mathbf{I}%
\end{array}%
\right) ,\ \ \mathbf{M}^{i}=\left( 
\begin{array}{cc}
%TCIMACRO{\TeXButton{Theta}{\boldmath{\Theta}}}%
%BeginExpansion
\boldmath{\Theta}%
%EndExpansion
^{i} & m_{i}\widetilde{\mathbf{d}}_{i} \\ 
-m_{i}\widetilde{\mathbf{d}}_{i} & m_{i}\mathbf{I}%
\end{array}%
\right)
\end{equation}%
where $%
%TCIMACRO{\TeXButton{Theta}{\boldmath{\Theta}}}%
%BeginExpansion
\boldmath{\Theta}%
%EndExpansion
^{\mathrm{b}}$ and $%
%TCIMACRO{\TeXButton{Theta}{\boldmath{\Theta}}}%
%BeginExpansion
\boldmath{\Theta}%
%EndExpansion
^{i}$ are the inertia tensors of the main body and the $i$-th rotor w.r.t.
the RFR, and $\mathbf{d}_{\mathrm{b}}$, $\mathbf{d}_{i}$ are the position
vectors to the COM w.r.t. the RFR. The total kinetic energy of the satellite
is%
\begin{equation}
T\left( \mathbf{V}_{\mathrm{b}},\dot{%
%TCIMACRO{\TeXButton{phi}{\mathbold{\varphi}}}%
%BeginExpansion
\mathbold{\varphi}%
%EndExpansion
}\right) =\frac{1}{2}\mathbf{V}_{\mathrm{b}}^{T}\mathbf{M}^{\mathrm{b}}%
\mathbf{V}_{\mathrm{b}}+\frac{1}{2}\sum_{i=1}^{3}\mathbf{V}_{i}^{T}\ \mathbf{%
M}^{i}\ \mathbf{V}_{i}=\frac{1}{2}\mathbf{V}_{\mathrm{b}}^{T}\mathbf{M}^{%
\mathrm{b}}\mathbf{V}_{\mathrm{b}}+\frac{1}{2}\sum_{i=1}^{3}\left( \mathbf{V}%
_{\mathrm{b}}+\overline{\mathbf{V}}_{i}\left( \dot{%
%TCIMACRO{\TeXButton{phi}{\mathbold{\varphi}}}%
%BeginExpansion
\mathbold{\varphi}%
%EndExpansion
}\right) \right) ^{T}\mathbf{M}^{i}\ \left( \mathbf{V}_{\mathrm{b}}+%
\overline{\mathbf{V}}_{i}\left( \dot{%
%TCIMACRO{\TeXButton{phi}{\mathbold{\varphi}}}%
%BeginExpansion
\mathbold{\varphi}%
%EndExpansion
}\right) \right) .
\end{equation}

\paragraph{Hamel Equations}

In the following, the matrix form of the equations will be presented, for
simplicity. The structure coefficients on the semi-direct product group $%
SE\left( 3\right) $ give rise to the matrix form of the adjoint operator 
\cite{Selig2005,RSPA2021}%
\begin{equation}
\mathbf{ad}_{\mathbf{V}_{\mathrm{b}}}=\left( 
\begin{array}{cc}
\widetilde{%
%TCIMACRO{\TeXButton{w}{\mathbold{\omega}}}%
%BeginExpansion
\mathbold{\omega}%
%EndExpansion
}_{\mathrm{b}} & \mathbf{0} \\ 
\widetilde{\mathbf{v}}_{\mathrm{b}} & \widetilde{%
%TCIMACRO{\TeXButton{w}{\mathbold{\omega}}}%
%BeginExpansion
\mathbold{\omega}%
%EndExpansion
}_{\mathrm{b}}%
\end{array}%
\right)  \label{adV}
\end{equation}%
so that the Lie bracket is $\mathbf{ad}_{\mathbf{X}}\mathbf{Y}=[\mathbf{X},%
\mathbf{Y}]$. The Euler-Poincar\'{e} equations (\ref{HEorig1}) are%
\begin{equation}
\frac{d}{dt}\frac{\partial T}{\partial \mathbf{V}_{\mathrm{b}}}-\mathbf{ad}_{%
\mathbf{V}_{\mathrm{b}}}^{T}\frac{\partial T}{\partial \mathbf{V}_{\mathrm{b}%
}}=\mathbf{L}\dot{\mathbf{V}}_{\mathrm{b}}+\ \sum_{i=1}^{3}\mathbf{M}^{i}%
\dot{\overline{\mathbf{V}}}_{i}-\mathbf{ad}_{\mathbf{V}_{\mathrm{b}}}^{T}%
%TCIMACRO{\TeXButton{Big}{\Big}}%
%BeginExpansion
\Big%
%EndExpansion
(\mathbf{LV}_{\mathrm{b}}+\sum_{i=1}^{3}\mathbf{M}^{i}\overline{\mathbf{V}}%
_{i}%
%TCIMACRO{\TeXButton{Big}{\Big}}%
%BeginExpansion
\Big%
%EndExpansion
)
\end{equation}%
with locked mass matrix $\mathbf{L}=\mathbf{M}^{\mathrm{b}}+\sum_{i=1}^{3}%
\mathbf{M}^{i}$. Written explicitly yields the instructive form%
\begin{eqnarray}
\bar{%
%TCIMACRO{\TeXButton{Theta}{{\bm{\Theta}}}}%
%BeginExpansion
{\bm{\Theta}}%
%EndExpansion
}\dot{%
%TCIMACRO{\TeXButton{w}{\mathbold{\omega}}}%
%BeginExpansion
\mathbold{\omega}%
%EndExpansion
}+\widetilde{%
%TCIMACRO{\TeXButton{w}{\mathbold{\omega}}}%
%BeginExpansion
\mathbold{\omega}%
%EndExpansion
}\bar{%
%TCIMACRO{\TeXButton{Theta}{{\bm{\Theta}}}}%
%BeginExpansion
{\bm{\Theta}}%
%EndExpansion
}%
%TCIMACRO{\TeXButton{w}{\mathbold{\omega}}}%
%BeginExpansion
\mathbold{\omega}%
%EndExpansion
-\bar{m}(\dot{\widetilde{\mathbf{v}}}+\widetilde{%
%TCIMACRO{\TeXButton{w}{\mathbold{\omega}}}%
%BeginExpansion
\mathbold{\omega}%
%EndExpansion
}\widetilde{\mathbf{v}})\mathbf{d+}\sum_{i=1}^{3}\left( 
%TCIMACRO{\TeXButton{theta}{{\bm{\theta}}}}%
%BeginExpansion
{\bm{\theta}}%
%EndExpansion
^{i}\ddot{\varphi}^{i}+\left( \widetilde{%
%TCIMACRO{\TeXButton{w}{\mathbold{\omega}}}%
%BeginExpansion
\mathbold{\omega}%
%EndExpansion
}%
%TCIMACRO{\TeXButton{theta}{\boldmath{\theta}}}%
%BeginExpansion
\boldmath{\theta}%
%EndExpansion
^{i}-m_{i}\widetilde{\mathbf{v}}\mathbf{a}_{i}\right) \dot{\varphi}%
^{i}\right) &=&\mathbf{0}  \notag \\
\bar{m}\left( \dot{\mathbf{v}}+\widetilde{%
%TCIMACRO{\TeXButton{w}{\mathbold{\omega}}}%
%BeginExpansion
\mathbold{\omega}%
%EndExpansion
}\mathbf{v}+(\dot{\widetilde{%
%TCIMACRO{\TeXButton{w}{\mathbold{\omega}}}%
%BeginExpansion
\mathbold{\omega}%
%EndExpansion
}}+\widetilde{%
%TCIMACRO{\TeXButton{w}{\mathbold{\omega}}}%
%BeginExpansion
\mathbold{\omega}%
%EndExpansion
}\widetilde{%
%TCIMACRO{\TeXButton{w}{\mathbold{\omega}}}%
%BeginExpansion
\mathbold{\omega}%
%EndExpansion
})\mathbf{d}\right) -\sum_{i=1}^{3}m_{i}\left( \mathbf{a}_{i}\ddot{\varphi}%
^{i}+\widetilde{%
%TCIMACRO{\TeXButton{w}{\mathbold{\omega}}}%
%BeginExpansion
\mathbold{\omega}%
%EndExpansion
}\mathbf{a}_{i}\dot{\varphi}^{i}\right) &=&\mathbf{0}
\end{eqnarray}%
with $\bar{%
%TCIMACRO{\TeXButton{Theta}{{\bm{\Theta}}}}%
%BeginExpansion
{\bm{\Theta}}%
%EndExpansion
}:=%
%TCIMACRO{\TeXButton{Theta}{{\bm{\Theta}}}}%
%BeginExpansion
{\bm{\Theta}}%
%EndExpansion
^{\mathrm{b}}+\sum_{i=1}^{3}%
%TCIMACRO{\TeXButton{Theta}{{\bm{\Theta}}}}%
%BeginExpansion
{\bm{\Theta}}%
%EndExpansion
^{i}$ and $%
%TCIMACRO{\TeXButton{theta}{{\mathbold{\theta}}}}%
%BeginExpansion
{\mathbold{\theta}}%
%EndExpansion
^{i}:=%
%TCIMACRO{\TeXButton{Theta}{{\bm{\Theta}}}}%
%BeginExpansion
{\bm{\Theta}}%
%EndExpansion
^{i}\mathbf{e}_{i}$ as above, and $\mathbf{a}_{i}:=\widetilde{\mathbf{d}}_{i}%
\mathbf{e}_{i}$, where $\mathbf{d}:=(m_{\mathrm{b}}\mathbf{d}_{\mathrm{b}%
}+\sum_{i=1}^{3}\mathbf{d}_{i}m_{i})/\bar{m}$ is the position vector of the
total COM measured in the RFR.

The Euler-Lagrange equations (\ref{HEorig2}) are found immediately as%
\begin{equation}
\frac{d}{dt}\frac{\partial T}{\partial \dot{\varphi}^{i}}=\Theta _{ii}^{i}%
\ddot{\varphi}^{i}+\dot{%
%TCIMACRO{\TeXButton{w}{{\mathbold{\omega}}}}%
%BeginExpansion
{\mathbold{\omega}}%
%EndExpansion
}^{T}%
%TCIMACRO{\TeXButton{theta}{{\mathbold{\theta}}}}%
%BeginExpansion
{\mathbold{\theta}}%
%EndExpansion
^{i}-m_{i}\mathbf{a}_{i}^{T}\dot{\mathbf{v}}\text{ \ \ (no summation over }i%
\text{)}.
\end{equation}%
Clearly, if $\mathbf{d}=\mathbf{0}$, i.e. the RFR $\mathcal{F}_{\mathrm{b}}$
is located at the total COM, these equations are equivalent to those in (\ref%
{EP1SO3}-\ref{ELSO3}) when modeling the system on $G=SO\left( 3\right)
\times {\mathbb{R}}^{3}$.

In matrix form, the motion equations are%
\begin{equation}
\left( 
\begin{array}{ccccc}
\bar{%
%TCIMACRO{\TeXButton{Theta}{{\bm{\Theta}}}}%
%BeginExpansion
{\bm{\Theta}}%
%EndExpansion
} & \bar{m}\widetilde{\mathbf{d}} & 
%TCIMACRO{\TeXButton{theta}{{\mathbold{\theta}}}}%
%BeginExpansion
{\mathbold{\theta}}%
%EndExpansion
^{1} & 
%TCIMACRO{\TeXButton{theta}{{\mathbold{\theta}}}}%
%BeginExpansion
{\mathbold{\theta}}%
%EndExpansion
^{2} & 
%TCIMACRO{\TeXButton{theta}{{\mathbold{\theta}}}}%
%BeginExpansion
{\mathbold{\theta}}%
%EndExpansion
^{3} \\ 
-\bar{m}\widetilde{\mathbf{d}} & \bar{m}\mathbf{I} & -m_{1}\mathbf{a}_{1} & 
-m_{2}\mathbf{a}_{2} & -m_{3}\mathbf{a}_{3} \\ 
%TCIMACRO{\TeXButton{theta}{{\mathbold{\theta}}}}%
%BeginExpansion
{\mathbold{\theta}}%
%EndExpansion
^{1^{T}} & -m_{1}\mathbf{a}_{1}^{T} & \Theta _{11}^{1} & \mathbf{0} & 
\mathbf{0} \\ 
%TCIMACRO{\TeXButton{theta}{{\mathbold{\theta}}}}%
%BeginExpansion
{\mathbold{\theta}}%
%EndExpansion
^{2^{T}} & -m_{2}\mathbf{a}_{2}^{T} & \mathbf{0} & \Theta _{22}^{2} & 
\mathbf{0} \\ 
%TCIMACRO{\TeXButton{theta}{{\mathbold{\theta}}}}%
%BeginExpansion
{\mathbold{\theta}}%
%EndExpansion
^{3^{T}} & -m_{3}\mathbf{a}_{3}^{T} & \mathbf{0} & \mathbf{0} & \Theta
_{33}^{3}%
\end{array}%
\right) \left( 
\begin{array}{c}
\dot{%
%TCIMACRO{\TeXButton{w}{\mathbold{\omega}}}%
%BeginExpansion
\mathbold{\omega}%
%EndExpansion
} \\ 
\dot{\mathbf{v}} \\ 
\ddot{\varphi}^{1} \\ 
\ddot{\varphi}^{2} \\ 
\ddot{\varphi}^{3}%
\end{array}%
\right) +\left( 
\begin{array}{c}
\ast \\ 
\ast \ast \\ 
\mathbf{0} \\ 
\mathbf{0} \\ 
\mathbf{0}%
\end{array}%
\right) =\mathbf{0}  \label{EOM2}
\end{equation}%
with $\ast :=\widetilde{%
%TCIMACRO{\TeXButton{w}{\mathbold{\omega}}}%
%BeginExpansion
\mathbold{\omega}%
%EndExpansion
}\bar{%
%TCIMACRO{\TeXButton{Theta}{{\bm{\Theta}}}}%
%BeginExpansion
{\bm{\Theta}}%
%EndExpansion
}%
%TCIMACRO{\TeXButton{w}{\mathbold{\omega}}}%
%BeginExpansion
\mathbold{\omega}%
%EndExpansion
-\bar{m}\widetilde{%
%TCIMACRO{\TeXButton{w}{\mathbold{\omega}}}%
%BeginExpansion
\mathbold{\omega}%
%EndExpansion
}\widetilde{\mathbf{v}}\mathbf{d}+\sum_{i=1}^{3}(\widetilde{%
%TCIMACRO{\TeXButton{w}{\mathbold{\omega}}}%
%BeginExpansion
\mathbold{\omega}%
%EndExpansion
}%
%TCIMACRO{\TeXButton{theta}{{\mathbold{\theta}}}}%
%BeginExpansion
{\mathbold{\theta}}%
%EndExpansion
^{i}-m_{i}\widetilde{\mathbf{v}}\mathbf{a}_{i})\dot{\varphi}^{i}$ and $\ast
\ast :=\bar{m}\left( \widetilde{%
%TCIMACRO{\TeXButton{w}{\mathbold{\omega}}}%
%BeginExpansion
\mathbold{\omega}%
%EndExpansion
}\mathbf{v}+\widetilde{%
%TCIMACRO{\TeXButton{w}{\mathbold{\omega}}}%
%BeginExpansion
\mathbold{\omega}%
%EndExpansion
}\widetilde{%
%TCIMACRO{\TeXButton{w}{\mathbold{\omega}}}%
%BeginExpansion
\mathbold{\omega}%
%EndExpansion
}\mathbf{d}\right) -\sum_{i=1}^{3}m_{i}\widetilde{%
%TCIMACRO{\TeXButton{w}{\mathbold{\omega}}}%
%BeginExpansion
\mathbold{\omega}%
%EndExpansion
}\mathbf{a}_{i}\dot{\varphi}^{i}$.

\paragraph{Eule-Lagrange Equations on a Trivial Principle Bundle}

The mass matrix in (\ref{EOM2}) is block-partitioned, according to (\ref{M}%
), with%
\begin{equation*}
\mathbf{L}=\left( 
\begin{array}{cc}
\bar{%
%TCIMACRO{\TeXButton{Theta}{{\bm{\Theta}}}}%
%BeginExpansion
{\bm{\Theta}}%
%EndExpansion
} & \bar{m}\widetilde{\mathbf{d}} \\ 
-\bar{m}\widetilde{\mathbf{d}} & \bar{m}\mathbf{I}%
\end{array}%
\right) ,\ \mathbf{K}=\left( 
\begin{array}{ccc}
%TCIMACRO{\TeXButton{theta}{{\mathbold{\theta}}}}%
%BeginExpansion
{\mathbold{\theta}}%
%EndExpansion
^{1} & 
%TCIMACRO{\TeXButton{theta}{{\mathbold{\theta}}}}%
%BeginExpansion
{\mathbold{\theta}}%
%EndExpansion
^{2} & 
%TCIMACRO{\TeXButton{theta}{{\mathbold{\theta}}}}%
%BeginExpansion
{\mathbold{\theta}}%
%EndExpansion
^{3} \\ 
-m_{1}\mathbf{a}_{1} & -m_{2}\mathbf{a}_{2} & -m_{3}\mathbf{a}_{3}%
\end{array}%
\right) ,\ \mathbf{S}=\left( 
\begin{array}{ccc}
\Theta _{11}^{1} & \mathbf{0} & \mathbf{0} \\ 
\mathbf{0} & \Theta _{22}^{2} & \mathbf{0} \\ 
\mathbf{0} & \mathbf{0} & \Theta _{33}^{3}%
\end{array}%
\right) .
\end{equation*}%
Therewith, the local connection, defining the locked velocity $\mathbf{V}_{%
\mathrm{loc}}=\mathbf{V}_{\mathrm{b}}+\mathcal{A}\dot{%
%TCIMACRO{\TeXButton{phi}{\mathbold{\varphi}}}%
%BeginExpansion
\mathbold{\varphi}%
%EndExpansion
}$ in (\ref{mechConn}), is%
\begin{equation}
\mathcal{A}=\mathbf{L}^{-1}\mathbf{K}=\left( 
\begin{array}{ccc}
\bar{\mathbf{m}}^{1} & \bar{\mathbf{m}}^{2} & \bar{\mathbf{m}}^{3}%
\end{array}%
\right) 
\end{equation}%
with column vectors $\bar{\mathbf{m}}^{i}:=\mathbf{L}^{-1}\mathbf{M}%
^{i}\left( 
\begin{array}{c}
\mathbf{e}_{i} \\ 
\mathbf{0}%
\end{array}%
\right) $. Explicit expressions for $\mathbf{L},\mathbf{K},\mathcal{A}$ are
given in the supplement \cite{Supplement}. The kinetic energy expressed with
the locked velocity is 
\begin{equation}
T\left( \mathbf{V}_{\mathrm{loc}},\dot{%
%TCIMACRO{\TeXButton{phi}{\mathbold{\varphi}}}%
%BeginExpansion
\mathbold{\varphi}%
%EndExpansion
}\right) =\frac{1}{2}\left( \mathbf{V}_{\mathrm{loc}}-\mathcal{A}\dot{%
%TCIMACRO{\TeXButton{phi}{\mathbold{\varphi}}}%
%BeginExpansion
\mathbold{\varphi}%
%EndExpansion
}\right) ^{T}\mathbf{M}^{\mathrm{b}}\left( \mathbf{V}_{\mathrm{loc}}-%
\mathcal{A}\dot{%
%TCIMACRO{\TeXButton{phi}{\mathbold{\varphi}}}%
%BeginExpansion
\mathbold{\varphi}%
%EndExpansion
}\right) +\frac{1}{2}\sum_{i=1}^{3}\left( \mathbf{V}_{\mathrm{loc}}+%
\overline{\mathbf{V}}_{i}\left( \dot{%
%TCIMACRO{\TeXButton{phi}{\mathbold{\varphi}}}%
%BeginExpansion
\mathbold{\varphi}%
%EndExpansion
}\right) -\mathcal{A}\dot{%
%TCIMACRO{\TeXButton{phi}{\mathbold{\varphi}}}%
%BeginExpansion
\mathbold{\varphi}%
%EndExpansion
}\right) ^{T}\mathbf{M}^{i}\ \left( \mathbf{V}_{\mathrm{loc}}+\overline{%
\mathbf{V}}_{i}\left( \dot{%
%TCIMACRO{\TeXButton{phi}{\mathbold{\varphi}}}%
%BeginExpansion
\mathbold{\varphi}%
%EndExpansion
}\right) -\mathcal{A}\dot{%
%TCIMACRO{\TeXButton{phi}{\mathbold{\varphi}}}%
%BeginExpansion
\mathbold{\varphi}%
%EndExpansion
}\right) .
\end{equation}%
The partial derivatives in (\ref{BHESym1}) and (\ref{BHESym2}) are found
(replacing $\Omega ^{\alpha }$ with $V_{\mathrm{loc}}^{\alpha }$ and $\dot{r}%
^{I}$ with $\dot{\varphi}^{i},i=\bar{I}$) as%
\begin{equation}
\frac{\partial T}{\partial \mathbf{V}_{\mathrm{loc}}}=\mathbf{L}\left( 
\mathbf{V}_{\mathrm{loc}}-\mathcal{A}\dot{%
%TCIMACRO{\TeXButton{phi}{\mathbold{\varphi}}}%
%BeginExpansion
\mathbold{\varphi}%
%EndExpansion
}\right) +\sum_{i=1}^{3}\mathbf{M}^{i}\overline{\mathbf{V}}_{i}=\mathbf{LV}_{%
\mathrm{loc}}  \label{Der2}
\end{equation}%
and $\frac{\partial T}{\partial \dot{%
%TCIMACRO{\TeXButton{phi}{\mathbold{\varphi}}}%
%BeginExpansion
\mathbold{\varphi}%
%EndExpansion
}}$ as in (\ref{Der1}). Consequently, the mass matrix becomes block diagonal
is in (\ref{MOmega}). The Hamel coefficients (\ref{gamma41}-\ref{gamma43})
are determined by the structure constants $c_{\beta \delta }^{\alpha }$ on $%
SE\left( 3\right) $. 
%TCIMACRO{\TeXButton{red}{\color[rgb]{0,0,0}}}%
%BeginExpansion
\color[rgb]{0,0,0}%
%EndExpansion
Again, the curvature $\mathcal{B}_{IJ}^{\alpha }=[\mathcal{A}_{I},\mathcal{A}%
_{J}]^{\alpha }$ does not vanishing because of the non-commutativity of $%
G=SE\left( 3\right) $, where $[\mathcal{A}_{I},\mathcal{A}_{J}]=[\bar{%
\mathbf{m}}^{\bar{I}},\bar{\mathbf{m}}^{\bar{J}}]=\mathbf{ad}_{\bar{\mathbf{m%
}}^{\bar{I}}}\bar{\mathbf{m}}^{\bar{J}}$ is the Lie bracket on $se\left(
3\right) $ in (\ref{adV}), i.e. screw product \cite{Selig2005}. 
%TCIMACRO{\TeXButton{black}{\color{black}}}%
%BeginExpansion
\color{black}%
%EndExpansion
The Euler-Lagrange equations (\ref{BHESym1}),(\ref{BHESym2}) are thus
determined explicitly. Finally, the inverse of the locked mass matrix
attains the closed form 
\begin{equation}
\mathbf{L}^{-1}=\left( 
\begin{array}{cc}
\mathbf{U} & -\mathbf{U}\widetilde{\mathbf{d}} \\ 
\widetilde{\mathbf{d}}\mathbf{U} & \frac{1}{\bar{m}}\mathbf{I}%
%TCIMACRO{\TeXButton{red}{\color[rgb]{0,0,0}}}%
%BeginExpansion
\color[rgb]{0,0,0}%
%EndExpansion
-\widetilde{\mathbf{d}}\mathbf{U}\widetilde{\mathbf{d}}%
%TCIMACRO{\TeXButton{black}{\color{black}}}%
%BeginExpansion
\color{black}%
%EndExpansion
\end{array}%
\right) ,\ \mathrm{with}\ \mathbf{U}=(\bar{%
%TCIMACRO{\TeXButton{Theta}{{\bm{\Theta}}}}%
%BeginExpansion
{\bm{\Theta}}%
%EndExpansion
}+\bar{m}\widetilde{\mathbf{d}}\widetilde{\mathbf{d}})^{-1}.
\end{equation}%
The satellite pose is obtained by solving the local kinematic reconstruction
equations $\dot{\mathbf{X}}=\mathbf{dexp}_{-\mathbf{X}}^{-1}\mathbf{V}_{%
\mathrm{b}}$ for the instantaneous screw coordinate vector $\mathbf{X}%
=\left( \mathbf{x},\mathbf{y}\right) \in {\mathbb{R}}^{6}\cong se\left(
3\right) $, see supplement \cite{Supplement}. The matrix form of the dexp
map on $SE\left( 3\right) $ also possesses a closed form \cite{RSPA2021}.

\section{%
%TCIMACRO{\TeXButton{red}{\color[rgb]{0,0,0}}}%
%BeginExpansion
\color[rgb]{0,0,0}%
%EndExpansion
Floating-Base Mechanical Systems with Symmetry and Conserved Momentum\label%
{secUnConDymConserved}}

\subsection{%
%TCIMACRO{\TeXButton{red}{\color[rgb]{0,0,0}}}%
%BeginExpansion
\color[rgb]{0,0,0}%
%EndExpansion
Hamel Equations, Lagrange--d'Alembert--Poincar\'{e} equations}

Conservation laws can be used to introduce a connection. For floating
systems with $G$-invariant Lagrangian $l(r^{I},\xi ^{\alpha },\dot{r}^{I})$,
the momentum $%
%TCIMACRO{\TeXButton{Pi}{\bm{\Pi}}}%
%BeginExpansion
\bm{\Pi}%
%EndExpansion
\in \mathfrak{g}^{\ast }$, in local bundle coordinates, is%
\begin{equation}
\Pi _{\alpha }=\frac{\partial l}{\partial \xi ^{\alpha }}=L_{\alpha \beta
}\xi ^{\beta }+K_{\alpha J}\dot{r}^{J}.  \label{Momentum}
\end{equation}%
Assuming that the initial momentum is zero, the momentum conservation $\Pi
_{\alpha }=0$ imposes non-holonomic \emph{dynamic constraints} $u^{\alpha
}=0,\alpha =1,\ldots ,\bar{m}$, which are expressed in terms of the
mechanical connection with%
\begin{equation}
u^{\alpha }=\xi ^{\alpha }+\mathcal{A}_{I}^{\alpha }(r^{I})\dot{r}^{I}.
\label{MomentumConst}
\end{equation}%
The Hamel equations are the 
%TCIMACRO{\TeXButton{red}{\color[rgb]{0,0,0}}}%
%BeginExpansion
\color[rgb]{0,0,0}%
%EndExpansion
Lagrange--d'Alembert--Poincar\'{e} equations 
%TCIMACRO{\TeXButton{black}{\color{black}}}%
%BeginExpansion
\color{black}%
%EndExpansion
(\ref{HE-kinBundle}), now with the curvature of the mechanical connection in
(\ref{mechConn}). The connection encodes dynamic constraints due to the
momentum conservation. If the initial momentum is zero, then the locked
velocity is also zero. Comparing $\xi ^{\alpha }=u^{\alpha }-\mathcal{A}%
_{I}^{\alpha }\dot{r}^{I}$, obtained from (\ref{MomentumConst}), with $\xi
^{\alpha }=\Omega ^{\alpha }-\mathcal{A}_{I}^{\alpha }\dot{r}^{I}$, obtained
from (\ref{mechConn}), shows that the equations (\ref{HE-kinBundle}) are
obtained from the equations (\ref{BHESym2}), in terms of the locked
velocity, when $\Omega ^{\alpha }$ is set to zero and the mechanical
connection is used:%
\begin{equation}
\frac{d}{dt}\frac{\partial \ell }{\partial \dot{r}^{I}}-\frac{\partial \ell 
}{\partial r^{I}}+\frac{\partial \ell }{\partial \Omega ^{\beta }}%
%TCIMACRO{\TeXButton{red}{\color[rgb]{0,0,0}}}%
%BeginExpansion
\color[rgb]{0,0,0}%
%EndExpansion
\mathcal{B}_{IJ}^{\beta }%
%TCIMACRO{\TeXButton{black}{\color{black}}}%
%BeginExpansion
\color{black}%
%EndExpansion
\dot{r}^{J}=Q_{I}  \label{HEMomCons}
\end{equation}%
with $\ell (r^{I},\Omega ^{\alpha },\dot{r}^{I}):=l(r^{I},\xi ^{\alpha
}:=\Omega ^{\alpha }-\mathcal{A}_{I}^{\alpha }\dot{r}^{I},\dot{r}^{I})$, 
%TCIMACRO{\TeXButton{red}{\color[rgb]{0,0,0}}}%
%BeginExpansion
\color[rgb]{0,0,0}%
%EndExpansion
and curvature $\mathcal{B}_{IJ}^{\beta }=\gamma _{IJ}^{\beta }$ given by the
Hamel coefficients in (\ref{gamma43})%
%TCIMACRO{\TeXButton{black}{\color{black}}}%
%BeginExpansion
\color{black}%
%EndExpansion
, where $\Omega ^{\alpha }$ is set to zero after taking the derivatives. The
system dynamics is thus described in terms of coordinates $r^{I}$ on the
base manifold (shape space). The motion in $G$ is determined as solution of
the kinematic reconstruction equations (\ref{KinRec}) with $\xi ^{\alpha }=-%
\mathcal{A}_{I}^{\alpha }\left( \mathbf{r}\right) \dot{r}^{I}$ 
%TCIMACRO{\TeXButton{red}{\color[rgb]{0,0,0}}}%
%BeginExpansion
\color[rgb]{0,0,0}%
%EndExpansion
defined by the dynamic constraints $\Pi _{\alpha }=0$%
%TCIMACRO{\TeXButton{black}{\color{black}}}%
%BeginExpansion
\color{black}%
%EndExpansion
.

As example, consider the satellite in Sec. \ref{secUnConBundle}.\ref%
{secSatellite}, with Lagrangian $\ell $ equal to the kinetic energy $T$.
According to (\ref{HEMomCons}), the Hamel equations in terms of the rotor
angles $\varphi ^{i},i=\bar{I}$ are given with 
\begin{equation}
\frac{d}{dt}\frac{\partial T}{\partial \dot{\varphi}^{\bar{I}}}-\frac{%
\partial T}{\partial V_{\mathrm{loc}}^{\beta }}\mathcal{B}_{IJ}^{\beta }\dot{%
\varphi}^{\bar{J}}=Q_{I}
\end{equation}%
with $\frac{\partial T}{\partial \dot{%
%TCIMACRO{\TeXButton{phi}{\mathbold{\varphi}}}%
%BeginExpansion
\mathbold{\varphi}%
%EndExpansion
}^{i}}$, and $\frac{\partial T}{\partial V_{\mathrm{loc}}^{\alpha }}$ in (%
\ref{Der1}) if $G=SO\left( 3\right) \times {\mathbb{R}}^{3}$, and with $%
\frac{\partial T}{\partial V_{\mathrm{loc}}^{\alpha }}$ in (\ref{Der2}) if $%
G=SE\left( 3\right) $. The components of the curvature are the Hamel
coefficients $\gamma _{IJ}^{\alpha }$ in (\ref{gamma43}) given with the
structure constants of the respective symmetry group $G$.

\subsection{%
%TCIMACRO{\TeXButton{red}{\color[rgb]{0,0,0}}}%
%BeginExpansion
\color[rgb]{0,0,0}%
%EndExpansion
Geometric Phase and Pseudo-Holonomic Motion}

The significance of the mechanical connection on the principal bundle is
that it reveals the \emph{geometric phase shift} (holonomy) $dg=-g\mathcal{A}%
d\mathbf{r}$, i.e. the motion in the fiber, as a result of the motion along
a closed curve in shape space (Rem. \ref{remKinControl}), which is
proportional to the curvature (here written for left-trivialization). This
is due to non-integrable condition imposed by the momentum conservation
(while for constrained systems this is due to non-holonomic kinematic
constraints, Rem. \ref{remKinControl}). Whether a closed path in shape space 
${\mathbb{V}}^{\bar{\delta}}$ leads to a closed path in $G$ is a question
arising in context of motion planning of space robots. Although for
non-holonomic systems, this is not possible globally, there may be
trajectories that show such cyclicity. This phenomenon was given the
attribute \emph{pseudo-holonomic}, and necessary conditions were reported in 
\cite{Mukherjee1994,Mukherjee1996} for planar space robots. 
%TCIMACRO{\TeXButton{red}{\color[rgb]{0,0,0}}}%
%BeginExpansion
\color[rgb]{0,0,0}%
%EndExpansion
This aspect was not treated in the literature for general space robots
performing spatial motions. In view of (\ref{int1}), it follows from the
mean value theorem that a necessary condition is the existence of a point $%
\mathbf{r}_{0}\in {\mathbb{V}}^{\bar{\delta}}$ within the area enclosed by
the closed path in shape space such that the curvature of the mechanical
connection vanishes, i.e. $\mathcal{B}\left( \mathbf{r}_{0}\right) =\mathbf{0%
}$. How this can be translated into cyclic 'pseudo-holomic' path planning is
topic of current research. As a simple example, a floating base robot
equipped with an arm comprising two revolute joints with parallel axes is
discusses in the supplement \cite{Supplement}. For this space robot, a
simple cyclic motion of the two joints leads to a pseudo-holonomic behavior
so that the base motion is also cyclic (zero geometric phase). That is,
along this path the base motion is a function of the arm motion, despite the
momentum conservation imposing a non-holonomic constraint.

\subsection{%
%TCIMACRO{\TeXButton{red}{\color[rgb]{0,0,0}}}%
%BeginExpansion
\color[rgb]{0,0,0}%
%EndExpansion
Non-Zero Momentum and the Dynamic Phase}

\begin{figure}[th]
\begin{center}
a)\hspace{-2ex}%
\includegraphics[draft=false,width=6.5cm]{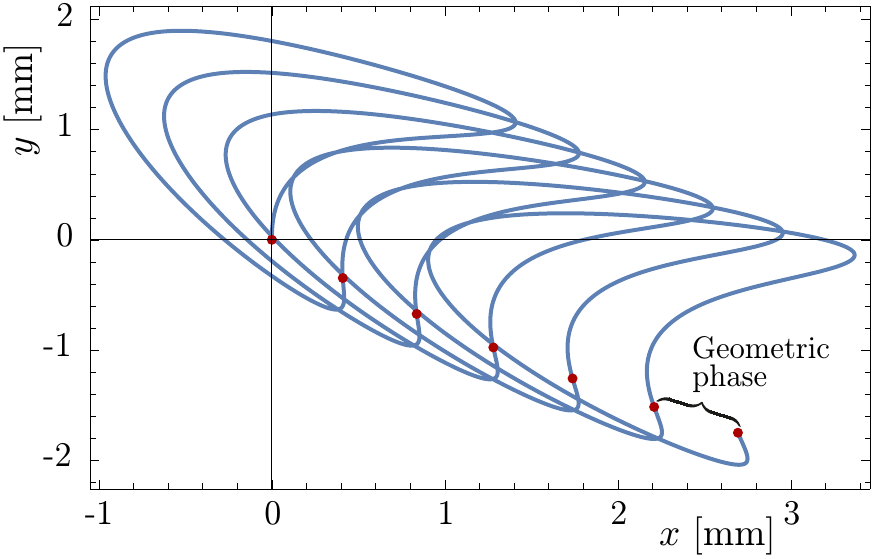}%
~~~~ b)\hspace{-2ex}%
\includegraphics[draft=false,width=6.5cm]{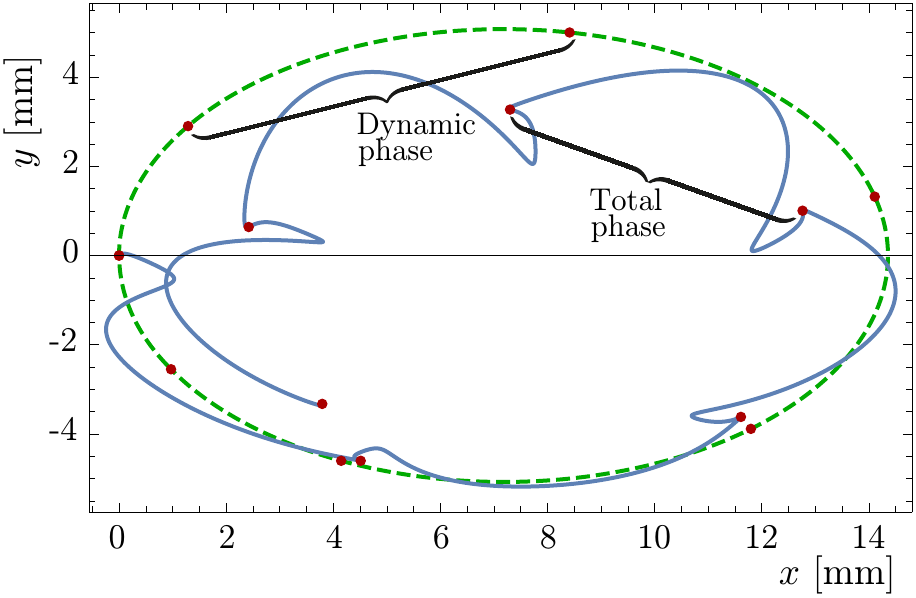}
\end{center}
\caption{a) Translation of base frame $\mathcal{F}_{\mathrm{b}}$ projected
onto the $x$-$y$-plane of $\mathcal{F}_{0}$ (Fig. \protect\ref{figCube}),
when the total momentum is zero, and $\mathcal{F}_{\mathrm{b}}$ and $%
\mathcal{F}_{0}$ initially coincide. Positions after full cycles (with 1 s)
of the rotor motion (geometric phase) are indicated. b) Base motion for
non-zero initial momentum (solid line), with positions after a full cycle of
rotor motions indicated. Shown separately (dashed line) is the motion only
due to the initial momentum, and the corresponding dynamic phase.}
\label{figSatellitePos}
\end{figure}

%TCIMACRO{\TeXButton{red}{\color[rgb]{0,0,0}}}%
%BeginExpansion
\color[rgb]{0,0,0}%
%EndExpansion
Equations (\ref{HEMomCons}) apply also when the total momentum is non-zero.
The centroidal momentum $%
%TCIMACRO{\TeXButton{Pi}{\bm{\Pi}}}%
%BeginExpansion
\bm{\Pi}%
%EndExpansion
_{\mathrm{G}}^{0}=\mathrm{const}$ is the conserved quantity, which is
related to its body-fixed representation by $%
%TCIMACRO{\TeXButton{Pi}{\bm{\Pi}}}%
%BeginExpansion
\bm{\Pi}%
%EndExpansion
^{0}\left( g\right) =\mathbf{Ad}_{g_{\mathrm{bG}}}^{-T}%
%TCIMACRO{\TeXButton{Pi}{\bm{\Pi}}}%
%BeginExpansion
\bm{\Pi}%
%EndExpansion
_{\mathrm{G}}^{0}$ (Sec. \ref{secUnConBundle}\ref{secDecoupling}). The net
change of group variables is determined by the extended reconstruction
equations%
\begin{equation}
dg=g\mathbf{L}^{-1}%
%TCIMACRO{\TeXButton{Pi}{\bm{\Pi}}}%
%BeginExpansion
\bm{\Pi}%
%EndExpansion
^{0}dt-g\widehat{\mathcal{A}d\mathbf{r}}  \label{dg}
\end{equation}%
that replace equations (\ref{KinRec}). Solving the reconstruction equations
for a full cycle along a closed path in ${\mathbb{V}}^{\bar{\delta}}$ yields
the total phase shift as in (\ref{int1}), but now with the additional term $g%
\mathbf{L}^{-1}%
%TCIMACRO{\TeXButton{Pi}{\bm{\Pi}}}%
%BeginExpansion
\bm{\Pi}%
%EndExpansion
^{0}$. The latter delivers the \emph{dynamic phase} which is intrinsically
due to the (initial) momentum, and leads to a symmetry breaking from $G$ to
the symmetry group that preserves the initial momentum. If (\ref{dg}) is
regarded as a control problem, this term is the drift vector field. As an
example, consider the satellite in Sec. \ref{secUnConBundle}\ref%
{secSatellite} with specific parameters. The rotation of the wheels is
prescribed as $\mathbf{r}\left( t\right) =\left( \pi \left( \cos \left( 2\pi
t\right) -1\right) ,\pi \sin \left( 2\pi t\right) ,\pi /2\sin \left( 4\pi
t\right) \right) $, which is periodic with cycle time 1\thinspace s. The
geometric and dynamic parameters, and animations can be found in the
supplementary material \cite{Supplement}. First assume zero total momentum.
The motion of the base (i.e. of base frame $\mathcal{F}_{\mathrm{b}}$
located at geometric center of the base body, as shown in Fig. \ref{figCube}%
) is found from the reconstruction equations (\ref{KinRec2}). Fig. \ref%
{figSatellitePos}a) shows the translation of $\mathcal{F}_{\mathrm{b}}$ in
the $x$-$y$-plane of $\mathcal{F}_{0}$ over 6\thinspace s time duration,
i.e. for six cycles of the rotor motion, starting at the origin. Indicated
is the position after each cycle, which corresponds to the translation
component of the geometric phase. The translation of $\mathcal{F}_{\mathrm{b}%
}$ is caused by the rotation about the total COM, which is not the origin of 
$\mathcal{F}_{\mathrm{b}}$. Fig. \ref{figSatellitePos}b) shows the
translation when the initial momentum is not zero. As an example, the
momentum is set to $%
%TCIMACRO{\TeXButton{Pi}{\bm{\Pi}}}%
%BeginExpansion
\bm{\Pi}%
%EndExpansion
^{0}=\mathbf{K}\dot{\mathbf{r}}\left( 0\right) $, which is the momentum
injected by the rotors when the base is at rest. This resembles the
situation where a satellite is released with spinning fly-wheels. The base
motion is caused by the turning rotors via the non-holonomic kinematics as
well as the dynamics due to the momentum, which determine the total phase.
For completeness, the translation that is generated by the conserved
momentum only when the rotors are rest, is shown as dashed line, which
yields the dynamic phase. 

\section{Remark on Classical Riemannian Geometry Formulations}

It should be recalled that geometric approaches to analytical dynamics of
discrete mechanical systems have a long history. They were originally
developed in the setting of Riemannian geometry for unconstrained systems
with Lagrangian that is quadratic in $\dot{q}^{a}$ \cite%
{Whittaker1988,SyngeSchild1978,Lurie1961}, where the configuration space ${%
\mathbb{V}}^{n}$ is treated as a Riemannian space with metric induced by the
Lagrangian. 
%TCIMACRO{\TeXButton{red}{\color[rgb]{0,0,0}}}%
%BeginExpansion
\color[rgb]{0,0,0}%
%EndExpansion
This was later extended to systems in non-holonomic quasi-velocities and
non-holonomically constrained systems, and Hamel's equations are viewed as
the Lagrange-d'Alembert equations on a configuration manifold whose tangent
space is defined by non-holonomic constraints (which are in classical
literature called non-holonomic tangent bundles). 
%TCIMACRO{\TeXButton{black}{\color{black}}}%
%BeginExpansion
\color{black}%
%EndExpansion
An overview of classical coordinate formulations can be found in \cite%
{Maisser1997,PapastavridisBook2002}, and using modern notations of
differential geometry in \cite{BulloLewis2004}. Only a few publications deal
with rheonomic systems and with systems where the Lagrangian is
non-quadratic in $\dot{q}^{a}$. Such systems are modeled in the $n+1$%
-dimensional event space ${\mathbb{V}}^{n}\times \mathbb{R}_{+}$, which is
treated as a Finsler space. Thus the corresponding model-based control
schemes are developed in event space \cite{FinslerCtrl}. Also in this
classical setting, the connection and its curvature play a key roll. For
example, consider an unconstrained holonomic system with Lagrangian $L\left(
q^{a},\dot{q}^{a}\right) :=\frac{1}{2}g_{ab}\left( q\right) \dot{q}^{a}\dot{q%
}^{b}$ defined by the kinetic energy. The mass matrix defines a Riemannian
metric with coefficients $g_{ab}$ on the configuration space ${\mathbb{V}}%
^{n}$. The system dynamics, on the holonomic tangent bundle $T{\mathbb{V}}%
^{n}$, is governed by the equations%
\begin{equation}
\frac{D}{dt}\dot{q}^{a}=Q^{a}
\end{equation}%
where $D\xi ^{a}=d\xi ^{a}+\Gamma _{bc}^{a}\xi ^{b}dq^{c}$ is the absolute
differential of a contravariant vector field $\xi ^{a}$, and the generalized
forces $Q^{a}=g^{ab}Q_{b}$. The Christoffel symbols of second kind $\Gamma
_{bc}^{a}$ define a natural affine connection, which is metric and symmetric
($\Gamma _{bc}^{a}=\Gamma _{cb}^{a}$), thus the configuration space ${%
\mathbb{V}}^{n}$ is torsion free. While this is a classical result, there is
a beautiful relation for the linearized equations, which is less known.
Denote with $\left( x^{a}\right) \in {\mathbb{R}}^{n}$ small perturbations
superposed to the nominal trajectory $q^{a}$, so that $q^{a}\left( t\right)
+x^{a}\left( t\right) $ is the perturbed trajectory. 
%TCIMACRO{\TeXButton{red}{\color[rgb]{0,0,0}}}%
%BeginExpansion
\color[rgb]{0,0,0}%
%EndExpansion
The linearized equations along the nominal trajectory $q^{a}$ are, in
covariant form, with the \emph{Riemann-Christoffel curvature }tensor $%
R_{cbda}$,%
\begin{equation}
g_{ab}\frac{D^{2}x^{b}}{dt^{2}}+\left( R_{cbda}\dot{q}^{c}\dot{q}^{d}-\nabla
_{b}Q_{a}\right) x^{b}=\Phi _{a}
\end{equation}%
the covariant derivative $\nabla _{b}Q_{a}=\frac{\partial Q_{a}}{\partial
q^{b}}-\Gamma _{ba}^{c}Q_{c}$, and small generalized forces $\Phi _{a}$ dual
to $x^{a}$. 
%TCIMACRO{\TeXButton{black}{\color{black}}}%
%BeginExpansion
\color{black}%
%EndExpansion
The curvature is hence a \emph{measure of stability} of the perturbed
dynamics.

\section{Conclusion}

%TCIMACRO{\TeXButton{red}{\color[rgb]{0,0,0}}}%
%BeginExpansion
\color[rgb]{0,0,0}%
%EndExpansion
The classical Hamel formulation is a generally applicable approach in
analytical mechanics for describing the dynamics of finite-dimensional
systems in terms of local coordinates, which can be extended to continua 
\cite{ShiZenkovDmitryBloch2020}. Frequently, the coordinate form of
equations that can be derived coordinate-free in the framework of geometric
mechanics are referred to as Hamel equations. The link between these
conceptually very different approaches has not been sufficiently addressed,
however. This link was established in this paper, where the key is to
identify the Hamel coefficients as the coefficients appearing in the
coordinate form of the reduced Euler-Lagrange equations, respectively the
Lagrange-Poincar\'{e} equations. Of particular significance are the local
curvature coefficients that are central in many aspects of control and
computational treatment of mechanical systems whose configuration space is a
non-linear manifold or a Lie-group. Therewith, a clear connection between
the equations governing the dynamics on a principle bundle, defined by a
connection originating from certain symmetries, and the original Hamel
formulation is established. In this context the choice of bundle coordinates
is crucial. As such the locked velocity, and the related concept of average
velocity, were discussed. The locked velocity leads to inertial decoupling,
which is important for control and computational investigations. This should
motivate further research into their use for deriving formulations with
improved efficiency. A problem that is increasingly receiving attention is
that the average velocity cannot be used for kinematic reconstruction \cite%
{MuleroMartinez2008,Garofalo2015,Saccon2017,Nava2018,Wenqian2021}. This
could, for instance, be addressed by means of holonomy minimizing gauge
transformations, i.e. introducing a frame that is not body-fixed nor aligned
with the inertia frame. As a geometric aspect of the motion of non-holomic
systems, it was discussed how the geometric phase leads to attitude change
of floating systems for instance, and that there may be pseudo-holonomic
motions. Since this is naturally covered by the geometric approach, it shall
motivate treating Hamel's formalism in a geometric setting. It remains to be
explored how Hamel's formulation can be extended to the general case when
constraints and Lagrangian possess (possibly complementary) symmetries, as
treated in \cite{MarsdenScheurle1993,BlochKrishnaprasadMarsdenMurray1996},
where a non-holonomic connection is introduced generalizing the kinematic
and mechanical connection. As a side-contribution, some differences and
inconsistencies of the definition of local curvature found in the literature
were identified, which is crucial when applying equations (\ref{LP1},\ref%
{LP2}). 
%TCIMACRO{\TeXButton{black}{\color{black}}}%
%BeginExpansion
\color{black}%
%EndExpansion

\section*{%
%TCIMACRO{\TeXButton{red}{\color[rgb]{0,0,0}}}%
%BeginExpansion
\color[rgb]{0,0,0}%
%EndExpansion
A. List of Symbols%
%TCIMACRO{\TeXButton{secSymbols}{\label{secSymbols}}}%
%BeginExpansion
\label{secSymbols}%
%EndExpansion
}

%TCIMACRO{\TeXButton{red}{\color[rgb]{0,0,0}}}%
%BeginExpansion
\color[rgb]{0,0,0}%
%EndExpansion
\begin{tabular}{lrl}
$n$ & -- & number of (generalized) coordinates \\ 
$\bar{m}$ & -- & i) number of Pfaffian constraints, ii) dimension of the
symmetry group $G$ \\ 
$\bar{\delta}=n-\bar{m}$ & -- & differential (instantaneous) DOF defined by
the $\bar{m}$ Pfaffian constraints \\ 
$a,b,c,\ldots $ & -- & indices running over all coordinates: $a=1,\ldots ,n$
\\ 
$I,K,L,\ldots $ & -- & indices $I=\bar{m}+1,\ldots ,n$ of i) independent
velocity, ii) shape coordinates \\ 
&  & (i.e. coordinates of the base manifold of the principle bundle) of a \\ 
&  & constrained or unconstrained system \\ 
$\alpha ,\beta ,\gamma ,\ldots $ & -- & indices $\alpha =1,\ldots ,\bar{m}$
of i) constraint equations, ii) dependent velocity \\ 
&  & coordinates, iii) canonical coordinates on the symmetry group $G$ \\ 
$i,j,k,l,\ldots $ & -- & indices $i,j,k,l=1,2,3$ of Cartesian vectors, e.g. $%
\mathbf{x}=(x^{i})\in {\mathbb{R}}^{3}$ \\ 
$q^{a},\ \mathbf{q}=\left( q^{a}\right) $ & -- & local coordinates,
generalized coordinates of unconstrained system \\ 
$r^{I},\ \mathbf{r}=(r^{I})$ & -- & independent (local) coordinates, $I=\bar{%
m}+1,\ldots ,n$ \\ 
$s^{\alpha },\ \mathbf{s}=\left( s^{\alpha }\right) $ & -- & dependent
(local) coordinates, $\alpha =1,\ldots ,\bar{m}$ \\ 
$\gamma _{bc}^{a},\gamma _{IJ}^{\alpha },\gamma _{\beta J}^{\alpha }$ & -- & 
Hamel coefficients \\ 
$u^{\alpha }$ & -- & i) quasi-velocities, ii) bundle coordinates, iii)
Pfaffian constraints \\ 
$\Omega ^{\alpha }$ & -- & local coordinates of the locked velocity \\ 
$\mathcal{A}_{I}^{\alpha },\mathcal{B}_{IJ}^{\alpha }$ & -- & coefficients
of the local connection and of the local curvature \\ 
$\varepsilon _{ijk}$ & -- & Levi-Civita symbol \\ 
$\delta _{ij}$ & -- & Kronecker delta symbol \\ 
$\widetilde{\mathbf{x}}$ & -- & skew symmetric matrix $\widetilde{\mathbf{x}}%
=(\varepsilon _{ikj}x^{k})$ associated to vector $\mathbf{x}=(x^{k})\in {%
\mathbb{R}}^{3}$ \\ 
$\mathbf{x}\times \mathbf{y}$ & -- & cross product of $\mathbf{x},\mathbf{y}%
\in {\mathbb{R}}^{3}$, can be written as $\widetilde{\mathbf{x}}\mathbf{y}$
\\ 
$\hat{\bm{\eta}}\in \mathfrak{g}$ & -- & element of Lie algebra $\mathfrak{g}
$, corresponding to vector $\bm{\eta}\in {\mathbb{R}}^{n}\cong \mathfrak{g}$
\\ 
$c_{\beta \lambda }^{\alpha }$ & -- & structure coefficients of the Lie
algebra $\mathfrak{g}$ \\ 
$\left[ X,Y\right] $ & -- & Lie bracket $\left[ X,Y\right] =c_{\beta \lambda
}^{\alpha }X^{\beta }Y^{\lambda }$ of $X,Y\in \mathfrak{g}$ \\ 
$\mathcal{F}_{0},\mathcal{F}_{\mathrm{b}}$ & -- & Inertial frame (IFR) $%
\mathcal{F}_{0}$, body-fixed frame $\mathcal{F}_{\mathrm{b}}$ \\ 
$SE\left( 3\right) $ & -- & special Euclidean group (rigid body motion
group) $SE\left( 3\right) =SO\left( 3\right) \ltimes {\mathbb{R}}^{3}$ \\ 
$SO\left( 3\right) $ & -- & special orthogonal group (rotation group)%
\end{tabular}

Ricci's summation convention: e.g. $B_{I}^{a}u^{I}=\sum_{I}B_{I}^{a}u^{I}$, $%
c_{\beta \lambda }^{\alpha }X^{\beta }Y^{\lambda }=\sum_{\beta
}\sum_{\lambda }c_{\beta \lambda }^{\alpha }X^{\beta }Y^{\lambda }$%
%TCIMACRO{\TeXButton{black}{\color{black}}}%
%BeginExpansion
\color{black}%
%EndExpansion

\bibliographystyle{IEEEtran}
\bibliography{HamelCoeff_ResearchGate}

\end{document}